\documentclass[11pt]{article}

\usepackage{graphicx}
\usepackage{epstopdf}

\usepackage{amsmath, amsfonts, amssymb, amsthm}
\usepackage{enumitem}
\usepackage[T1]{fontenc}    
\usepackage{url}            
\usepackage{booktabs}       
\usepackage{nicefrac}       
\usepackage{microtype}      
\usepackage{subfigure}

\usepackage{comment}

\usepackage{mathrsfs}

\usepackage{fullpage}
\usepackage[usenames,dvipsnames,svgnames,table]{xcolor}
\definecolor{darkgreen}{rgb}{0.0,0,0.9}
\usepackage[pagebackref,colorlinks=true,pdfpagemode=none,citecolor=OliveGreen,linkcolor=BrickRed,urlcolor=BrickRed,pdfstartview=FitH]{hyperref}
\usepackage{color}

\theoremstyle{plain}
\newtheorem{theorem}{Theorem}

\newtheorem{proposition}{Proposition}
\newtheorem{corollary}{Corollary}
\newtheorem*{theorem*}{Theorem}
\newtheorem*{lemma*}{Lemma}
\newtheorem*{proposition*}{Proposition}
\newtheorem*{corollary*}{Corollary}
\newtheorem{remark}{Remark}
\newtheorem*{remark*}{Remark}

\newtheorem*{note*}{Note}
\theoremstyle{definition}

\newtheorem*{definition*}{Definition}
\newtheorem{assumption}{Assumption}

\usepackage{xcolor}
\usepackage{dsfont}
\newcommand{\odima}[1]{{\color{violet}#1}}

\newcommand{\wt}{\widetilde}
\newcommand{\wh}{\widehat}
\newcommand{\R}{\mathds{R}}
\newcommand{\E}{\mathds{E}}
\newcommand{\Prob}{\mathds{P}}

\newcommand{\rank}{\textup{rank}}
\newcommand{\ind}{\mathds{1}}

\newcommand{\cN}{\mathcal{N}}
\newcommand{\cZ}{\mathcal{Z}}

\newcommand{\cY}{\mathcal{Y}}

\newcommand{\cH}{\mathcal{H}}
\newcommand{\G}{\boldsymbol{G}}
\renewcommand{\H}{\boldsymbol{H}}

\newcommand{\M}{\boldsymbol{J}}
\newcommand{\Mopt}{\M^{\boldsymbol\star}}
\newcommand{\bSigma}{\boldsymbol{\Sigma}}
\newcommand{\bPi}{\boldsymbol{\Pi}}
\newcommand{\Id}{\boldsymbol{I}}
\newcommand{\weakto}{\rightsquigarrow}
\newcommand{\Tr}{\textup{Tr}}
\newcommand{\proofstep}[1]{$\boldsymbol{{#1}^o}$}
\newcommand{\Var}{\textup{Var}}
\newcommand{\Bias}{\beta}
\renewcommand{\tilde}{\wt}
\renewcommand{\hat}{\wh}
\DeclareMathOperator*{\argmin}{argmin}

\newcommand{\tht}{\theta}
\newcommand{\rmin}{\texttt{r}}
\newcommand{\cP}{\mathcal{P}}
\newcommand{\cE}{\mathcal{E}}
\newcommand{\cC}{\mathcal{C}}
\newcommand{\piC}{\pi_{\cC}}
\newcommand{\bOmega}{\boldsymbol{\Omega}}

\renewcommand{\le}{\leq}
\renewcommand{\ge}{\geq}
\renewcommand{\leq}{\leqslant}
\renewcommand{\geq}{\geqslant}
\newcommand{\pinv}{\dagger}
\newcommand{\bA}{\boldsymbol{A}}
\newcommand{\bB}{\boldsymbol{B}}

\newcommand{\ue}{\boldsymbol{e}}

\newcommand{\learner}{\textit{Learner}}
\newcommand{\auditor}{\textit{Auditor}}

\newcommand{\Cov}{\textup{Cov}}

\newcommand{\cT}{\textsf{T}}
\newcommand{\sfC}{\textsf{\textup{C}}}

\newcommand{\cL}{\mathcal{L}}

\newcommand{\col}{\textup{col}}
\newcommand{\LR}{\textsf{LR}}
\newcommand{\MLR}{\textsf{MLR}}

\newcommand{\eff}{\textsf{eff}}
\newcommand{\reff}{r_{\eff}}

\newcommand{\dem}{\textup{dem}}
\newcommand{\rep}{\textup{rep}}

\newcommand{\Lin}{\textsf{\textup{Lin}}}
\newcommand{\Grad}{\textsf{\textup{Grad}}}
\if false
\author{%
  Author\thanks{Use footnote for providing further information
    about author (webpage, alternative address)---\emph{not} for acknowledging
    funding agencies.} \\
  Department of ....\\
  University of Southern California\\
  Los Angeles, LA 90007 \\
  \texttt{author1@usc.edu} \\
   \And
   Coauthor \\
   Department of ....\\
   University of Southern California \\
  Los Angeles, LA 90007 \\
   \texttt{author2@usc.edu} \\
   \And
   Adel Javanmard \\
   Marshall School of Business\\
   University of Southern California \\
  Los Angeles, LA 90007 \\
   \texttt{author3@usc.edu} \\
  \And
   Adel J \\
   Department of ....\\
   University of Southern California \\
  Los Angeles, LA 90007 \\
   \texttt{author4@usc.edu} \\
}
\fi

\title{
Near-Optimal Procedures for Model Discrimination  \\ with Non-Disclosure Properties
}
			
\newcommand*\samethanks[1][\value{footnote}]{\footnotemark[#1]}

\author{
Dmitrii M.~Ostrovskii\thanks{Equal contribution of the first two authors.}\,
\thanks{University of Southern California, Viterbi School of Engineering, 3650 McClintock Ave, Los Angeles, CA 90089, USA. Email: \texttt{dostrovs@usc.edu}.}
\hspace{0.45cm}
Mohamed Ndaoud\samethanks[1]\,
\thanks{ESSEC Business School, 3 avenue Bernard Hirsch, 95021 Cergy-Pontoise Cedex, France. Email: \texttt{ndaoud@essec.edu}.}
\hspace{0.45cm}
Adel Javanmard\thanks{University of Southern California, Marshall School of Business,
3670 Trousdale Pkwy, Los Angeles, CA 90089, USA. Email: \texttt{ajavanma@usc.edu}.}
\hspace{0.45cm}
Meisam Razaviyayn\thanks{University of Southern California, Viterbi School of Engineering, 3650 McClintock Ave, Los Angeles, CA 90089, USA. Email: \texttt{razaviya@usc.edu}.}
}
\date{}

\begin{document}

\maketitle

\begin{abstract}
Let $\theta_0,\theta_1 \in \R^d$ be the population risk minimizers associated to some loss $\ell: \R^d \times \mathcal{Z} \to \R$ and two distributions $\Prob_0,\Prob_1$ on $\mathcal{Z}$. The models $\theta_0$ and $\theta_1$ are unknown and the distributions $\Prob_0, \Prob_1$ can be accessed by drawing i.i.d samples from them. Our work is motivated by the following model discrimination question:
\begin{quote}
    {\em Given i.i.d.~samples from~$\Prob_0$ and\,~$\Prob_1$, what sample sizes are sufficient and necessary to distinguish between the two hypotheses~$\theta^* = \theta_0$ and~$\theta^* = \theta_1$ for given~$\theta^* \in \{\theta_0, \theta_1\}$?}
\end{quote}
Making the first steps towards answering this question in full generality, we first consider the case of a well-specified linear model with squared loss.
Here we provide matching upper and lower bounds on the sample complexity as given by~$\min\{1/\Delta^2, \sqrt{r}/\Delta\}$ up to a constant factor; here~$\Delta$ is a measure of separation between~$\Prob_0$ and~$\Prob_1$ and~$r$ is the rank of the design covariance matrix.
This bound is dimension-independent, and rank-independent for large enough separation. 
We then extend this result in two directions: (i) for the general parametric setup in asymptotic regime; (ii) for generalized linear models in the small-sample regime~$n \le r$ and under weak moment assumptions.
In both cases we derive sample complexity bounds of a similar form while allowing for model misspecification.
In fact, our testing procedures only access $\theta^*$ through a certain functional of empirical risk. In addition, the number of observations that allows to reach statistical confidence in our tests does not allow to ``resolve'' the two models -- that is, recover $\theta_0,\theta_1$ up to $O(\Delta)$ prediction accuracy. These two properties allow to use our framework in applied tasks where one would like to {\em identify} a prediction model, which can be proprietary, while guaranteeing that the model cannot be actually {\em inferred} by the agent performing identification.
\end{abstract}

\section{Introduction}
\label{sec:intro}
 
Statistical theory teaches us that testing is generally easier than estimation or prediction. 
This reasoning is exploited, for example, when proving minimax lower bounds in parametric estimation. Indeed, the minimax risk in such problems can often be lower-bounded in terms of the testing error in a hypothesis testing problem associated to the estimated signal or parameter~\cite{tsybakov_mono,nemirovski2000topics}.
More interestingly, in some situations one observes {\em quantifiable gaps} between the sample complexity of an estimation or prediction problem and that of the testing problem associated to the initial problem in the natural sense.\footnote{For this expository discussion, we define the sample complexity of a (binary) testing problem as the size of an i.i.d.~sample for which there exists a test with testing errors of both types at most~$0.05$.} 
For example, such a situation arises in certain detection-type problems where the goal is to detect the presence of a ground-truth signal in the background noise.
Signal detection can be relatively easy because it does not necessarily entail localizing the signal or estimating its direction; thus one does not have to deal with the complexity of the associated hypotheses spaces.
Along these lines, statistical and computational gaps between signal detection and estimation have been observed, for example, in sparse linear regression~\cite{donoho2016high,ingster2010detection} and in the spiked covariance model~\cite{berthet2013optimal}.

In this paper, we identify and study a class of testing problems
in which the ``testing-easier-than-estimation'' phenomenon can be {\em used for practical purposes}, by relying on the following observation:
\begin{align}
\label{key-idea}
\parbox{\dimexpr\linewidth-5em}{
{\em One might {\bf test hypotheses} about a parametric predictive model by observing its output, with provable guarantees of not being able to actually {\bf recover the model},
whenever the number of observations necessary for recovery is larger than the one sufficient for testing.}
}
\tag{\textit{\textbf{$\boldsymbol{\star}$}}}
\end{align}
Before we make this high-level observation more concrete, let us define the class of testing problems that we study in this paper.
This class is natural and interesting in its own right, not merely in the context of~\eqref{key-idea}.
Much to our surprise, it seems to not have been studied in the literature so far.

\paragraph{Problem formulation.}
Let~$z$ be a random observation in some space~$\cZ$, and let~$\ell(\cdot,z): \R^d \to \R$ be a random loss function associated to~$z$.
Our task is to distinguish between the two hypotheses
\begin{equation}
\label{eq:hypotheses}
\cH_0: \{\theta^* = \theta_0 \}, \quad 
\cH_1: \{\theta^* = \theta_1 \}
\end{equation}
where~$\theta_0,\theta_1$ minimize the population risks~$L_0(\cdot), L_1(\cdot)$ associated with two distributions~$\Prob_0, \Prob_1$ of~$z$: 
\begin{equation}
\label{eq:minimizers}
\theta_k = \argmin_{\theta \in \R^d} \left\{ L_k(\tht) := \E_{z \sim \Prob_k} [\ell(\tht,z)] \right\} \quad (k \in \{0,1\}).
\end{equation}
In~\eqref{eq:hypotheses},~$\theta^* \in \R^d$ is known; thus, hypotheses~\eqref{eq:hypotheses} can be understood as those about the unknown~$\theta_0, \theta_1$.
The loss~$\ell(\cdot,z): \R^d \to \R$ is also known, and is assumed strictly convex in~$\theta$, which guarantees the uniqueness of~$\theta_0,\theta_1$.
To construct a test, the statistician can generate i.i.d.~samples from~$\Prob_0$ and~$\Prob_1$:
\begin{equation}
\label{eq:samples}
Z^{(0)} := \left(z^{(0)}_{1}, ..., z^{(0)}_{\vphantom{0}n_0}\right) \sim \Prob_0^{\otimes n_0}, \quad 
Z^{(1)} := \left(z^{(1)}_{1}, ..., z^{(1)}_{\vphantom{1}n_1}\right) \sim \Prob_1^{\otimes n_1}.
\end{equation}
The performance of a test~$\wh T: (\theta^*, Z^{(0)}, Z^{(1)}, \ell) \mapsto \{0,1\}$ can be measured by the largest of two error probabilities~$\Prob_{\cH_0}[\wh T=1], \Prob_{\cH_1}[\wh T=0]$. 
We characterize the sample complexity for this problem in terms of~$n_0, n_1$ that are necessary and sufficient to guarantee fixed error probabilities of both types.\\

At first glance, the problem formulation in~\eqref{eq:hypotheses}--\eqref{eq:samples} might seem artificial.
However, the testing procedures we develop to address it turn out to be useful in a number of practical tasks 
where, conforming to~\eqref{key-idea}, one aims at {\em identifying} a statistical prediction model while at the same time providing a guarantee that the model cannot be {\em recovered} by a testing agent with due accuracy.\footnote{We shall specify the precise meaning of ``recovery with due accuracy'' in the next section.} 
Such guarantees  are central for our testing procedures when looking at them from the viewpoint of applications.  
We shall refer to such guarantees as the {\em non-disclosure property} of a testing procedure.    


In the next section, we give a high-level explanation of the mechanism behind our proposed testing procedures and discuss in more detail how the idea~\eqref{key-idea} allows to guarantee the non-disclosure property for them.
For convenience of the subsequent discussion, we shall abstract out the specific details arising in testing procedures due to slight variations in the setup.
Instead, we shall look at all these procedures through the lens of a {\em unified testing protocol} that summarizes their common properties and allows to take advantage of the non-disclosure properties in several applications.


%
%
\subsection{Testing protocol  and non-disclosure property}
\label{sec:protocol}

\paragraph{Access through empirical prediction scores.}
The mechanism behind the non-disclosure property of our testing procedures relies on the following fact:
\begin{quote}
{\em While in the formulation~\eqref{eq:hypotheses}--\eqref{eq:samples} the data~$(\theta^*,Z^{(0)}, Z^{(1)},\ell)$ is formally assumed to be known to the testing agent, 
our {\bf actual testing procedures} do not require direct access to these data.
Rather, the data is accessed through \mbox{a pair of {\bf empirical prediction scores}.}}
\end{quote}
Formally, the empirical prediction scores of~$\theta^*$ for two samples~$Z^{(0)},Z^{(1)}$ are the local Newton decrements (see~\cite{nesterov2013introductory}) of the empirical risks
\[
\wh L_0(\theta) := \frac{1}{n_0}\sum_{i=1}^{n_0} \ell(\theta,z_i^{(0)}), 
\quad 
\wh L_1(\theta) := \frac{1}{n_1}\sum_{i=1}^{n_1} \ell(\theta,z_i^{(1)})
\] 
at~$\theta^*$, that is
\begin{equation}
\label{eq:newton-decrements}
\| \wh \H_0(\theta^*)^{{\dagger}/{2}} \nabla \wh L_0(\theta^*) \|^2, \quad
\| \wh \H_1(\theta^*)^{{\dagger}/{2}} \nabla \wh L_1(\theta^*) \|^2
\end{equation}
where~$\wh \H_k(\theta)^{{\dagger}}$ ($k \in \{0,1\}$) is the generalized inverse of the empirical risk Hessian~$\wh\H_k(\theta) := \nabla^2 \wh L_k(\theta)$, and~$\wh \H_k(\theta)^{{\dagger/2}}$ is the (unique) positive-semidefinite square root of~$\wh \H_k(\theta)^{{\dagger}}$.
The Newton decrements~\eqref{eq:newton-decrements} quantify how well~$\theta^*$ fits each of the two samples in terms of the local affine-invariant gradient norm.
Our testing procedures can be run without directly communicating the data~$(\theta^*, Z^{(0)}, Z^{(1)}, \ell)$ to the testing agent, but only granting access to the Newton decrements.\footnote{This is a simplification for the initial discussion. In fact, in every setting of interest we start with a basic test that requires additional parameters of~$\Prob_0,\Prob_1$ such as noise variances or effective ranks; we then construct an {\em adaptive} test that uses estimates obtained from~$Z^{(0)}, Z^{(1)}$ or hold-out samples. 
We shall address these considerations later on.}
It is especially important here that~$\theta^*$ does not have to be directly communicated: as we shall see in Section~\ref{sec:applications}, in many applications~$\theta^*$ specifies the ``proprietary'' prediction model which must remain undisclosed to the testing agent.

\paragraph{Testing protocol.}
The testing protocol which shall be presented next is motivated by the following considerations.
In applications to be discussed in Section~\ref{sec:applications},~$\theta^*$ specifies a prediction model trained on a large dataset. 
Direct access to such model, which allows to compute the empirical prediction scores~\eqref{eq:newton-decrements}, is a prerogative of its owner, possibly a private corporation. 
We shall refer to this entity as \learner{}. 
On the other hand, testing the hypotheses in~\eqref{eq:hypotheses} is a task of a {different entity}, called \auditor{}, who cannot directly access~$\theta^*$. We discuss several applications in details in Section~\ref{sec:applications}, but it would be useful to briefly consider an example here to motivate the testing protocol. A natural right for users of a platform is the ``right to be forgotten''. The users may request to remove all or part of their data from the platform database e.g. due to privacy concerns. Nonetheless, the platform may have incentive to ignore such requests e.g., for a better predictive model or to avoid retraining its model. In this example, $\theta_0$ and $\theta_1$ correspond to the models trained with and without the datapoints to be removed, and an auditor would like to test which model is actually used by the platform.
Motivated by such applications, we require two natural constraints on the testing protocol:
\begin{enumerate}[label={(\em \alph*)}]
\item \label{req:i}
We would like to protect \learner{}'s model~$\theta^*$ from recovery by~\auditor{} (in the exact sense to be defined in the next section) while still allowing \auditor{} to run the test and distinguish between the two hypotheses in~\eqref{eq:hypotheses}, that is, conclude which of the two samples~$Z^{(0)}, Z^{(1)}$ was generated by the distribution corresponding to~$\theta^*$.
\item \label{req:ii}
In addition, we would like to protect the complementary model\vspace{-0.2cm}
\[
\bar\theta := \theta_0 + \theta_1 - \theta^*
\vspace{-0.2cm}
\]
from recovery by either~\learner{} or \auditor{}. 
\end{enumerate}
One way of respecting both these two requirements is by granting~\auditor{} access only to the bare minimum of information that suffices to run the test, namely the two prediction scores~\eqref{eq:newton-decrements} (this guarantees~\ref{req:i}) and by choosing the sample sizes~$n_0, n_1$ small enough for not allowing to recover~$\theta_0, \theta_1$ from the data~$(Z^{(0)}, Z^{(1)}, \ell)$ (this allows to guarantee~\ref{req:ii} by exploiting~\eqref{key-idea}).
This approach can be implemented in a formal {\bf testing protocol} that involves three parties: \learner{} (owner of~$\theta^*$), \auditor{} (or testing agent), and {\em sampling oracles} for~$\Prob_0,\Prob_1$ (see Figure~\ref{fig:protocol}):
\begin{figure}[t!]
    \centering
    \includegraphics[width=0.7\textwidth]{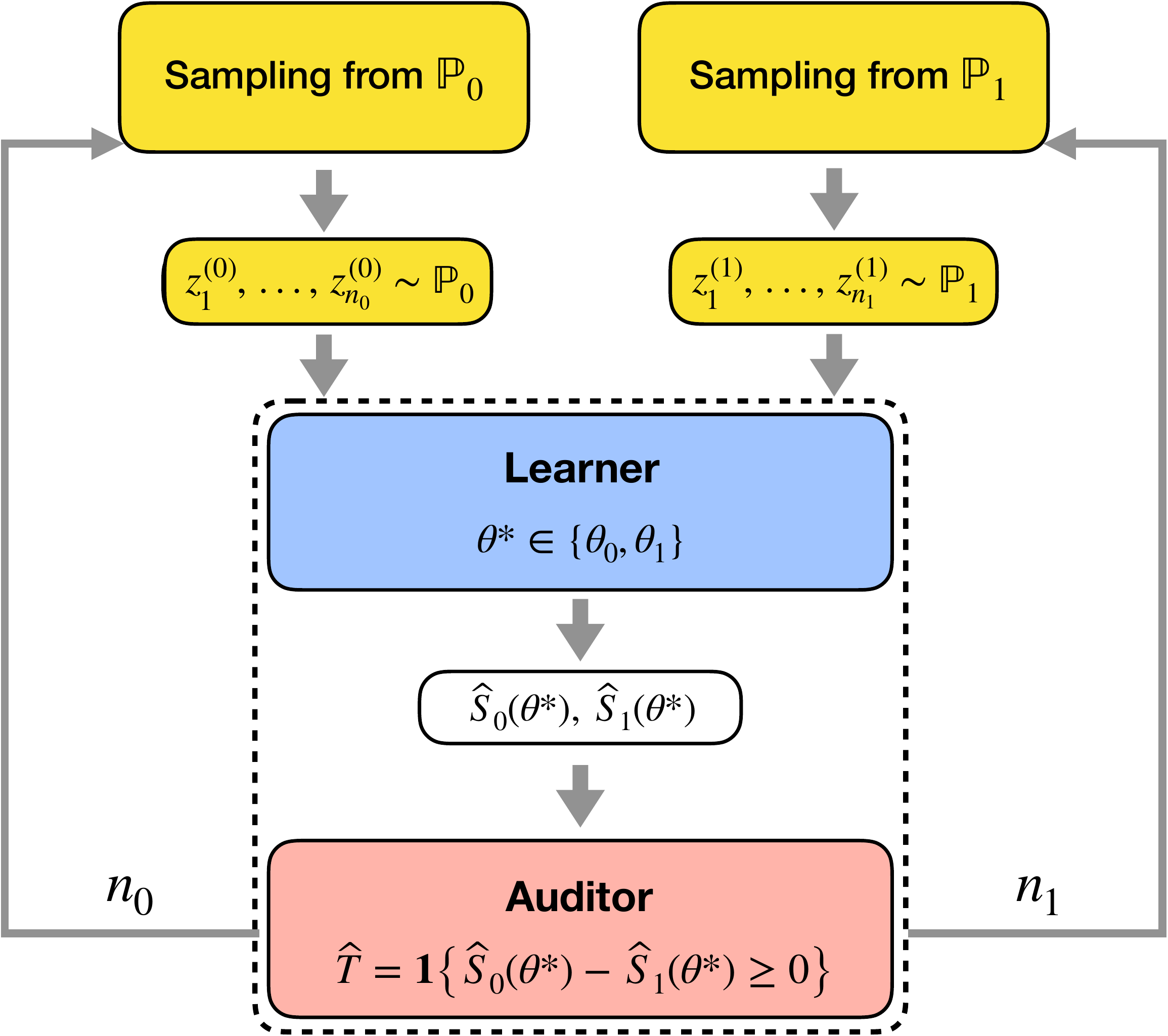}
    \caption{Testing protocol used in the applications of our framework. 
    \learner{} and \auditor{} may be unified into a single entity, which we depict via the dashed box.}
    \label{fig:protocol}
\end{figure}
\begin{enumerate}
\item 
\learner{} is granted direct access to the ground-truth model~$\theta^*$ and local access to the loss~$\ell(\theta,z)$ at~$\theta = \theta^*$ which allows to compute the gradient and Hessian~$\nabla\ell(\theta^*,z), \nabla^2 \ell(\theta^*,z)$ over~$\theta$ for any~$z \in \cZ$. In practice, \learner{} may also have access to the ground-truth distribution, i.e., the one among~$\Prob_0, \Prob_1$ corresponding to~$\theta^*$; however, this is never used in the protocol. 
\learner{} is also granted restricted access to both distributions~$\Prob_0, \Prob_1$ through {\em sampling oracles}. Specifically, the sampling oracle corresponding to~$\Prob_0$ (resp.,~$\Prob_1$) can form an i.i.d.~sample $Z^{(0)}$ (resp.,~$Z^{(1)}$) of the requested size~$n_0$ (resp.,~$n_1$). 
\item  
\auditor{} comes up with two sample sizes~$n_0, n_1$ and communicates them to {\em sampling oracles} who form the two samples~$Z^{(0)},Z^{(1)}$ of the requested sizes. These samples are sent to~\learner{}. 
\item
Using the access modalities described above, \learner{} computes a pair of statistics~$\wh S_0(\theta^*), \wh S_1(\theta^*)$ based on the empirical prediction scores  for the two samples (cf.~\eqref{eq:newton-decrements}). While it is~\auditor{} who gets to decide on the sample sizes, it is only~\learner{}, not~\auditor{}, who accesses the samples.
\item 
\learner{} then communicates the two statistics~$\wh S_0(\theta^*), \wh S_1(\theta^*)$ back to~\auditor{}, who simply compares them and accordingly chooses the hypothesis~$\cH_k \in \{\cH_0, \cH_1\}$ ``with the best fit.'' 
\end{enumerate}

\paragraph{Non-disclosure mechanism.}
Above we described our preferred variant of the testing protocol with minimal access modality for \auditor{}.
In practice, \auditor{} might enjoy a more favorable access modality, say, have access to~$(Z^{(0)}, Z^{(1)}, \ell)$ in addition to~$\wh S_0(\theta^*), \wh S_1(\theta^*)$, or even have access to the full data~$(\theta^*, Z^{(0)}, Z^{(1)}, \ell)$ if \auditor{} and \learner{} are united into a single entity (say, as two divisions in a single technology company). 
%
%
Observe, however, that we can still fulfill both requirements~\ref{req:i}, \ref{req:ii} in the first of these cases, and requirement~\ref{req:ii} in the second case, if we guarantee that the data
\[
\left( Z^{(0)},Z^{(1)}, \ell, \wh S_0(\theta^*), \wh S_1(\theta^*) \right)
\]
allows to recover neither~$\theta_0$ nor~$\theta_1$. 
Furthermore, intuitively these data contain only a little more statistical information about~$\theta_0$ and~$\theta_1$ than the data
\[
\left( Z^{(0)},Z^{(1)}, \ell \right)
\]
---after all,~$\wh S_0(\theta^*), \wh S_1(\theta^*)$ are just two random measurements, whereas each sample~$Z^{(k)}$ contains~$n_k$ measurements. (We refer to Appendix~\ref{app:fisher} for a formal discussion through the Fisher information.)
Let us now complete the description of the non-disclosure mechanism by specifying the precise meaning of ``impossibility to recover~$\theta_0,\theta_1$ from the data~$\left( Z^{(0)},Z^{(1)}, \ell \right)$.'' 
%


To this end, assume that~$\theta_0, \theta_1$ are separated ``in terms of excess risk,'' namely
\begin{align}
\Delta : = \min \left\{ L_{0}(\theta_1) - L_0(\theta_0), L_{1}(\theta_0) - L_1(\theta_1) \right\} > 0.
\end{align}
In Section~\ref{sec:contrib}, we derive sample complexity estimates for confident testing in~\eqref{eq:hypotheses}--\eqref{eq:samples}, namely the required sample size such that both types of testing error are bounded by arbitrarily small, fixed value say 0.05. As we will see the derived sample complexity bounds for confident testing are controlled by  the separation~$\Delta$, and are orderwise (in terms of $\Delta$) smaller than the complexity of estimating~$\theta_0,\theta_1$ from the data~$(Z^{(0)}, Z^{(1)}, \ell)$ up to excess risk~$O(\Delta)$, i.e., the required sample size so that
\begin{equation}
\label{eq:estimation-goal}
\min \left\{\E_{\Prob_0}[L_0(\wh \theta_0)] - L_0(\theta_0),\; \E_{\Prob_1}[L_1(\wh \theta_1)] - L_1(\theta_1)\right\} = O(\Delta).
\end{equation}
Therefore, one can choose the sample sizes~$n_0, n_1$ that suffice for distinguishing between the two hypotheses in~\eqref{eq:hypotheses}, but do {\em not} allow to recover~$\theta_0$ or~$\theta_1$ with the prediction accuracy better than~$\Delta$, i.e., better than the one we get by simply using~$\theta_1$ instead of~$\theta_0$ and vice versa.

As a result, given full data $(\theta^*, Z^{(0)}, Z^{(1)}, \ell)$, \learner{} cannot recover the complementary model~$\bar\theta$
with ``prediction accuracy'' better than~$\Delta$.
Similarly, even when granted access to~$(Z^{(0)}, Z^{(1)},\ell)$,~\auditor{} cannot recover~\learner{}'s model~$\theta^*$ with prediction accuracy better than~$\Delta$. 

\begin{remark}
\label{rem:multiple-testing}
Our testing framework and the protocol in Figure~\ref{fig:protocol} can be easily extended to the multiple-testing setup, where one has to distinguish between~$m$ hypotheses corresponding to~$\Prob_{0},..., \Prob_{m-1}$. To this end, it suffices to increase the number of sampling oracles and replace the decision rule by~$\wh T \in \argmin_{\; k \in \{0,..., m-1\}} \wh S_k(\theta^*)$. Our theory can be easily extended to this multiple testing setup.
\end{remark}

Next we discuss some practical applications of our testing framework. 

\subsection{Applications}
\label{sec:applications}
\paragraph{Applications related to ``the right to be forgotten.''}

Ubiquitous collection and storage of large volumes of user data by internet corporations poses societal risks such as potential for unfair/discriminatory outcomes and privacy violations.
To confront such tendencies, governments and consumer protection agencies have established guidelines and regulations such as~\cite{GDPR1,GDPR2,CCPA}. 
In a nutshell, these documents require from organizations to provide users with ``the right to be forgotten,'' i.e., removal of users' data upon request. 
In the simplest scenario, this means removal of the part of a big dataset pertaining to a single user or a group of users. 
However, the corporate entity might be incentivized to ignore or violate such a removal request for various reasons, such as potential degradation of the model performance after the data removal procedure or a high computational cost of retraining the model.\footnote{While techniques for removing the impact of data points from a model without retraining it have been proposed~\cite{guo2019certified}, their analysis is only tight in simple cases, and lacks generalization guarantees.}
Thus, users might be interested in verifying whether the corporate entity complied with their deletion requests. Ideally, such verification is to be done by a third party, without compromising users' data nor revealing the trained model to the third party. 

The testing protocol in Figure~\ref{fig:protocol} can be used in such an application.
Here we view~$\Prob_0$ as the original dataset and~$\Prob_1$ the dataset after receiving (potentially many) data removal requests. 
(Both datasets are assumed to be large, hence are modeled by population distributions.) 
The respective population risk minimizers~$\theta_0, \theta_1$ correspond to the model trained with or without the datapoints to be removed. 
Under~$\cH_1$,~\learner{} indeed removed the data, whereas~$\cH_0$ corresponds to the request violation; in both cases,~$\theta^*$ is the {\em actual} model trained by~\learner{}.~\auditor{}, which can be a third party or an internal agent in a technological platform, can verify if the data removal indeed took place by asking~\learner{} to evaluate the empirical prediction scores~$\theta^*$ on the ``snapshots''~$Z^{(0)}, Z^{(1)}$ of the full datasets, with and without the datapoints requested for deletion.
Such verification does not require to provide~\auditor{} with direct access to neither the trained model~$\theta^*$ nor the full datasets~$\Prob_0,\Prob_1$, and protects~$\theta^*$ from being learned by~\auditor{} as discussed in Section~\ref{sec:protocol}. 

Note that in this application, it would be problematic to assume that both sampling oracles are readily available.
Instead, it is reasonable to assume that \auditor{} has access to the dataset of deletion requests (let us call it~$Q$), and can also ask \learner{} for a subsample of their data.
Now, we claim that this allows~\auditor{} to emulate both sample oracles (for~$\Prob_0$ and~$\Prob_1$) assuming that \learner{} behaves rationally, i.e., is interested in demonstrating {\em compliance with deletion requests}. 
\begin{itemize}
\item 
Indeed, first observe that~$\Prob_0 = (1-\delta)\Prob_1 + \delta Q$, where~$\delta$ is the share of deletion requests; thus, assuming the knowledge of~$\delta$, a sampling oracle for~$\Prob_0$ can be emulated given those for~$\Prob_1$ and~$Q$. Observe also that~\auditor{} has full access to~$Q$, and thus can also sample from it.
\item 
On the other hand, a rational \learner{}, in the sense defined above, will sample from~$\Prob_1$ under either hypothesis: clearly, this is the case under~$\cH_1$, but this is also the case under~$\cH_0$, since in this case~\learner{} would like to feign the execution of deletion requests.\footnote{In fact, if any datapoint from~$Q$ is included into~\learner{}'s subsample,~\auditor{} can immediately detect \learner{}'s non-compliance
by finding a duplicate of this datapoint in~$Q$. This is possible thanks to \auditor's full access to~$Q$.} This grants \auditor{} a sampling oracle for~$\Prob_1$, and thus for~$\Prob_0$. 
\end{itemize}

\paragraph{Verification of fair representation of minority groups in training dataset.}
In the context of fair machine learning~\cite{datta2015automated, sweeney2013discrimination, bolukbasi2016man, machinebias2016}, a common task is to ensure that a prediction model does not lead to discriminatory outcomes against minority groups. 
Such discriminatory outcomes may be caused by an imbalanced representation of sub-populations in the dataset used in training; 
for example, the minority sub-populations will suffer higher test error due to their small share in the training dataset.
In this context, our testing protocol can be used to verify whether the training data has equal representation of sub-populations.
To this end, a data platform plays the role of \learner{} and has access to the trained model~$\theta^*$ while \auditor{} functions as a verification agent;
the sampling oracles can be implemented by collecting (small) datasets among separate subpopulations and appropriately mixing them.
More precisely, assume there are two different subpopulations that we shall call Democrats and Republicans. 
The null-hypothesis corresponds to a balanced dataset with an equal representation of both of them, and the alternative is an imbalanced dataset skewed by a margin of~$p-1/2 > 0$ towards democrats.
Denoting~$\Prob_{\dem}$ and~$\Prob_{\rep}$ the two populations, we have that~$\Prob_0 = \tfrac{1}{2} \Prob_{\dem} + \tfrac{1}{2} \Prob_{\rep}$ and~$\Prob_1 = p \Prob_{\dem} + (1-p) \Prob_{\rep}$. The sampling oracles can then be implemented by sampling from~$\Prob_{\dem}$ and~$\Prob_{\rep}$, tossing a fair coin for~$\Prob_0$, and a biased coin for~$\Prob_1$.
We note that the resulting protocol can be viewed as the fundamental building block of the \textit{multiparty computation framework} proposed in \cite{jagielski2018differentially, veale2017fairer, kilbertus2018blind}.

\paragraph{Applications related to user identification.}
Finally, consider the task of identifying a user whose data has been used to train a predictive model. 
Here, there are~$m$ datasets~$\Prob_1, ..., \Prob_m$ pertaining to~$m$ users of a data platform. 
In a similar vein to the previous example, \auditor{} must infer which of those~$m$ datasets has been used by the~\learner{} to train the model. 
This can be done by measuring the prediction scores of the trained model on the snapshots of~$m$ datasets.
Note that this application requires a multiple-hypothesis extension of the framework as discussed in Remark~\ref{rem:multiple-testing}.\\


We shall now overview our technical contributions pertaining to the problem defined in~\eqref{eq:hypotheses}--\eqref{eq:samples}.

\subsection{Overview of technical contributions, literature and paper organization}
\label{sec:contrib}

We make important first steps towards solving the general case of the problem summarized in~\eqref{eq:hypotheses}--\eqref{eq:samples}.

In Section~\ref{sec:gaussian}, we focus on the case of well-specified linear regression with random design, where we have an almost complete picture. 
In other words, here we assume that~$z = (x,y)$ with~$x \in \R^d$ and~$y \in \R$, and the distributions~$\Prob_0, \Prob_1$ are given by
\begin{equation}
\label{eq:intro-gaussian}
\Prob_k: x \sim \cN(0,\bSigma_k), \; y|x \sim \cN(x^\top \theta_{k}^{\vphantom 2},\sigma_k^2) \;\; \text{for} \;\; k \in \{0,1\}.
\end{equation}
We first assume the noise variances to be known, which is equivalent to~$\sigma_0^2 = \sigma_1^2 = 1$ by rescaling; however, other parameters $\theta_0,\theta_1,\bSigma_0,\bSigma_1$ are assumed unknown. 
Writing the two samples in~\eqref{eq:samples} in a concise form as~$(X^{(0)},Y^{(0)})$ and~$(X^{(1)}, Y^{(1)})$ with~$X^{(k)} \in \R^{n_k \times d}$,~$Y^{(k)} \in \R^{n_k}$ for~$k \in \{0,1\}$,
our approach is based on the high-level idea of comparing the squared norms of residuals at~$\theta^*$ for the two samples, that is~$\|Y^{(0)} - X^{(0)} \theta^*\|^2$ and~$\|Y^{(1)} - X^{(1)} \theta^*\|^2$. We further adjust this idea as follows:
\begin{itemize}
\item
Instead of using the residuals directly, we first project them onto the respective ``signal spaces'', i.e., the column spaces of~$X^{(0)}$ and~$X^{(1)}$, thus passing to~$\|\bPi_{X^{(0)}}[Y^{(0)} - X^{(0)} \theta^*]\|^2$ and~$\|\bPi_{X^{(1)}}[Y^{(1)} - X^{(1)} \theta^*]\|^2$. (Here~$\bPi_A := A(A^{\top}A)^{\dagger}A^{\top}$ is the projector onto~$\col(A)$, the column space of~$A$.)
\item
Instead of comparing the squared norms of projected residuals, we compare the {\em deviations} of the squared norms from their conditional expectations under~$\Prob_0,\Prob_1$. 
We exploit the fact that these expectations are observed and given by the ranks of empirical covariance matrices~$\wh \bSigma_0, \wh \bSigma_1$.
\end{itemize}
Assuming~$n_0 = n_1$ for the sake of simplicity, these adjustments result in a test with sample complexity
\begin{equation}
\label{eq:intro-complexity}
O\bigg(\min\bigg\{\frac{1}{\Delta^2}, \frac{\sqrt{\max\{r_0, r_1\}}}{\Delta} \bigg\} \bigg),
\end{equation}
where~$r_k = \rank(\bSigma_k)$ and~$\Delta = \min \{\Delta_0, \Delta_1\}$ is the smallest of the two (squared) prediction distances
\[
\Delta_k := \|\bSigma_k^{1/2}(\theta_1 - \theta_0)\|^2 \quad (k \in \{0,1\}).
\] 
Note that the sample complexity bound~\eqref{eq:intro-complexity} depends on~$\max\{r_0, r_1\}$ rather than the ambient dimension~$d$, and even this dependency vanishes when~$\Delta \ge 1/\sqrt{\max\{r_0, r_1\}}$.
We later show that bound~\eqref{eq:intro-complexity} is near-optimal in the minimax sense, up to the replacement of~$\max\{r_0,r_1\}$ by~$\min\{r_0, r_1\}$. 

Finally, we consider the case of unknown and unequal noise variances. 
The challenge here is that the conditional expectations of~$\|\bPi_{X^{(0)}}[Y^{(0)} - X^{(0)} \theta^*]\|^2$ and~$\|\bPi_{X^{(1)}}[Y^{(1)} - X^{(1)} \theta^*]\|^2$ are not given explicitly anymore, and have to be inferred from observations. 
Fortunately, this turns out to be possible: we construct an {\em adaptive} test that performs essentially the same as if~$\sigma_0^2, \sigma_1^2$ were known.\\

In Section~\ref{sec:asymp}, we revisit the general~$M$-estimation setup.
Here we adjust the test introduced in the linear model setup by replacing the quantities~~$\|\bPi_{X^{(0)}}[Y^{(0)} - X^{(0)} \theta^*]\|^2$ and~$\|\bPi_{X^{(1)}}[Y^{(1)} - X^{(1)} \theta^*]\|^2$ with {\em empirical prediction scores}, i.e., the Newton decrements~\eqref{eq:newton-decrements} of empirical risks~$\wh L_0(\theta^*), \wh L_1(\theta^*)$ over the two samples. 
As the case of a linear model falls under this more general scenario, the {lower} bound on the sample complexity established in that case, and nearly matching~\eqref{eq:intro-complexity}, still applies.
However, deriving a matching {upper} bound proves to be challenging in the general case.
Hence, we focus on the asymptotic regime of large sample sizes and small prediction distances, namely~$n_0, n_1 \to \infty$ and~$\Delta_{0},\Delta_1 \to 0$ with 
finite~$n_k\Delta_k \to \lambda_k$.
In this regime, a local Gaussian approximation applies~\cite{Lecam1986}, and we are in a position to use the central limit theorem in our analysis.
In terms of the derived dependency of the testing risk from~$\lambda$, the obtained results nearly match the lower bound, but also take into account the effect of model misspecification.
We also demonstrate that the second-order information plays a critical role: the risk for a natural test in which one measures the gradients of empirical risks without using their Hessians can be {\em arbitrarily larger} than for our test.\\

In Section~\ref{sec:glm}, we extend our theory in another direction, now focusing on the {\em small-sample regime} 
\begin{equation}
\label{eq:small-sample-regime-intro}
n_k \le r_k \quad (k \in \{0,1\}).
\end{equation}
While estimation is clearly impossible in this regime unless~$\Delta \gg 1$, the testing problem is still well-posed. 
Indeed, according to~\eqref{eq:intro-complexity}, in the simple case~$n_0 = n_1$ the sample complexity is at most~$O(\max\{r_0,r_1\})$ unless~$\Delta$ is very small (more precisely, unless~$\Delta \ll 1/\sqrt{\max\{r_0,r_1\}}$). 
Thus, for linear models condition~\eqref{eq:small-sample-regime-intro} does not restrict us from constructing tests with near-optimal complexity bounds, provided that~$r_0,r_1$ are of the same order and~$\Delta$ is not too small. 
A similar conclusion for generalized linear models (GLMs) follows directly from the upper bound on sample complexity derived in~Section~\ref{sec:glm}. 
Namely, for GLMs under~\eqref{eq:small-sample-regime-intro} we again use a Newton-decrement based statistic, and prove an analogue of bound~\eqref{eq:intro-complexity} under weak moment assumptions and moderate level of misspecification. 
As a byproduct, we extend the results of Sec.~\ref{sec:gaussian} to linear models with heavy-tailed noise distribution.
From a technical viewpoint, small sample size leads to the diagonalization of the Newton decrements, with two advantages as a result: (i) simpler statistical analysis; (ii) faster computation of the test.\\ 

In Section~\ref{sec:simul}, we present the results of numerical experiments. 
They show that our approach has better empirical performance compared to those that do not use the second-order information and are based only on the values or gradients of the empirical risk.

\paragraph{Related work.}
Our testing problem, as formulated in~\eqref{eq:hypotheses}--\eqref{eq:samples}, appears to be new. 
The closest line of research concerns hypothesis testing in linear regression. Related literature is very rich and it would be hard to thoroughly overview it here.
Still, one important direction is signal testing in linear regression~\cite{baraud2002non,comminges2013minimax,mukherjee2020minimax}. Variants of this problem were considered when the signal is sparse and the noise either known~\cite{carpentier2019minimax} or unknown~\cite{carpentier2020estimation}. Other works test, in the same flavor, sparsity of the signal~\cite{carpentier2019optimal} or some component of the signal \cite{bradic2018testability}.  
We encourage an interested reader to study these recent works and references therein for further details. 

Let us, however, emphasize the differences between the classical setup of parametric testing and our setup in~\eqref{eq:hypotheses}-\eqref{eq:samples}.
In the classical parametric testing setup, one is given a sample from a parametric distribution with unknown mean $\theta^*$ and we are asked to choose $\theta^*$ from $\{\theta_0,\theta_1\}$. 
In contrast, in~\eqref{eq:hypotheses}-\eqref{eq:samples} one is given a parameter $\theta^*$ and two samples with respective means $\{\theta_0,\theta_1\}$ to choose from. While we only have access to one of the two parameters, we have access to more samples than a single one. 
Thus, the classical setup is sample-focused, while ours is parameter-focused.

Finally, it should be mentioned that our problem shares some similarities with the general two-sample testing problems (see, e.g.,~\cite{cai2014two,cai2013two,fromont2012kernels} and references therein), where the statistician is given two independent samples and is asked to determine whether they come from the same distribution---in other words, to test~$\cH_0: \Prob_0 = \Prob_1$ against~$\cH_1: \Prob_0 \ne \Prob_1$. Here the main distinction from our setup is the absence of the knowledge of~$\theta^* \in \{\theta_0, \theta_1\}$. This knowledge turns out to be instrumental in the model discrimination problem, allowing to circumvent the estimation of~$\theta_0,\theta_1$ and thus ultimately leading to confident testing im the small-sample regime~\eqref{eq:small-sample-regime-intro}. To the best of our knowledge, analogous phenomena have not been observed in the context of two-sample testing.

\paragraph{Notations.} 
We use~$C,c,c',c_1,...$ for generic positive constants. 
We write~$g = O(f)$ or~$f = \Omega(g)$ to state that~$g(t) \le Cf(t)$ uniformly over all possible values of~$t$.
Notation~$f \ll g$ (or~$g \gg f$) is the negation of~$g = O(f)$. 
We let~$[n] := \{1, 2, ..., n\}$ for~$n \in \mathds{N}$. We use~$\| \cdot \|$ for the~$\ell_2$-norm of a vector and the operator norm (the largest eigenvalue) of a positive-semidefinite (PSD) matrix. We use bold capitals for PSD matrices and regular capitals for rectangular matrices.
$\M^\dagger$ is the generalized inverse of a PSD matrix~$\M$, i.e., the matrix with the same eigenbasis as~$\M$ and inverse non-zero eigenvalues.
$\M^{\dagger/2}$ is the PSD square root of~$\M^\dagger$. 
We denote with~$\col(A)$ the column space of~$A \in \R^{n \times d}$, i.e., the linear span of the columns of~$A$.
$A \in \R^{n \times d}$ is said to have full column rank when~$\dim(\col(A)) = n$.
We let~$\|u\|_A := \|\bPi_{A} u\|^2$, where~$\bPi_A := A(A^{\top}A)^{\dagger}A^{\top}$ is the projector on~$\col(A)$.
For $r\leq d$, we denote with~$\Id_r \in \R^{d \times d}$ the diagonal matrix with $1$ on the first $r$ diagonal coordinates and $0$ otherwise.
We use ``Matlab notation'' for matrix concatenation:~$[A,B]$ (resp.,~$[A;B]$) is the horizontal (resp., vertical) concatenation of~$A$ and~$B$ with compatible dimensions. 
\section{Well-specified linear models}
\label{sec:gaussian}

In this section, we consider well-specified linear models -- in other words, the case of linear prediction (as defined in~Sec.~\ref{sec:intro}) with~$\ell_z(\theta) = (y - x^\top \theta)^2$ and Gaussian conditional distributions of~$y|x$ such that
\begin{equation}
\label{eq:linear-model}
y = x^\top \theta_{k} + \varepsilon, \;\; \varepsilon \sim \cN(0,1) \quad (k \in \{0,1\})
\end{equation}
with~$\varepsilon$ independent from~$x$. 
For convenience, let us denote the corresponding conditional distribution of $y|x$~(i.e.,~$\cN(x^\top \theta_k, 1)$) by~$\Prob_{k}[\cdot|x]$. 
In order to highlight the intuition behind our approach, we begin with a fixed-design scenario, so that all randomness of $y$ stems from the noise term $\varepsilon$.
Later on, we address the random-design scenario, thus returning to our canonical problem formulation as in~Sec.~\ref{sec:intro}. 
This requires an additional step of marginalizing over~$x$ in the analysis, and results in the appearance of the ranks~$r_0, r_1$ of {\em population covariances}~$\bSigma_0, \bSigma_1$ in the error bound and sample complexity estimate.
Another simplification in~\eqref{eq:linear-model} is the implicit assumption of the uniform noise variance under both models, i.e.,~$\varepsilon \sim \cN(0,1)$ instead of~$\varepsilon \sim \cN(0, \sigma_k^2)$ for~$k \in \{0,1\}$.  
When variances~$\sigma_0^2 ,\sigma_1^2$ are different but {\em known}, we can simply rescale the two samples dividing~$(X^{(k)},Y^{(k)})$ by~$\sigma_k$; 
however, this is impossible without knowledge of~$\sigma_0$ and~$\sigma_1$. 
Therefore, in Section~\ref{sec:gaussian-adaptive} we shall construct a test which is adaptive to the unknown noise variances. 

\subsection{Basic test and statistical guarantee}
\label{sec:gaussian-upper}

Let~$(x_i{}^{(k)}, y_i{}^{(k)}),$ for~$k \in \{0,1\}$ and~$i \in [n_k]$, be two samples with~$y_1{}^{(k)}, ..., y_{n_k}{}^{(k)}$ distributed independently according to~$\Prob_k[\cdot|x_i{}^{(k)}]$ as in~\eqref{eq:linear-model}. 
Writing them in a matrix form as~$X^{(k)} \in \R^{n_k \times d}$ and~$Y^{(k)} \in \R^{n_k}$,  
we first consider the test 
\begin{equation}
\label{eq:test_1}
\hat{T} = \ind \left\{ \| Y^{(0)} - X^{(0)} \theta^*\|^2_{X^{(0)}} - \rank(\wh\bSigma_0) \geq \| Y^{(1)} - X^{(1)} \theta^* \|^2_{X^{(1)}} - \rank(\wh \bSigma_1) \right\}.
\tag{$\textsf{\textup{Lin}}$}
\end{equation}
(Recall that~$\|u\|_A := \|\bPi_{A} u\|^2$ for~$u \in \R^n$, where~$\bPi_A := A(A^{\top}A)^{\dagger}A^{\top}$ is the projector onto the column space~$\col(A)$ of~$A \in \R^{n \times d}$.)
Test~\eqref{eq:test_1} combines two ideas.
First, the residuals are projected onto their respective signal subspaces~$\col(X^{(0)})$ and~$\col(X^{(1)})$. 
This reduces the noise variance without affecting the signal magnitude. 
Second, instead of directly comparing the squared norms of (projected) residuals, we compare their deviations from the corresponding sample covariance ranks. 
A direct calculation shows that the sample covariance rank is precisely the conditional (on the design) expectation of the corresponding residual under the matching hypothesis. 
Thus, subtracting the ranks debiases the squared norms and improves statistical performance of the resulting test.

In order to quantify the statistical performance of~\eqref{eq:test_1}, we introduce two measures of separation:
\begin{equation}
\label{def:Delta-empirical}
\wh\Delta_k = \|\wh\bSigma_k^{1/2}(\theta_1 - \theta_0)\|^2 \quad (k \in \{0,1\})
\end{equation}
where~$\wh \bSigma_k = \tfrac{1}{n_k} {X^{(k)\top}} X^{(k)}$ is the corresponding sample covariance matrix.
Essentially,~$\wh\Delta_k$ controls the separation between~$\Prob_0$ and~$\Prob_1$ from the viewpoint of the $k$-th model.
Now our first result follows.

\begin{proposition}
\label{prop:upper-fixed-design}
The type I error probability of test~\eqref{eq:test_1} is bounded as
\begin{equation}
\label{eq:upper-fixed-design}
\Prob_{\cH_0}[\hat{T} = 1|X^{(0)}, X^{(1)}] \le C \exp \bigg(-c n_1 \wh\Delta_1  \min\bigg\{1, \frac{n_1\wh\Delta_1}{\max\{\rank(\wh \bSigma_0), \rank(\wh \bSigma_1)\}} \bigg\}\bigg)\,,
\end{equation}
for some constants~$c,C>0$. The type II error admits a similar bound with the replacement~$1 \mapsto 0$.
\end{proposition}
We now sketch the proof of this result. The full proof, as well as the proofs of subsequent results, is deferred to appendix.
Denote~$\wh r_k := \rank(\wh \bSigma_k)$, and let~$\xi^{(0)} , \xi^{(1)}$ be the additive noise vectors:
\[
\xi^{(k)} := Y^{(k)} - X^{(k)} \theta_k \sim \cN(0,\Id_{n_k}).
\]
Observing that~$\| X^{(1)}(\theta_0 - \theta_1)\|^2 = n_1\wh\Delta_1$ we bound the type-I error probability as follows:
\begin{align*}
&\Prob_{\cH_0}[\hat{T} = 1 | X^{(0)}, X^{(1)}] \\
&= \Prob \big[ \| \xi^{(0)} \|^2_{X^{(0)}} - \wh r_0  + \wh r_1 - \| \xi^{(1)} \|^2_{X^{(1)}} + 2\langle X^{(1)}(\theta_0 - \theta_1), \xi^{(1)} \rangle \geq n_1\wh \Delta_1 \big]  \\ 
&\leq \Prob \bigg[ \| \xi^{(0)} \|^2_{X^{(0)}} - \wh r_0 \hspace{-0.07cm} \geq \hspace{-0.07cm} \frac{n_1 \wh\Delta_1}{3}\bigg] 
	+ \Prob \bigg[ \wh r_1 - \| \xi^{(1)} \|^2_{X^{(1)}} \hspace{-0.07cm} \geq \hspace{-0.07cm} \frac{n_1 \wh\Delta_1}{3}\bigg] 
	+ \Prob \bigg[ \langle X^{(1)}(\theta_0 - \theta_1), \xi^{(1)} \rangle \geq \frac{n_1 \wh\Delta_1}{6} \bigg].
\end{align*}
The first two terms are controlled through the standard chi-squared tail bounds (see, e.g.,~\cite[Lem.~1]{laurent2000adaptive}). 
Moreover,~$\langle X^{(1)}(\theta_0 - \theta_1), \xi^{(1)} \rangle \sim \cN(0, n_1\wh \Delta_1)$ conditionally on~$X^{(1)}$,
and the last deviation probability can be controlled via the standard Gaussian tail bound.
Combining these results yields~\eqref{eq:upper-fixed-design}.
\qed

\begin{remark}
\label{rem:heavy-tails}
Proposition~\ref{prop:upper-fixed-design} can be extended to the case where the additive noise is subgaussian. 
We here avoid such an extension, as in Section~\ref{sec:glm} we establish more general results that hold for generalized linear models under weak moment assumptions, and allow for model misspecification.
\end{remark}

Our next goal is to extend Proposition~\ref{prop:upper-fixed-design} to the random-design setup, conforming to the scenario of repeated i.i.d.~observations in Section~\ref{sec:intro}. 
We specify two distributions~$\Prob_{0},\Prob_1$ of~$(x,y)$ as follows:
\begin{equation}
\label{eq:gaussian-distributions}
\begin{aligned}
\Prob_k: x \sim \cN(0,\bSigma_k), \; y|x &\sim \cN(x^\top \theta_{k},1) \quad (k \in \{0,1\}).
\end{aligned}
\end{equation}
Furthermore, we let~$(X^{(k)},Y^{(k)})$ represent the i.i.d.~sample from~$\Prob_k$ with size~$n_k$ ($k \in \{0,1\}$).
Finally, we let~$r_k = \rank(\bSigma_k)$ for~$k \in \{0,1\}$.
Our test is still~\eqref{eq:test_1}, but we now additionally observe that
\[
\wh r_k = \rank(\wh \bSigma_{k}) = \min \{n_k, r_k\} \quad \text{\em a.s.}
\]
(More generally, this is the case when the marginal distribution of~$x$ has density on its support). 
Next we characterize the sample complexity of reliable testing in terms of the {\em population separations}~$\Delta_{0},\Delta_{1}$,
\begin{equation}
\label{def:Delta-population}
\Delta_k := \E_k [\wh\Delta_k] = \|\bSigma_k^{1/2}(\theta_1 - \theta_0)\|^2 \quad (k \in \{0,1\}),
\end{equation}
i.e., the squared Mahalanobis distances between~$\theta_0,\theta_1$ associated with  the two covariance matrices. 
\begin{theorem}
\label{th:upper-random-design}
In the random-design setup specified above, the type I error of test~\eqref{eq:test_1} satisfies 
\begin{equation}
\label{eq:upper-random-design}
\Prob_{\cH_0}[\hat{T} = 1] \le C \exp \bigg(-c n_1\Delta_1 \min \bigg\{ 1, \frac{n_1\Delta_1}{\bar r} \bigg\} \bigg) + C\exp(-cn_1),
\end{equation}
\begin{equation}
\label{eq:r-bar}
\text{where} \;\;
\bar r := \max\{\min\{n_0, r_0\}, \min \{n_1,r_1\} \}.
\end{equation}
The type II error probability admits a similar bound with the replacement~$1 \mapsto 0$.
\end{theorem}
In the proof of Theorem~\ref{th:upper-random-design}, we first condition on~$X^{(0)}, X^{(1)}$ and repeat the analysis leading to Proposition~\ref{prop:upper-fixed-design}. 
Then we add a marginalization step, which leads to the result. 
The term~$\exp(-cn_1)$ in~\eqref{eq:upper-random-design} appears in this step.
When~$\Delta_1 \le c$, this term is dominated, and the guarantee is essentially the same as in the fixed-design scenario, except that~$\wh\Delta_1$ gets replaced with its expectation~$\Delta_k$.

Next we examine Theorem~\ref{th:upper-random-design} in more detail and discuss its implications.
For the sake of simplicity, until the end of this section we assume~$n_0 = n_1 [:= n]$ and~$\max\{\Delta_{0}, \Delta_{1}\} \le 1$.
Clearly, the latter assumption is rather mild as it only excludes very simple problems which are of little interest anyway.
We also define~$r :=  \max\{r_0, r_1\}$, so that~$\bar r = \min\{r, n\}$.

First, with some algebra we translate~\eqref{eq:upper-random-design} 
into the sample complexity bound~\eqref{eq:intro-complexity} announced in Section~\ref{sec:contrib}. 
More precisely, we establish the following result.
\begin{corollary}
\label{cor:complexity-gaussian}
Assuming~$\max\{\Delta_0,\Delta_1\} \le 1$, the sample complexity of distinguishing between two linear models, i.e.,~$n = n_0 = n_1$ that guarantees~$\max\{\Prob_{\cH_0}[\hat{T} = 1],\Prob_{\cH_1}[\hat{T} = 0]\} \le \delta < 1/2$, is
\begin{equation}
\label{eq:complexity-gaussian}
n = O\left(\min\left\{\frac{\log(1/\delta)}{\Delta^2}, \frac{\sqrt{r \log(1/\delta)}+\log(1/\delta)}{\Delta} \right\} \right),
\end{equation}
where~$\Delta = \min\{\Delta_0, \Delta_1\}$ and $r= \max\{r_0,r_1\}$. 
This reduction is tight: conversely,~when~$\max \{\Delta_0,\Delta_1\} \le 1$, bounding the testing errors of both types simultaneously by~$\delta$, as per~\eqref{eq:upper-random-design}, requires~$n$ to be at least as in~\eqref{eq:complexity-gaussian}.
\end{corollary}
Corollary~\ref{cor:complexity-gaussian} is proved in Appendix~\ref{app:upper-linear}.
From this result we see that the sample complexity of testing with fixed confidence is~$O(\min\{1/\Delta^2,\sqrt{r}/\Delta\})$.
In particular, in the {\em well-separated regime} 
\begin{equation}
\label{eq:high-dim-range}
\Delta = \Omega\left(\frac{1}{\sqrt{r}}\right),
\end{equation}
the sample complexity becomes rank-independent, scaling as~$n = \Omega(1/\Delta^{2})$. 
An intuitive explanation of this result is that testing does {\em not} require to estimate the population risk minimizers~$\theta_0,\theta_1$---which would inevitably result in a rank-dependent complexity term.
In fact, our analysis in the proofs of Proposition~\ref{prop:upper-fixed-design} and~Theorem~\ref{th:upper-random-design} shows that the~$O(1/\Delta^2)$ complexity term appears because of the random fluctuations of a centered chi-squared noise around its mean~$\bar r = \min\{r,n\}$, which does not exceed the population rank~$r$. 
On the procedural level, we achieve this effect in~\eqref{eq:test_1} by projecting the residuals onto the column spaces~$\col(X^{(0)}),\col(X^{(1)})$ with dimensions~$\wh r_k \le r_k$,~$k \in \{0,1\}$.

In what follows next, we discuss the sample complexity result~\eqref{eq:complexity-gaussian} in a somewhat broader context. 

\subsection{Discussion of implications}
\label{sec:discussion}

\paragraph{Key implication: model non-disclosure.}
In the context of the testing protocol discussed in~Section~\ref{sec:applications}, 
test~\eqref{eq:test_1} provides a guarantee that the model~$\theta^*$ trained by~\learner{} cannot be inferred by \auditor{} (in the variant of the protocol where the two entities are separated).
Indeed, as we discussed in Section~\ref{sec:applications}, we cannot hope to estimate~$\theta^*$ up to prediction error~$\Delta$ (in the covariance matrix corresponding to the true hypothesis) without estimating at least one of the two parameters~$\theta_k$ ($k \in \{0,1\}$) with prediction error~$\Delta_k$, i.e., such that~$\| \bSigma_k{}^{1/2}(\wh\theta_k - \theta_k)\|^2 \le \Delta_k$. 
When~$\theta_k$ is arbitrary, this is known to require~$\Omega(r_k/\Delta_k)$ sample size for any estimator (see, e.g.,~\cite{tsybakov_mono}).
In the ``typical'' situation where~$\Delta_0, \Delta_1$, as well as~$r_0,r_1$, are of the same order, this lower bound on sample complexity reduces to~$\Omega(r/\Delta)$, which orderwise (in terms of $\Delta$) dominates the upper bound~$\min\{1/\Delta^2, \sqrt{r}/\Delta\}$ in Corollary~\ref{cor:complexity-gaussian}.
Thus, one can choose the sample sizes that suffice for confident testing, yet do not allow to recover~$\theta_0, \theta_1$ up to~$O(\Delta)$ prediction accuracy.

In fact, this conclusion extends beyond the case of well-specified linear models to other scenarios.
Indeed, for general~$M$-estimators, the lower bound for the complexity of estimating~$\theta_k$ (with respect to the excess risk, i.e., such that~$L_k(\wh \theta_k) - L_k(\theta_k) \le \Delta_k$) is generally as bad as~$\Omega(r_k/\Delta_k)$: this is already so in the case of linear models, in which~$L_k(\wh \theta_k) - L_k(\theta_k) = \tfrac{1}{2}\| \bSigma_k^{1/2}(\wh\theta_k - \theta_k)\|^2$.
In particular, estimation with prediction accuracy~$\Delta_k \le 1$ is generally impossible when~$n_k \le r_k$. 
On the other hand, in Section~\ref{sec:glm} we show that the rank-independent bound on the sample complexity of testing,
\[
O\left(\frac{1}{\Delta^2}\right),
\] 
extends to generalized linear models under~\eqref{eq:small-sample-regime-intro} provided that~$\Delta = \Omega(1/\sqrt{\min\{r_0,r_1\}})$ (see the discussion after Theorem~\ref{th:glm-fixed-prob}). In this regime~$1/\Delta^2 \ll r/\Delta$, i.e., testing is way easier than recovery.
Complementary to this result, in Section~\ref{sec:asymp} we consider the case of general~$M$-estimators in the large sample size regime (with~$n \to \infty, \Delta \to 0$, and fixed ranks)
and show that testing is possible whenever
\[
n = O\left(\frac{\sqrt{\reff}}{\Delta}\right)
\]
where~$\reff$ is the largest of the two effective ranks (under~$\Prob_0$ and~$\Prob_1$) reducing to~$r$ in the case of well-specified models (see~Proposition~\ref{prop:asymp} for the exact statement).
Typically, we expect~$\reff = O(r)$, and the sample complexity of testing is~$O(\sqrt{r}/\Delta)$, again way smaller than the complexity of recovery.

\paragraph{``Plug-in'' interpretation of test~\eqref{eq:test_1}.} 
There is an interesting alternative interpretation of test~\eqref{eq:test_1}. 
Observe that the normal equation for the least-squares estimate~$\wh \theta_k$ of~$\theta_k$ reads
\begin{equation}
\label{eq:normal-equations}
\wh\bSigma_k \wh \theta_k = \tfrac{1}{n_k} X^{(k)\top} Y^{(k)}.
\end{equation}
This allows to rewrite the statistic in~\eqref{eq:test_1}: omitting the subscript~$k$ and superscript~$(k)$ for brevity, 
\begin{align}
\|Y - X\theta^*\|^2_X  
= (Y - X\theta^*)^\top \bPi_X (Y - X \theta^*)
&= (X^\top Y - X^\top X \theta^*)^\top (X^\top X)^{\pinv} (X^\top Y - X^\top X \theta^*) \notag\\
&= n^2 (\wh\theta - \theta^*)^\top \wh\bSigma (X^\top X)^{\pinv} \wh\bSigma (\wh\theta - \theta^*) \notag\\
&= n (\wh\theta - \theta^*)^\top \wh\bSigma^{\vphantom\pinv} \wh \bSigma^{\pinv} \wh\bSigma^{\vphantom\pinv} (\wh\theta - \theta^*) \notag\\
&= n (\wh\theta - \theta^*)^\top \wh\bSigma (\wh\theta - \theta^*) \notag\\
&= n \| \wh\bSigma^{1/2} (\wh\theta - \theta^*) \|^2.
\label{eq:residuals-to-plug-in}
\end{align}
As a result, we recast test~\eqref{eq:test_1} as
\begin{equation}
\label{eq:test_1-as-plug-in}
\wh T = \ind\big\{ n_0 \big\|\wh\bSigma_0^{1/2}(\theta^* - \wh\theta_0)\big\|^2 - \wh r_0 \ge n_1 \big\|\wh\bSigma_1^{1/2}(\theta^* - \wh\theta_1)\big\|^2 - \wh r_1 \big\}.
\end{equation}
Now, assume (w.l.o.g.) that~$\bSigma_0,\bSigma_1$ are both full-rank, and let~$n_k \ge r_k$ ($k \in \{0,1\}$) so that~$\wh\theta_0, \wh \theta_1$ are defined unambiguously as the unique solutions to~\eqref{eq:normal-equations}.
In this regime, we can interpret~\eqref{eq:test_1-as-plug-in} as a test that constructs plug-in estimates of the squared empirical prediction distances~$\|\wh\bSigma_0{}^{1/2}(\theta^* - \theta_0)\|^2$ and~$\|\wh \bSigma_1{}^{1/2}(\theta^* - \theta_1)\|^2$ of~$\theta^*$ to~$\theta_0$ and to~$\theta_1$, then rescales and debiases these estimates conditionally on~$X^{(0)}, X^{(1)}$ under the matching hypothesis, and compares the results. 
However, equations~\eqref{eq:normal-equations} characterize the solutions~$\wh\theta_0, \wh\theta_1$ to the corresponding least-squares problems even in the small-sample regime, where these solutions are not unique, let alone concentrate around~$\theta_0, \theta_1$. 
Yet, the squared norms used in~\eqref{eq:test_1-as-plug-in} are still well-defined due to~\eqref{eq:residuals-to-plug-in}. 
In this sense, the ``plug-in'' interpretation (cf.~\eqref{eq:test_1-as-plug-in}) of test~\eqref{eq:test_1} is limited and can be misleading.

However, the ``plug-in'' interpretation~\eqref{eq:test_1-as-plug-in} is still useful since it provides an alternative explanation of the terms~$\wh r_0, \wh r_1$ appearing in~\eqref{eq:test_1-as-plug-in}. 
Indeed, instead of~\eqref{eq:test_1-as-plug-in} one could think of using a simpler test
\begin{equation}
\label{eq:pred-dist-comp}
\wh T = 
\ind\big\{ \big\|\wh\bSigma_0^{1/2}(\theta^* - \wh\theta_0)\big\|^2 \ge \big\|\wh\bSigma_1^{1/2}(\theta^* - \wh\theta_1)\big\|^2\big\},
\end{equation}
i.e., directly compare the plug-in estimates of the squared prediction distances~$\big\|\bSigma_k^{1/2}(\theta^* - \theta_k)\big\|^2$ ($k \in \{0,1\}$) one of which vanishes under the matching hypothesis. 
However, this approach does not take into account the additional noise arising due to the random fluctuations of~$\wh \theta_k$ around~$\theta_k$.
To illustrate this issue, it suffices to consider an idealized test with the population covariances, namely
\begin{equation}
\label{eq:pred-dist-comp-ideal}
\wh T = 
\ind\big\{ \big\|\bSigma_0^{1/2}(\theta^* - \wh\theta_0)\big\|^2 \ge \big\|\bSigma_1^{1/2}(\theta^* - \wh\theta_1)\big\|^2\big\}.
\end{equation}
Note that~$\E[\|\bSigma_k^{1/2} (\wh \theta_k - \theta_k)\|^2] = {r_k}/{n_k}$; thus, under~$\cH_0$ the expectation of the left-hand side in \eqref{eq:pred-dist-comp-ideal} is~$r_0/n_0$ and of the right-hand side is~$\Delta_1 + r_1/n_1$. 
As a result, in the case~$n_0 = n_1 [= n]$, the type I error for~\eqref{eq:pred-dist-comp-ideal} cannot be controlled when~$n < (r_0 - r_1)/\Delta_1$; similarly, the type II error cannot be controlled when~$n < (r_1 - r_0) / \Delta_0$; overall, we cannot control at least one of these errors when~$n < |r_1 - r_0| / \Delta$. 
Hence, when the ranks~$r_0, r_1$ are significantly different, and when~$n_0 = n_1$, the sample complexity for~\eqref{eq:pred-dist-comp-ideal} is at least~$\Omega(r/\Delta)$ -- way larger than the upper bound~$O(\min\{\sqrt{r}/\Delta, 1/\Delta^2\})$ for~\eqref{eq:test_1}. 

It is also clear that the above-described issue cannot be solved by choosing the sample sizes proportionally to the ranks (when this is allowed, and when the ranks are known) whenever~$\min\{r_0,r_1\}$ is way smaller than~$\max\{r_0,r_1\}$. Indeed, while such an approach would guarantee that~$r_0/n_0 = r_1/n_1$ and thus correct for the additional bias due to~$r_0 \ne r_1$, it would also result in increased magnitude of fluctuations, compared to the choice~$n_0 = n_1$, as the sample size corresponding to the smallest rank must also be very small.
Finally, when~$r_0 = r_1$ and~$n_0 = n_1$, there is no issue with additional bias to begin with, and in this case~\eqref{eq:pred-dist-comp} reduces to~\eqref{eq:test_1}.

\paragraph{Necessity of projection.}
Consider the small sample regime~$n \le \min\{r_0, r_1\}$ (cf.~\eqref{eq:small-sample-regime-intro}). 
In this regime,~$X^{(0)}$ and~$X^{(1)}$ have full column ranks~($n_0$ and~$n_1$ respectively)
and test~\eqref{eq:test_1} simplifies to
\begin{equation}
\label{eq:test_1-small-sample}
\wh T = \ind \{ \| Y^{(0)} - X^{(0)}\theta^*\|^2 - n_0 \geq \| Y^{(1)} - X^{(1)} \theta^* \|^2 - n_1\}.
\end{equation}
This simple observation has several implications.

\begin{itemize}
\item
In the well-separated regime (cf.~\eqref{eq:high-dim-range}), the sample complexity upper bound for~\eqref{eq:test_1}, cf.~\eqref{eq:complexity-gaussian}, falls into the range~$n \le r$.
When~$r_0$ and~$r_1$ are of the same order, confident testing in the well-separated regime thus only requires the sample size~$n = O(\min\{r_0,r_1\})$.
One can then run~\eqref{eq:test_1} with such sample sizes to begin with, thus effectively reducing it to~\eqref{eq:test_1-small-sample} and avoid the computational burden of performing projections. 
\item
On the other hand, the projection step is crucial for {\em ill-separated} problems, i.e., when~$\Delta \ll 1/\sqrt{r}$. 
Indeed, in this case the sample complexity bound in~\eqref{eq:complexity-gaussian} falls {\em beyond} the range~$n \le r$ and test~\eqref{eq:test_1-small-sample}, which does not project onto $\col(X^{(0)})$ and $\col(X^{(1)})$, becomes suboptimal.
More precisely, one can easily verify (by mimicking the decomposition in the proof of Proposition~\ref{prop:upper-fixed-design}) that the sample complexity for test~\eqref{eq:test_1-small-sample}
is~$O(\max\{1/\Delta^2, 1/\Delta\})$, i.e.,~$O(1/\Delta^2)$ unless~$\Delta \gg 1$, whereas for~\eqref{eq:test_1} it is~$O(\sqrt{r}/\Delta) \ll 1/\Delta^2$. 
\item
In the general~$M$-estimation setup (as defined in~\eqref{eq:hypotheses}-\eqref{eq:samples}), 
we naturally generalize test~\eqref{eq:test_1-small-sample} to
\begin{equation}
\label{eq:test_naive}
\hat{T} = \ind\big\{ n_0 \wh L_{0}(\theta^*) - n_0 L_0(\theta_0) \geq n_1 \wh L_{1}(\theta^*) - n_1 L_1(\theta_1)\big\},
\tag{$\textsf{\textup{Val}}$}
\end{equation}
where~$\wh L_{k}(\theta)$ and~$L_k(\theta)$ are, respectively, the empirical and population risks~($k \in \{0,1\}$). 
An immediate problem with this generalization is that the population risks~$L_k(\theta_k)$ are generally unobservable, unless in the exceptional case of a well-specified linear model where~$L_k(\theta_k) = 1$. 
Thus, test~\eqref{eq:test_naive} generally cannot be implemented.

Meanwhile, in Sections~\ref{sec:asymp}-\ref{sec:glm} we shall see that test~\eqref{eq:test_1} does generalize beyond the setup of a well-specified linear model,
as the corresponding debiasing terms (replacing~$\wh r_0, \wh r_1$ in~\eqref{eq:test_1}) can be estimated from observations. 
Moreover, the reduction of~\eqref{eq:test_1} to~\eqref{eq:test_naive} for small sample sizes is also limited to the case of a well-specified linear model.
This is because the terms generalizing~$\| Y^{(0)} - X^{(0)} \theta^*\|^2_{X^{(0)}}$ and~$\| Y^{(1)} - X^{(1)} \theta^* \|^2_{X^{(1)}}$ in~\eqref{eq:test_1} will stem from the second-order Taylor approximation of empirical risks rather than the empirical risks themselves.
\end{itemize}

In the next section, we consider the ``heterogeneous'' scenario in which the noise variance depends on the hypothesis and is unknown. 
We construct an adaptive test that performs essentially as good as in the case of known variances, albeit this requires an extra assumption in small-sample regime.

\subsection{Adaptation to unknown noise levels}
\label{sec:gaussian-adaptive}

Recall that test~\eqref{eq:test_1}, in its precise form, is designed for the case of unit variance of additive noise under both hypotheses, cf.~\eqref{eq:linear-model}. 
If this assumption is not met, i.e.,~$\Prob_k[y|x]$ reads
\begin{equation}
\label{eq:linear-model-var}
y = x^\top \theta_{k} + \sigma_k \varepsilon, \;\; \varepsilon \sim \cN(0,1) \quad (k \in \{0,1\}),
\end{equation}
instead of~\eqref{eq:linear-model}, then test~\eqref{eq:test_1} can be replaced with the following one:
\begin{equation}
\label{eq:test_1-var}
\hat{T} = \ind \left\{ \frac{1}{\sigma_0^2} \| Y^{(0)} - X^{(0)} \theta^*\|^2_{X^{(0)}} - \rank(\wh\bSigma_0) \geq 
\frac{1}{\sigma_1^2} \| Y^{(1)} - X^{(1)} \theta^* \|^2_{X^{(1)}} - \rank(\wh \bSigma_1) \right\}.
\tag{$\textsf{\textup{Lin-Var}}$}
\end{equation}
Since~\eqref{eq:test_1-var} reduces to~\eqref{eq:test_1} when each sample is rescaled by~$1/\sigma_k$, the results in Section~\ref{sec:gaussian-upper} generalize for this test, with~$\wh\Delta_k$ (and correspondingly~$\Delta_k$) rescaled by~$1/\sigma_k^2$. 
However, such a reduction requires the knowledge of both noise variances~$\sigma_0^2, \sigma_1^2$. 
Next we construct an adaptive counterpart of test~\eqref{eq:test_1-var} that admits similar statistical guarantees without such knowledge.
To that end, we construct some estimates~$\wh \sigma_0^2, \wh \sigma_1^2$ of~$\sigma_0^2,\sigma_1^2$ and plug them in~\eqref{eq:test_1-var}, 
thus arriving at
\begin{equation}
\label{eq:test_1-adapt}
\hat{T} = \ind \left\{ \frac{1}{\wh\sigma_0^2} \| Y^{(0)} - X^{(0)} \theta^*\|^2_{X^{(0)}} - \rank(\wh\bSigma_0) \geq 
\frac{1}{\wh\sigma_1^2} Y^{(1)} - X^{(1)} \theta^* \|^2_{X^{(1)}} - \rank(\wh \bSigma_1) \right\}.
\tag{$\textsf{\textup{Lin-Var$^{+}$}}$}
\end{equation}
Assuming~$n_0 = n_1 = n$ for the sake of simplicity, we construct~$\wh\sigma_k^2$ ($k \in \{0,1\}$) in two different ways, depending on how large is~$n$ compared to~$\wh r_k = \rank(\wh\bSigma_k)$. 
Namely, if~$n$ is large enough, then we construct estimates with sufficient accuracy without additional assumptions;
however, in the small-sample regime we in addition assume that~$y$ can be resampled conditionally on~$x$ (once).
\begin{assumption}
\label{ass:resampling}
One can sample~$\wt Y^{(k)}$ such that~$\wt Y^{(k)}, Y^{(k)}$ are i.i.d.~conditionally on~$X^{(k)}$ ($k \in \{0,1\}$).
\end{assumption}
Admittedly, this assumption is rather strong. However, we only needed it in the small-sample regime.
We now present the variance estimates (we omit the subscript in~$X^{(k)}$ and~$Y^{(k)}$ for brevity).
\begin{itemize}
\item {\bf Case~$n \ge 2 \wh r_k$} (note that this is equivalent to~$n \ge 2r_k$ since~$\wh r_k = \min\{n,r_k\}$ {\em a.s}). 
Here we put
\begin{equation}
\label{eq:var-estimate-simple}
\wh\sigma_k^2 = \frac{\| Y - X \theta^*\|^2 - \| Y - X \theta^*\|^2_{X}}{n - \wh r_k}.
\end{equation}
The idea is to estimate the noise variance from the component of the residual which is orthogonal to the signal subspace and thus only contains the noise. 
\item {\bf Case~$n \le 2\wh r_k$} (similarly, this is equivalent to~$n \le 2r_k$). 
Here we leverage Assumption~\ref{ass:resampling}
and let
\begin{equation}
\label{eq:resampling-simple}
\wh\sigma_k^2 
= \frac{\|{\wh \bSigma_{k}}^{\pinv/2} X^\top [Y - \wt Y]\|^2}{2n\wh r_k}.
\end{equation}
Here the idea is to project onto the signal subspace while eliminating the signal via resampling.
Note that whenever~$n \le r_k$, estimate~\eqref{eq:resampling-simple} reduces to~$\tfrac{1}{2n} \|Y - \wt Y \|^2$, as in this case~$\wh r_k = n$ and
\[
\frac{\|\wh \bSigma_{k}^{\pinv/2} X^\top (Y - \wt Y)\|^2}{2n\wh r_k} = 
\frac{\|\wh \bSigma_{k}^{\pinv/2} X^\top (Y - \wt Y)\|^2}{2n^2} 
= \frac{\|\Pi_{X}[Y - \wt Y]\|^2}{2n} 
= \frac{\|Y - \wt Y\|^2}{2n}
\]
where the last step uses that~$X$ has full column rank.
Note also that, strictly speaking, we do not need access to~$\wt Y$ itself -- rather, to the right-hand side of~\eqref{eq:resampling-simple} or to~$\|Y- \wt Y\|^2$ when~$n \le r_k$.
\end{itemize}
The key property in both these constructions is that the normalized estimate~$\wh\tau_k := \wh\sigma_k^2/\sigma_k^2$ satisfies
\begin{equation}
\label{eq:var-ratio-tail-bound}
\Prob\left[\left|\wh\tau_k - 1 \middle| X \right| \ge t\right] 
\leq
2\exp\left(-c n \min\{t^2, t\}\right).
\end{equation}
Indeed, for~\eqref{eq:var-estimate-simple} the conditional distribution of~$\wh\tau_k$ given~$X^{(k)}$ is~$\chi_{n-\wh r_k}^2/(n-\wh r_k)$, hence by the standard~$\chi^2$ tail bound (see, e.g.,~\eqref{eq:tail-bounds-chi-sq-right}--\eqref{eq:tail-bounds-chi-sq-left} in the appendix) we have
\[
\Prob\left[\left|\wh\tau_k - 1 \right| \ge t \middle| X^{(k)} \right] 
\leq 
2\exp\left(-2c(n-\wh r_k) \min\{t^2, t\}\right) 
\leq
2\exp\left(-c n \min\{t^2, t\}\right).
\]
On the other hand, in the case of~\eqref{eq:resampling-simple} we have~$\wh\tau_k \sim \chi_{\wh r_k}^2/\wh r_k$.
Since~$2\wh r_k > n$, we again arrive at~\eqref{eq:var-ratio-tail-bound}.

Applying~\eqref{eq:var-ratio-tail-bound} to~$\wh\tau_k$ for~$k \in \{0,1\}$ we can show that the resulting test~\eqref{eq:test_1-adapt} admits similar statistical guarantees to those for the oracle test~\eqref{eq:test_1-var}.
To this end, we first extend Proposition~\ref{prop:upper-fixed-design}.
\begin{proposition}
\label{prop:upper-fixed-design-adapt}
The type I error probability of test~\eqref{eq:test_1-adapt} is bounded as
\begin{equation}
\label{eq:upper-fixed-design-adapt}
\Prob_{\cH_0}[\hat{T} = 1|X^{(0)}, X^{(1)}] \le C \exp \bigg(-\frac{cn \wh\Delta_1}{\sigma_1^2} \min\bigg\{1, \frac{n\wh\Delta_1}{\sigma_1^2 \wh r} \bigg\}\bigg) + C\exp(-cn)
\end{equation}
where~$\wh r = \max\{\rank(\wh \bSigma_0), \rank(\wh \bSigma_1) \}$.
The type II error admits the bound with the replacement~$1 \mapsto 0$.
%
\end{proposition}
The full proof of this result is deferred to the appendix, and here we sketch it.
We have
\begin{align}
&\; 
\Prob_{\cH_0}[\hat{T} = 1|X^{(0)}, X^{(1)}] \notag\\
&= \Prob \bigg[ \frac{1}{\wh \tau_0} \|\xi^{(0)}\|^2_{X^{(0)}} - \wh r_0  + \wh r_1 - \frac{1}{\wh \tau_1}\| \xi^{(1)} \|^2_{X^{(1)}} + \frac{2}{\wh \tau_1 \sigma_1} \langle X^{(1)} (\theta_1 - \theta_0), \xi^{(1)} \rangle \geq \frac{n\wh \Delta_1}{\wh\tau_1 \sigma_1^2} \bigg]  \notag\\ 
&= \Prob \bigg[ \frac{\wh \tau_1}{\wh \tau_0} \|\xi^{(0)}\|^2_{X^{(0)}} - \wh\tau_1 \wh r_0  + \wh\tau_1 \wh r_1 - \|\xi^{(1)}\|^2_{X^{(1)}} + \frac{2}{\sigma_1} \langle X^{(1)} (\theta_1 - \theta_0), \xi^{(1)} \rangle \geq \frac{n\wh \Delta_1}{\sigma_1^2} \bigg]  \notag\\
%
&\leq \Prob \bigg[ \wh r_1 - \|\xi^{(1)}\|^2_{X^{(1)}}  \geq \frac{n \wh\Delta_1}{5 \sigma_1^2}\bigg] 
	+ \Prob \bigg[\frac{1}{\sigma_1} \langle X^{(1)} (\theta_1 - \theta_0), \xi^{(1)} \rangle \geq \frac{n \wh\Delta_1}{10 \sigma_1^2} \bigg] \notag\\
&\quad 
	+ \Prob \bigg[ \frac{\wh\tau_1}{\wh \tau_0} \left( \|\xi^{(0)}\|^2_{X^{(0)}} - \wh r_0 \right) \geq \frac{n \wh\Delta_1}{5 \sigma_1^2}\bigg]
	+ \Prob \bigg[ (\wh \tau_1 - 1) \wh r_1 \ge \frac{n \wh\Delta_1}{5 \sigma_1^2} \bigg] 
	+ \Prob \bigg[ \wh \tau_1 \bigg(\frac{1}{\wh \tau_0} - 1\bigg) \wh r_0 \ge \frac{n \wh\Delta_1}{5 \sigma_1^2} \bigg]
\notag
\end{align}
where all probabilities are conditional on~$X^{(0)},X^{(1)}$. 
The first two terms in the right-hand side can be bounded via the standard chi-squared and Gaussian tail bounds as in the proof of Proposition~\ref{prop:upper-fixed-design} since these terms do not contain~$\wh\tau_k$.
Now, it turns out that the remaining three terms also admit the bound from the right-hand side of~\eqref{eq:upper-fixed-design-adapt}. 
To obtain this result, we separate the factors containing~$\wh\tau_0, \wh\tau_1$ and use~\eqref{eq:var-ratio-tail-bound} to bound by~$C\exp(-cn)$ the probability for these factors to deviate from a constant. 
The resulting~$(\wh\tau_0,\wh\tau_1)$-free deviations are bounded as in~Proposition~\ref{prop:upper-fixed-design} for the first term in the last line,
and by Theorem~\ref{th:upper-random-design} for the last two terms.
\qed\\

Finally, repeating the steps in the proof of Theorem~\ref{th:upper-random-design} -- namely, conditioning over design and using that~$\Prob \left[ n_1\wh \Delta_1 \le {n_1 \Delta_1}/{2} \right] \le \exp(-cn_1)$ -- we extend~Proposition~\ref{prop:upper-fixed-design-adapt} to the random design setup. 

\begin{corollary}
\label{cor:upper-random-design-adapt}
Let~$\Prob_k$ ($k \in \{0,1\}$) correspond to~$x \sim \cN(0,\bSigma)$ and~$y|x \sim \cN(x^\top \theta^{\vphantom 2}_{k}, \sigma_k^2)$ as in the premise of Theorem~\ref{th:upper-random-design}.
Then test~\eqref{eq:test_1-adapt} satisfies 
\begin{equation}
\label{eq:upper-random-design-adapt}
\Prob_{\cH_0}[\hat{T} = 1] \le C \exp \bigg(-\frac{c n \Delta_1}{\sigma_1^2} \min\bigg\{1, \frac{n\Delta_1}{\bar r \sigma_1^2}\bigg\}\bigg) + C\exp(-cn),
\end{equation}
where~$\bar r$ is as in Theorem~\ref{th:upper-random-design} (cf.~\eqref{eq:r-bar}).
The type II error admits the bound with the replacement~$1\mapsto 0$.
\end{corollary}

\subsection{Lower bound}
\label{sec:gaussian-lower}

Our next result demonstrates that test~\eqref{eq:test_1} is near-optimal in the minimax sense.\vspace{-0.1cm}
\begin{theorem}
\label{th:lower-gaussian}
Let~$\Delta^2 > 0$ and integers~$r_0, r_1, n \ge 28$ be given.
Consider distributions~$\Prob_0, \Prob_1$ given by 
\[
\begin{aligned}
\Prob_k: x \sim \cN(0,\bSigma_k), \; y|x \sim \cN(x^\top \theta_{k},1) \quad (k \in \{0,1\}),
\end{aligned}
\]
parametrized by unknown~$\theta_0, \theta_1 \in \R^d$ with~$d \ge \max\{r_0, r_1\}$,
with covariances~$\bSigma_0 = \Id_{r_0}$ and~$\bSigma_1 = \Id_{r_1}$. 
Observing~$\theta^* \in \{\theta_0, \theta_1\}$ and two independent samples~$(X^{(0)},Y^{(0)}) \sim \Prob_0^{\otimes n}$ and~$(X^{(1)}, Y^{(1)}) \sim \Prob_1^{\otimes n}$, 
consider testing~$\cH_0: \theta_0 = \theta^*$ against~$\cH_1: \theta_1 = \theta^*$. 
Then
\begin{equation}
\label{eq:lower-gaussian}
\underset{\wh T}{\inf} \underset{(\theta_0,\theta_1) \in \Theta(\Delta^2)}{\sup} 
\left\{ \Prob_{\cH_0}[\wh T = 1] + \Prob_{\cH_1}[\wh T = 0] \right\} 
\ge C \exp \bigg(-c n \Delta \min\bigg\{1, \frac{n\Delta}{\textup{\rmin}}\bigg\}\bigg),
\end{equation}
where the~infimum is over all measurable maps~$T:(\theta^*, X^{(0)}, Y^{(0)}, X^{(1)}, Y^{(1)}) \to \{0,1\}$, the supremum is over the set~
$
\Theta(\Delta) := \{(\theta_0,\theta_1) \in \R^{2d}: \max\{\Delta_0, \Delta_1\} \geq \Delta\} 
$
with~$\Delta_{0}, \Delta_{1}$ as defined in~\eqref{def:Delta-population}, 
and
\[
\textup{\rmin} = \min\{n, r_0, r_1\}.
\]
\end{theorem}
Theorems~\ref{th:upper-random-design} and~\ref{th:lower-gaussian} together show that test~\eqref{eq:test_1} is optimal up to the replacement of~$\min\{\Delta_0,\Delta_1\}$ with~$\max\{\Delta_0,\Delta_1\}$ 
and~$\max\{r_0,r_1\}$ with $\min\{r_0,r_1\}$. 
Removing this gap remains an open problem.\\

The main technical challenge in the proof of Theorem~\ref{th:lower-gaussian} stems from not observing the residuals at the complementary to~$\theta^*$ model
$
\bar \theta = \theta_0 + \theta_1 -  \theta^*.
$
This leads to a {\em composite} hypothesis testing problem. 
To circumvent this challenge, we first consider the simplified situation where one observes the residuals corresponding to {\em both~$\theta^*$ and~$\bar\theta$}, rather than only~$\theta^*$, but one is not told ``which one is which.'' In other words, one observes the {\em unordered} pair of residuals.
The testing problem then reduces to inferring the  right {ordering} of~$\bar\theta$ and~$\theta^*$, i.e., testing two {\em simple} hypotheses~$(\theta^*, \bar\theta) = (\theta_0, \theta_1)$ and~$(\theta^*, \bar\theta) = (\theta_1, \theta_0)$. In this situation, the Neyman-Pearson lemma~\cite{lehmann2006testing} ensures that the likelihood ratio test is optimal. It turns out that this reduction already suffices to match the~$\rmin$-independent first term under exponent in~\eqref{eq:lower-gaussian}. 
However, it fails to capture the dependency on~$\rmin$ which is reflected in the second term.
In order to do that, we rely on a constrained Bayesian approach, putting a truncated Gaussian prior on~$\bar{\theta}$ in the same spirit as in \cite{ndaoud2019interplay,ndaoud2018sharp}.\\ 

In the next section, we generalize our approach to the general~$M$-estimation scenario described in Section~\ref{sec:intro}.
The key insight is that the projected residuals used in test~\eqref{eq:test_1} are in fact instantiations of the Newton decrements of empirical risks over the two samples.
This allows to generalize the test by using the Newton decrements in the general case. 
In addition, we have to appropriately modify the debiasing terms with which the Newton decrements are compared. 
This is done by replacing the covariance rank with the {\em  local effective rank} -- the trace of the Hessian-standardized local Fisher matrix,
which is equal to identity in the well-specified case (see the preliminary part of Section~\ref{sec:asymp} for details).
On the positive side, such modification allows to effortlessly handle misspecified models: 
the argument to bound the error probability remains similar, and we only pay for model misspecification through a (controlled) deterioration of the derived error bounds.
(This occurs due to the modified Fisher matrix not being identity in the misspecified case, see the statement of Proposition~\ref{prop:asymp} for more details.)
However, unlike the actual model ranks, the effective ranks are usually unknown, and our testing procedure must be adaptive to them to begin with. 
We address this challenge via an estimator of the effective rank, which might be interesting in a wider context.
\section{General asymptotic theory}
\label{sec:asymp}

\paragraph{Reminder of the setup.}
We now revisit the general~$M$-estimation setup introduced in Section~\ref{sec:intro}. 
Our goal is to characterize the sample complexity of confident testing in the asymptotic large-sample regime, with prediction distances vanishing as~$n_0,n_1 \to \infty$. 
Before we specify the regime and present the proposed testing procedures, let us recap the setup.
Hypothesis~$\cH_k$~($k \in \{0,1\}$) specifies the (unique) unconstrained minimizer~$\theta_k$ of the population risk~$L_k(\theta) := \E_k[\ell(\theta,z)]$ via the distribution~$\Prob_k$ of~$z$, where we use the shorthand~$\E_k[\cdot] := \E_{z \sim \Prob_k}[\cdot]$. In what follows, we require the loss function~$\ell(\cdot,z)$ to be strictly convex and twice continuously differentiable for all possible values of~$z$.
Given two samples~$Z^{(0)} \sim \Prob_0^{\otimes n_0}, Z^{(1)} \sim \Prob_1^{\otimes n_1}$ (cf.~\eqref{eq:samples}) we introduce the pair of empirical risks: 
\begin{equation}
\label{def:emp-risks}
\wh L_0(\theta) := \frac{1}{n_0}\sum_{i=1}^{n_0} \ell(\theta, z_i^{(0)}), 
\quad
\wh L_1(\theta) := \frac{1}{n_1}\sum_{i=1}^{n_1} \ell(\theta, z_i^{(1)}).
\end{equation}

In order to facilitate the reader, we shall now remind the main notions related to~$M$-estimators. 
We use a simplified notation:~$L(\theta)$, without any subscript, is the population risk corresponding to an abstract distribution~$\Prob$ (later on~$\Prob$ will be either~$\Prob_0$ or~$\Prob_1$).
We similarly omit the subscript for all matrix functionals of~$L(\theta)$ (see~\eqref{eq:matrices-population}--\eqref{def:fisher-matrix} below). 
Finally, we let~$\theta^\circ$ be a minimizer of~$L(\cdot)$ unique due to the strict convexity of~$\ell_z(\cdot)$. 
In other words, under~$\cH_k$ we have~$L(\theta) \equiv L_k(\theta)$ and~$\theta^\circ = \theta_k$. 

\paragraph{Background on~$M$-estimators.}
In the context of~$M$-estimators (see, e.g.,~\cite{ostrovskii2018finite,marteau2019beyond,spokoiny2012parametric}),
it is convenient to define the local Fisher information matrix and the population risk Hessian:
\begin{equation}
\label{eq:matrices-population}
\begin{aligned}
\G(\theta) &:= \Cov[\nabla \ell(\theta,z)], \quad
\H(\theta) := \E[\nabla^2 \ell(\theta,z)] = \nabla^2 L(\theta). 
\end{aligned}
\end{equation}
Hereinafter the loss is differentiated in the first argument.
\begin{itemize}
\item
Recall the well-known fact (see,~e.g.,~\cite{lehmann2006theory}): the identity
\begin{equation}
\label{eq:bartlett-identity}
\G(\theta^\circ) = \H(\theta^\circ),
\end{equation}
also called {Bartlett's identity}~\cite{bartlett1953approximate}, holds in the {\em well-specified case}, 
i.e., when~$\ell_z(\theta) = -\log P_{\theta}(z)$ with~$P_{\theta} \in \cP := \{P_{\theta}(\cdot),\; \theta \in \R^d\}$, 
where~$\cP$ is a family of probability densities~$P_\theta(\cdot)$ that is sufficiently regular (see, e.g.,~\cite{lehmann2006theory} for more details on regularity), and~$z$ is generated by a distribution with density from~$\cP$; this density is then necessarily~$P_{\theta^\circ}$ {\em a.s.}~by the non-negativity of the Kullback-Leibler divergence.\footnote{More generally, we can have~$\theta \in \Theta$ with~$\Theta \subseteq \R^d$ open and with non-empty interior, as we only need that~$\nabla L(\theta^\circ) = 0$.}
Note also that this was the case in the setup of Section~\ref{sec:gaussian-upper}, with~$\ell(\theta,z) = \tfrac{1}{2} (y-x^\top \theta)^2$  under the unit noise variance assumption~\eqref{eq:linear-model}, 
and in Remark~\ref{rem:scale-misspec} we shall revisit the case of linear models with unequal variances (cf.~\eqref{eq:linear-model-var}) in this context.
\item
Now, under~\eqref{eq:bartlett-identity}, the Hessian-standardized Fisher matrix~$\M(\theta)$, as given by
\begin{equation}
\label{def:fisher-matrix}
\M(\theta) := \H(\theta)^{\dagger/2} \, \G(\theta) \, \H(\theta)^{\dagger/2},
\end{equation}
satisfies~$\M(\theta^\circ) = \Id_{r^\circ}$, where~$r^\circ = \rank(\H(\theta^\circ))$.
In fact, we can redefine~$r_k := \rank(\H_k(\theta_k))$ consistently with~Sec.~\ref{sec:gaussian} since~$\H(\theta) \equiv \bSigma$ under~\eqref{eq:linear-model}, and~$\H_k^{\vphantom{2}}(\theta) \equiv \sigma_k^2 \bSigma_k^{\vphantom{2}}$ under~\eqref{eq:linear-model-var}.
In the remainder of this section, we follow this convention and refer to~$r^o$ simply as the {\em model rank.}
\item
Under model misspecification (including the case where~$\ell(\theta,z)$ does not correspond to {any} log-likelihood), we can view~$\M(\theta^\circ)$ as a distortion of~$\Id_{r^\circ}$;
then~$\Tr[\M(\theta^\circ)]$ is close to the model rank when the level of misspecification is moderate. 
More generally, we can consider the~{\em \hbox{$\ell_p$-effective} rank}~$\Tr[\M^p(\theta^\circ)]$ with~$p \in [1,\infty]$ as a generalization of the model rank.
Typically we have~$\|\M(\theta^\circ)\| \ge 1$, and~$\ell_p$-effective rank grows with~$p \ge 1$ becoming more sensitive to misspecification.
Alternatively, we might quantify the misspecification with~$(\tfrac{1}{r^\circ}\Tr[\M^p(\theta^\circ)])^{1/p}$, i.e., the ratio of the Schatten~$p$-norms of~$\M(\theta^\circ)$ and~$\Id_{r^\circ}$. 
Note that this measure also grows with~$p \ge 1$ (by H\"older's inequality), and reduces to the operator norm~$\|\M(\theta^\circ)\|$ when~$p \to \infty$. 
\end{itemize}




\subsection{Basic test via Newton decrements}
The key insight allowing to generalize test~\eqref{eq:test_1} to non-linear models is as follows.
Omitting~$k \in \{0,1\}$ for brevity, test~\eqref{eq:test_1} was based on the quantities~$\|Y - X\theta^{*}\|^2_X$ with~$Z = (X,Y) \sim \Prob^{\otimes n}$. 
Moreover, observe that in the setup of Section~\ref{sec:gaussian-upper} one has~$\nabla \wh L(\theta) = \tfrac{1}{n} X^{\top}(Y-X\theta)$ and
$
\wh \H(\theta) = \tfrac{1}{n} X^\top X.
$
As such, we can express~$\|Y - X\theta^{*}\|^2_X$ in terms of the gradient~$\nabla \wh L(\theta^*)$ and Hessian~$\wh\H(\theta^*)$ of empirical risk over the corresponding sample:
\[
\|Y - X\theta^{*}\|^2_X 
= 
\|\bPi_{X}(Y - X\theta^{*})\|^2
= 
\|(X^\top X)^{\pinv/2} X^{\top}(Y - X\theta^{*})\|^2
= 
n\| \wh \H(\theta^*)^{\dagger/2} \nabla \wh L(\theta^*) \|^2.
\]
The quantity in the right-hand side is the rescaled Newton decrement of empirical risk~$\wh L(\cdot)$ at~$\theta^*$.
This naturally leads us to replacing the quadratic terms in~\eqref{eq:test_1} with the corresponding Newton decrements~$\| \wh \H_k(\theta^*)^{\dagger/2} \nabla \wh L_k(\theta^*) \|^2$ over the two samples~($k \in \{0,1\}$). 
However, we also have to adjust the debiasing terms. 
To this end, examining the proof of Proposition~\ref{prop:upper-fixed-design}, we note that empirical rank~$\wh r_k = \rank(\wh \bSigma_k)$ in~\eqref{eq:test_1} appeared as the conditional expectation of the squared  norm of the corresponding projected residual~$\|Y^{(k)}-X^{(k)}\theta^*\|_{X^{(k)}}^2$ under the matching hypothesis~$\cH_k$. 
In the setup of Sec.~\ref{sec:gaussian},~$\rank(\wh\bSigma_k) \stackrel{a.s.}{\to} \rank(\bSigma_k)$ as~$n_k \to \infty$, and~$\rank(\bSigma_k)$ is the {total} expectation of the squared norm of the projected residual over the whole observation~$z=(x,y)$. 
Returning to the~$M$-estimation setup with~$n_0,n_1 \to \infty$, this reasoning leads us to the test given by
\begin{equation}
\label{eq:test-asymp}
\wh T = \ind \big\{
n_0 \| \wh \H_0(\theta^*)^{{\dagger}/{2}} \nabla \wh L_0(\theta^*) \|^2 - \Tr[\M_0(\theta^*)]
\ge
n_1 \| \wh \H_1(\theta^*)^{{\dagger}/{2}} \nabla \wh L_1(\theta^*) \|^2 - \Tr[\M_1(\theta^*)]
\big\}.
\tag{$\textsf{\textup{Asymp}}$}
\end{equation}
We now explain the mechanism behind this test in more detail.
First, one can easily check that 
\begin{equation}
\label{eq:trace-explicit}
\Tr[\M_k(\theta^*)] = n_k \E \big[ \big\| \H_k(\theta^*)^{{\dagger}/{2}}  \big[ \nabla \wh L_k(\theta^*) - \nabla L_k(\theta^*) \big] \big\|^2 \big],
\end{equation}
and~$\nabla L_k(\theta_k) = 0$.
Moreover,~$\wh \H_k(\theta^*)$ by the law of large numbers converges to~$\H_k(\theta^*)$ in probability as~$n_k \to \infty$.
Then Slutsky's theorem (see,~e.g.,~\cite{lehmann2004elements}) allows to replace the Newton decrement with its simpler counterpart~$n_k \| \H_k(\theta^*)^{{\dagger}/{2}} \nabla \wh L_k(\theta^*) \|^2$ in which only the gradient is random. Thus, under~$\cH_0$ (resp.,~$\cH_1$) the expectation of the left-hand side (resp., right-hand side) of the inequality in~\eqref{eq:test-asymp} vanishes, while the expectation of the right-hand side (resp., left-hand side) is positive.
More precisely, the latter two expectations are given correspondingly by~$n_1 \bar\Delta_1$ and~$n_0\bar\Delta_0$ with~$\bar\Delta_0$ and~$\bar \Delta_1$ defined as
\begin{equation}
\label{eq:distances-asymp}
\begin{aligned}
\bar\Delta_0 &:= \| \H_0(\theta_{1})^{\dagger/2} \nabla L_0(\theta_{1}) \|^2,\quad
\bar\Delta_1 &:= \| \H_1(\theta_{0})^{\dagger/2} \nabla L_1(\theta_{0}) \|^2.
\end{aligned}
\end{equation}
Essentially,~$\bar \Delta_0,\bar\Delta_1$ are two (asymmetric) measures of separation between~$\Prob_0$ and~$\Prob_1$. They generalize, in a consistent manner, the prediction distances~$\Delta_0,\Delta_1$ defined in~\eqref{def:Delta-population} for the case of linear models.

\paragraph{Statistical guarantee.}
Next we show that the error probabilities for~\eqref{eq:test-asymp} are controlled via the products~$n_0 \bar\Delta_0$ and~$n_1\bar\Delta_1$ in a similar way as for test~\eqref{eq:test_1} in the setup of Theorem~\ref{th:upper-random-design}. 
Specifically, we consider the asymptotic regime in which~$n_k \to \infty$ and~$\bar\Delta_{k} \to 0$ such that~$n_k \bar\Delta_k \to \lambda_k$ for some~$\lambda_k > 0$ ($k \in \{0,1\}$). 
In such regime, Theorem~\ref{th:upper-random-design} bounds the type I error probability for~\eqref{eq:test_1} as
\begin{equation}
\label{eq:gaussian-error-reminder}
O\left(\exp\bigg(-c \min\bigg\{\lambda_1,\frac{\lambda_1^2}{\max\{r_0,r_1\}}\bigg\}\bigg)\right).
\end{equation}
In the result presented next, this dependency is generalized for non-linear and misspecified models.
\begin{proposition}
\label{prop:asymp}
Let~$n_0, n_1 \to \infty$ and~$\theta_1-\theta_0 \to 0$ at such a rate that~$n_k \bar\Delta_k \to \lambda_k > 0$ for~$k \in \{0,1\}$.
Assume that the map~$\M_{0}(\theta)$ (resp.,~$\M_{1}(\theta)$) is continuous at~$\theta_0$ (resp.,~$\theta_1$). 
Then test~\eqref{eq:test-asymp} satisfies
\begin{equation}
\label{eq:prob-asymp}
\Prob_{\cH_0}[\wh T = 1] \to P_{\infty} \le C \exp\bigg(-c \min \bigg\{
\frac{\lambda_1}{\max\big\{\|\M_0(\theta_0)\|, \|\M_1(\theta_1)\|\big\}},  \frac{\lambda_1^2}{\max\big\{\Tr[\M_0^2(\theta_0^{\vphantom{2}})], \Tr[\M_1^2(\theta_1^{\vphantom{2}})]\big\}} 
\bigg\}\bigg).
\end{equation}
The type II error admits a similar bound with the replacements~$1 \mapsto 0$ and~$0 \mapsto 1$.
\end{proposition}
We now discuss this result. 
In a well-specified scenario, we recover~\eqref{eq:gaussian-error-reminder} since~$\M_k(\theta_k) = \Id_{r_k}$. 
Under misspecification the bound is adjusted:~$r_k$ is replaced with the corresponding~$\ell_2$-effective rank~$\Tr[\M_k^2(\theta_k^{\vphantom{2}})]$, and the linear in~$\lambda_1$ term is discounted by the largest (i.e., worst) operator norm. 

\begin{remark}
\label{rem:gradient-test}
{ Under each hypothesis, one of the two population risks vanishes at~$\theta^*$. 
Hence, one may wonder why not to use the simpler {gradient-norm test}
\begin{equation}
\label{eq:test_grad}
    \wh T = \ind\big\{ 
 \|\nabla\wh L_{0}(\theta^*)\|^2 - \E [\|\nabla \wh L_{0}(\theta_0)\|^2] \geq \| \nabla\wh L_{1}(\theta^*)\|^2 - \E [\|\nabla \wh L_{1}(\theta_1)\|^2]
 \big\}.
\tag{$\textsf{\textup{Grad}}$}
\end{equation}
However, we demonstrate, both rigorously (in Appendix~\ref{app:grad}) and in simulations (in Sec.~\ref{sec:simul}), that already for linear models,~\eqref{eq:test_1} can {\em arbitrarily} outperform~\eqref{eq:test_grad} when the design is ill-conditioned. This is not surprising, as the statistic in~\eqref{eq:test_grad} is not affine-invariant (i.e., it changes under reparametrizations~$\theta \mapsto A\theta + b$); in contrast, tests~\eqref{eq:test_1} and~\eqref{eq:test-asymp} use affine-invariant statistics.}
\end{remark}

\begin{remark}
\label{rem:scale-misspec}
It is instructive to apply the result of Proposition~\ref{prop:asymp} to the case of linear models with different noise variances (cf. Section~\ref{sec:gaussian-adaptive}).
Test~\eqref{eq:test-asymp} uses the unscaled loss~$\ell(\theta,z) = \tfrac{1}{2}(y-x^\top\theta)^2$ which corresponds to the assumption of unit noise variances~($\sigma_0^2 = \sigma_1^2 = 1$) while the actual variance is~$\sigma_0^2$ or~$\sigma_1^2$ depending on the hypothesis. 
Then~$\M_k(\theta_k) = \sigma_k^2 \Id_{r_k}^{\vphantom 2}$ and the right-hand side of~\eqref{eq:prob-asymp} reads
\[
C \exp \bigg(-\frac{c n \Delta_1}{\max\{\sigma_0^2,\sigma_1^2\}} \min\bigg\{1, \frac{n\Delta_1}{\bar r \max\{\sigma_0^2, \sigma_1^2\}}\bigg\}\bigg).
\]
Whenever~$\max\{ \sigma_0^2, \sigma_1^2\} > 1$, this bound becomes worse than~\eqref{eq:upper-random-design-adapt} for test~\eqref{eq:test_1-adapt} (cf.~Corollary~\ref{cor:upper-random-design-adapt}):~$\Delta_1$ is discounted by the largest of two variances.
This is because test~\eqref{eq:test_1-adapt} circumvents the issue of model misspecification: by estimating~$\sigma_0^2, \sigma_1^2$ in~\eqref{eq:linear-model-var} it manages to mimic the oracle test~\eqref{eq:test_1-var},
but running~\eqref{eq:test_1-var} is equivalent to running~\eqref{eq:test_1} run on the rescaled samples~$\tfrac{1}{\sigma_0}(X^{(0)},Y^{(0)})$ and~$\tfrac{1}{\sigma_1}(X^{(1)}, Y^{(1)})$ that both have unit noise variance, i.e., to a well-specified scenario.
Meanwhile, test~\eqref{eq:test_1} suffers from the effect of misspecification unless~$\max\{\sigma_0^2,\sigma_1^2\} \le 1$ (i.e., unless the actual distribution is less noisy than the model assumes). 
This comparison illuminates the virtue of adaptation in Section~\ref{sec:gaussian-adaptive}.
\end{remark}

\begin{remark}
\label{rem:traces-wellspec}
To be applied, test~\eqref{eq:test-asymp} requires the knowledge of the traces~$\Tr[\M_0(\theta^*)],\Tr[\M_1(\theta^*)]$.
In the well-specified case, these traces can be inferred from observations, and~\eqref{eq:test-asymp} can be applied.
Indeed, in this case, under~$\cH_0$ (w.l.o.g.) one has
\[
\begin{aligned}
\Tr[\M_0(\theta^*)] &= \Tr[\M_0(\theta_0)] = \Tr[\Id_{r_0}] = r_0, \\
\Tr[\M_1(\theta^*)] &= \Tr[\M_1(\theta_0)] = r_1 + o(1),
\end{aligned}
\]
where~$r_k = \rank[\H_k(\theta_k)]$ is the corresponding model rank, and~$\Tr[\M_1(\theta_0) - \M_1(\theta_1)] = o(1)$ vanishes as~$\theta_1 - \theta_0 \to 0$ by the continuity of~$\M_1(\cdot)$. 
On the other hand, test~\eqref{eq:test-asymp} is robust to additive perturbation of the test statistic by~o(1): inspecting the proof of Proposition~\ref{prop:asymp} we see that guarantee~\eqref{eq:prob-asymp} remains valid if the statistic in~\eqref{eq:test-asymp} is perturbed by an additive~$o(1)$ term.
Thus, the two traces in~\eqref{eq:test-asymp} can be safely replaced with~$r_0,r_1$ while preserving~\eqref{eq:prob-asymp}.
Finally, we can replace~$r_0, r_1$ with the quantities~$\rank(\wh\H_0(\theta^*)), \rank(\wh\H_1(\theta^*))$ which are observable in our regime of interest here.
Indeed, under~$\cH_0$ (w.l.o.g.), we have that~$\rank(\wh \H_0(\theta^*)) = \rank(\wh \H_0(\theta_0)) \stackrel{a.s.}{\to} \rank(\H_0(\theta_0)) = r_0$ and
\[
\begin{aligned}
\rank(\wh \H_1(\theta^*)) = \rank(\wh \H_1(\theta_0)) \stackrel{a.s.}{\to} \rank(\H_1(\theta_0));
\end{aligned}
\]
meanwhile~$\rank(\H_1(\theta_0)) - r_1 \to 0$ as~$\theta_1 - \theta_0 \to 0$ under the mild assumption that~$\H_1(\cdot)$ is continuous.
\end{remark}

In the misspecified case, the argument in Remark~\ref{rem:traces-wellspec} is no longer valid; 
the terms~$\Tr[\M_k(\theta^*)]$ ($k \in \{0,1\}$) in~\eqref{eq:test-asymp} cannot be simply replaced with the ranks of empirical risk Hessians, and~\eqref{eq:test-asymp} cannot be applied.
Next we present a technique for estimating these traces, which leads to an adaptive counterpart of test~\eqref{eq:test-asymp} with essentially the same guarantee as for~\eqref{eq:test-asymp}. 

\subsection{Adaptive test}
\label{sec:asymp-adapt}



We now construct a counterpart of test~\eqref{eq:test-asymp} adaptive to the terms~$\cT_k := \Tr[\M_0(\theta^*)]$, $k~\in \{0,1\}$.
The new test has the form
\begin{equation}
\label{eq:test-asymp-adapt}
\begin{aligned}
\wh T = \ind \big\{
n_0 \| \wh \H_0(\theta^*)^{{\dagger}/{2}} \nabla \wh L_0(\theta^*) \|^2 - \wh\cT_0
\ge
n_1 \| \wh \H_1(\theta^*)^{{\dagger}/{2}} \nabla \wh L_1(\theta^*) \|^2 - \wh\cT_1 \big\},
\end{aligned}
\tag{$\textsf{\textup{Asymp$^{+}$}}$}
\end{equation}
where~$\wh\cT_0,\wh\cT_1$ estimate~$\cT_0, \cT_1$ as follows (cf.~\eqref{eq:trace-explicit}):
\begin{equation}
\label{eq:trace-estimate}
\wh\cT_k 
:= 
\frac{n_k}{2} 
\big\| \wh\H_k(\theta^*)^{{\dagger}/{2}}  \big[ \nabla \wh L_k(\theta^*) - \wt \nabla \wh L_k(\theta^*) \big] \big\|^2 
	\quad (k \in \{0,1\}).
\end{equation}
Here~$\wt \nabla \wh L_0(\theta^*)$ is an independent copy of~$\nabla \wh L_0(\theta^*)$ computed using an independent copy~$\wt Z^{(k)}$ of~$Z^{(k)}$ (note that such a copy can always be obtained by sample splitting).
For the resulting test~\eqref{eq:test-asymp-adapt} we have the following result. 
\begin{theorem}
\label{th:asymp-adapt}
Under the premise of Proposition~\ref{prop:asymp}, test~\eqref{eq:test-asymp-adapt} with estimates~\eqref{eq:trace-estimate} satisfies the same bound~\eqref{eq:prob-asymp} as test~\eqref{eq:test-asymp} up to constant factors.
\end{theorem}
We now sketch the proof of this result. Defining~$\cT_k = \Tr[\M_k(\theta^*)]$ for~$k \in \{0,1\}$, let
\[
\wh S := n_0 \| \wh \H_0(\theta^*)^{{\dagger}/{2}} \nabla \wh L_0(\theta^*) \|^2 - \cT_0 + \cT_1
-
n_1 \| \wh \H_1(\theta^*)^{{\dagger}/{2}} \nabla \wh L_1(\theta^*) \|^2
\] 
be the statistic such that~$\wh T = \ind\{\wh S \ge 0\}$ for test~\eqref{eq:test-asymp}. 
Recall that
$
\E_{\cH_0}[\wh S] = -n_1 \bar \Delta_1 \to  -\lambda_1$
and~$
\E_{\cH_1}[\wh S] = n_0 \bar \Delta_0 \to \lambda_0$
in the asymptotic regime of interest (i.e., the one from the premise of Proposition~\ref{prop:asymp}).
Thus, in this regime the limiting probability of type I error for~\eqref{eq:test-asymp} is upper-bounded by~$\Prob_{\cH_0} [\wh W_1 \ge \lambda_1],$ where~$\wh W_1 := \wh S + \lambda_1$ satisfies~$\E_{\cH_0}[\wh W_1] \to 0$. 
Working under~$\cH_0$ (w.l.o.g.), in order to prove Proposition~\ref{prop:asymp} we proceed by decomposing~$\wh W_1$ as follows (see Appendix~\ref{app:upper-general}):
\begin{equation}
\label{eq:centered-asymp-statistic} 
\begin{aligned}
\wh W_1 
= \left( n_0 \| \wh \H_0(\theta_0)^{{\dagger}/{2}} \nabla \wh L_0(\theta_0) \|^2 - \cT_0 \right)
	&- \left(n_1 \| \wh \H_1(\theta_0)^{{\dagger}/{2}} [\nabla \wh L_1(\theta_0) - \nabla L_1(\theta_0)]\|^2 - \cT_1\right) 
	  + \zeta,
\end{aligned}
\end{equation}
where the first two terms have limiting centered chi-squared type distributions, namely those of the centered squares of~$\cN(0,\M_0(\theta_0))$ and~$\cN(0,\M_1(\theta_0))$, and~$\zeta$ corresponds to the cross-term with limiting centered normal distribution. 
Then~\eqref{eq:prob-asymp} followed by bounding from above~$\Prob_{\cH_0}[\wh W_1 \ge \lambda_1]$ 
with the sum of probabilities for each of the three terms to exceed~$\lambda_1/3$, and by combining the corresponding deviation bounds. 
In the present situation, the statistic associated to~\eqref{eq:test-asymp-adapt} reads~$\wh S + \cT_0 - \wh \cT_0 + \wh \cT_1 - \cT_1,$
thus the limiting type I error probability for~\eqref{eq:test-asymp-adapt} is at most
\[
\Prob_{\cH_0}[\wh W_1 \ge \lambda_1/3] + \Prob_{\cH_0} [ -(\wh \cT_0 - \cT_0) \ge \lambda_1/3 ] + \Prob_{\cH_0} [\wh \cT_1 - \cT_1 \ge \lambda_1 / 3].
\]
Now, the first term admits the same bound as the one we derived for~$\Prob_{\cH_0}[\wh W_1 \ge \lambda_1]$ (since the deviation bounds on which we relied are preserved when scaling the deviation by a constant factor). 
On the other hand, the terms~$\wh \cT_0 - \cT_0$ and~$\wh \cT_1 - \cT_1$ have the same limiting distributions as the first two terms in~\eqref{eq:centered-asymp-statistic} (cf.~\eqref{eq:trace-estimate}), hence they can be controlled via the same deviation bounds as before.
\qed\\
%

Next we shall focus on small-sample regime~\eqref{eq:small-sample-regime-intro} and 
demonstrate that a slight modification of test~\eqref{eq:test-asymp} allows to handle generalized linear models under~\eqref{eq:small-sample-regime-intro} and weak moment assumptions.



\section{Generalized linear models}
\label{sec:glm}

The next part of our theory focuses on generalized linear models (GLMs) in the small-sample regime
\begin{equation}
\label{eq:small-sample-regime-glm}
n_k \le r_k \quad (k \in \{0,1\}).
\end{equation}
This assumption is not too strong: in Section~\ref{sec:gaussian-upper} we have seen that the sample complexity of confident testing typically admits an even smaller bound, unless in the case of very poor separation (cf.~Corollary~\ref{cor:complexity-gaussian} and the subsequent discussion), and we aim at proving a similar upper complexity bound for generalized linear models (so the reasoning would still apply). 
On the other hand, working in this regime allows us to diagonalize the Newton decrements as in this case there is no explicit projection. 
This, in turn, simplifies statistical analysis of the test (and also its computation). 
We shall now specify the setup and introduce auxiliary quantities to facilitate the presentation of results. 
 
\paragraph{Background on GLMs.}
Our observation is~$z = (x,y)$ with~$x \in \R^d$ and~$y \in \cY \subseteq \R$ where~$\cY$ is a label space (e.g.,~$\R$ or~$\{-1,1\}$). 
Although we expect that our analysis generalizes to general single-index models, here we focus on the case of canonical exponential family for~$\cY$ (see, e.g.,~\cite{mccullagh1989generalized}).
In other words, we assume that the loss is of the form~$\ell(\theta,z) = \phi(x^\top \theta, y)$ with~$\phi(\eta,y)$ given by
\begin{equation}
\label{eq:glm-exp-family}
\phi(\eta, y) = -\eta y + a(\eta) - b(y);
\end{equation}
here~$a(\eta) = \log[\int_{\cY} \exp({\eta y+b(y)}) dy]$ is the {\em cumulant} for the canonical distribution that corresponds to the (conditional) density~$P(y|\eta) = \exp({-\phi(\eta,y)})$. 
One can easily verify, directly or using~\cite[Sec.~2.1]{ostrovskii2018finite}, that~$a(\eta)$ is twice differentiable and strictly convex (unless~$P(y|\eta)$ is {\em a.s.}~deterministic).
For example, one has~$a(\eta) = \tfrac{1}{2}\eta^2$ in the setup of~Section~\ref{sec:gaussian} (linear models with~$\sigma_0^2 = \sigma_1^2 = 1$) and~$a(\eta) = \log(e^{\eta} + e^{-\eta})$ for the logistic loss with~$\cY = \{\pm 1\}$.
Now, defining the mapping~$\theta \mapsto \eta$ by
\[
\eta = \eta_x(\theta) := x^\top \theta,
\]
the gradient and Hessian of~$\ell(\cdot,z)$ can be expressed as
\begin{equation}
\label{eq:glm-loss-derivatives}
\begin{aligned}
\nabla \ell(\theta,z) &= (a'(\eta) - y) x, \quad
\nabla^2 \ell(\theta,z) = a''(\eta) xx^\top.
\end{aligned}
\end{equation}
Note that~$\nabla^2 \ell(\theta,z)$ does not depend on~$y$ under the canonical exponential family assumption~\eqref{eq:glm-exp-family}.
This property simplifies our analysis due to more straightforward conditioning. 
However, we anticipate that our results can be extended to general single-index models -- for example, by following~\cite{marteau2019beyond}.



\paragraph{Auxiliary quantities.}
We shall now define some functionals of~$\Prob_{0}, \Prob_1$ related to the exponential family structures~\eqref{eq:glm-exp-family}; this will facilitate the presentation of our results and subsequent discussions.
When defining these functionals, we shall use the same conventions as in Section~\ref{sec:asymp}, i.e., omit the subscript~$k \in \{0,1\}$ and let~$\theta^\circ$ be the population risk minimizer corresponding to an abstract distribution~$\Prob$ of~$z$.
Later on, we shall add the subscript~$k \in \{0,1\}$ and let~$\theta^\circ = \theta_k$ under~$\cH_k$.
\begin{itemize}
\item
Conditional \textit{relative label variance}~$\nu(\theta|x)$ and {\em normalized label kurtosis}~$\kappa(\theta|x)$:
\begin{equation}
\label{eq:conditional-moments}
\begin{aligned}
\nu(\theta|x) &:= \frac{\E[(y - \E[y|x])^2|x]}{a''(\eta_x(\theta))}, \quad
\kappa(\theta|x) := \frac{\E[(y - \E[y|x])^4|x]}{a''(\eta_x(\theta))^2}, \\
\end{aligned}
\end{equation}
and their marginals over~$x$:
$\nu(\theta) := \E[\nu(\theta|x)]$ and~$\kappa(\theta) := \E[\kappa(\theta|x)]$.
We define a pair with subscript~$k \in \{0,1\}$ for each of the above quantities, with expectations~$\E_k[\cdot] := \E_{(x,y) \sim \Prob_k}[\cdot]$. 

\item 
Conditional {\em squared misspecification bias}~$\Bias(\theta|x)$:
\begin{equation}
\label{eq:misspec-bias}
\Bias(\theta|x) := \frac{(a'(\eta_{x}(\theta^\circ) - \E[y|x])^2}{a''(\eta_{x}(\theta))}
\end{equation}
and its marginal~$\Bias(\theta) := \E[\Bias(\theta|x)]$ over~$x$. 
As before, we define~$\Bias_k(\theta|x)$ and~$\Bias_k(\theta)$ for~$k \in \{0,1\}$.
\end{itemize}

The names of the quantities defined above stem from the well-known fact (see, e.g.,~\cite[Sec.~2.1]{ostrovskii2018finite}) that, under~\eqref{eq:glm-exp-family}, the derivatives~$a^{(p)}(\eta)$ of the cumulant ($p \in \mathds{N}$) are equal to the central moments of~$y$ according to the canonical distribution.
Thus, in the well-specified case 
one has~$\nu(\theta^\circ|x) = 1$ and $\beta(\theta|x) = 0$ at any~$\theta \in \R^d$ {\em a.s.}~over~$x$. 
Thus, the conditional bias~$\beta(\theta^\circ|x)$ and the excess conditional variance~$\nu(\theta^\circ|x)-1$ are two local measures of model misspecification.\footnote{The two are of different nature:~$\nu(\theta^*|x)$ is sensitive to ``scale of noise'' whereas~$\beta(\theta^*|x)$ to the ``location'' of the optimal parameter. 
For example, in a linear model with misspecified variance of the noise,~$\beta(\theta^*|x) = 0$ but~$\nu(\theta^*|x) \ne 1$.} 
Similarly,~$\kappa(\theta|x)$ is the ratio of the actual (conditional) fourth central moment of~$y$ and its squared variance according to the canonical distribution with the substitution~$\eta = x^\top \theta$. 
Note that~$\kappa(\theta|x) \ge \nu^2(\theta|x)$, and~$\kappa(\theta^\circ|x)$ is the actual conditional kurtosis of~$y$ in the well-specified case, i.e., according to the density~$e^{-\phi(y|x^\top \theta^\circ)}$.
\\ 

Below we let~$r_k := \rank(\bSigma_k) \equiv \rank(\H_k(\theta))$; due to~$a''(\eta) > 0$ the identity holds  whenever the marginal distribution of~$x$ under~$\Prob_k$ has density on its support, which we assume from now on.


\paragraph{Basic test construction.}
The test statistic we are about to present is analogous to the one used in test~\eqref{eq:test-asymp} in Section~\ref{sec:asymp}. 
However, now we focus on the {\em small-sample regime} (cf.~\eqref{eq:small-sample-regime-glm})
and adjust the statistic accordingly.
First, we replace the effective ranks~$\Tr[\M_k(\theta^*)]$ in~\eqref{eq:test-asymp} with the scaled variances~$n_k \nu_k(\theta^*)$ that are typically smaller. 
(For example, in the setup of Section~\ref{sec:gaussian-upper} 
this amounts to replacing~$r_k$ with~$\min\{n_k,r_k\}$, cf.~\eqref{eq:test_1}.
Second, we express the Newton decrements explicitly without matrix inversion. 
To this end, we first rewrite~\eqref{eq:glm-loss-derivatives} as
\[x
\nabla \ell_{z_i}(\theta) = \rho_i(\theta) x_i(\theta), \quad
\nabla^2 \ell_{z_i}(\theta) = x_i(\theta) x_i(\theta)^\top,
\] 
where we define
\[
\rho_i(\theta) := \frac{a'(\eta_{x_i}(\theta)) - y_i}{\sqrt{a''(\eta_{x_i}(\theta))}}, 
\quad
x_i(\theta) := \sqrt{a''(\eta_{x_i}(\theta))} x_i.
\]
Here~$\rho_i(\theta)$ can be understood as local residuals and~$x_i(\theta)$ as local predictors, in both cases rescaled by the standard deviation according to the local canonical distribution.
In this notation, the Newton decrement corresponding to~$\Prob_k$ reads
\[
\| \wh \H_k(\theta)^{{\pinv}/{2}} \nabla \wh L_k(\theta) \|^2 =  \|[\rho_1^{(k)}(\theta); ...; \rho_{n_k}^{(k)}(\theta)]\|_{X^{(k)}(\theta)}^2
\]
where on the right-hand side we measure the norm of the local residual vector (in~$\R^{n_k}$) projected onto the column space of the local predictor matrix~$X^{(k)}(\theta) := [x_1^{(k)}(\theta); ...; x_{n_k}^{(k)}(\theta)]\in \R^{n_k \times d}$. 
Observe that, under~\eqref{eq:small-sample-regime-glm},~$X^{(k)}(\theta)$ has full column rank, thus the corresponding Newton decrement reduces to the sum of squared local residuals. 
This reasoning leads to the following test:
\begin{equation}
\label{eq:test-glm}
\wh T = \ind\left\{
\sum_{i=1}^{n_0} 
	\big[ \rho_i^{(0)}(\theta^*) \big]^2 - n_0 \nu_0(\theta^*) 
\geq \sum_{i=1}^{n_1} 
	\big[ \rho_i^{(1)}(\theta^*) \big]^2 - n_1 \nu_1(\theta^*) \right\},
\tag{$\textsf{\textup{GLM}}$}
\end{equation}
However, this test requires the knowledge of the local relative variances~$\nu_0(\theta^*),\nu_1(\theta^*)$ at least one of which 
can hardly be assumed to be known even in the well-specified case, as it uses the distribution~$\bar\Prob$ (unless in the case of well-specified linear models where~$\nu(\theta) \equiv 1$).
Thus, test~\eqref{eq:test-glm} is not practical.
It is nonetheless instructive to study its statistical properties first. 
Later on we shall construct its adaptive counterpart that performs essentially as good as if~$\nu_0(\theta^*),\nu_1(\theta^*)$ were known. (This shall, however, require Assumption~\ref{ass:resampling} as in the case of small-sample variance estimator~\eqref{eq:resampling-simple} in Section~\ref{sec:gaussian-adaptive}.)


\paragraph{Modified notion of separation.}
In order to characterize the sample complexity for test~\eqref{eq:test-glm}, 
we shall first extend the notion of model separation to GLMs.
Namely, we replace~\eqref{def:Delta-empirical} and~\eqref{def:Delta-population} with
\begin{equation}
\label{def:Delta-glm}
\wh\Delta_k(x) := \frac{[a'(\eta_x(\theta_1)) - a'(\eta_x(\theta_{0}))]^2}{a''(\eta_x(\theta_{1-k}))}
\;\; \text{and} \;\;
\Delta_k := \E_{(x,y) \sim \Prob_k}[\wh\Delta_k(x)]
\quad (k \in \{0,1\}).
\end{equation}
For well-specified linear models (cf.~Section~\ref{sec:gaussian}), we have that~$a'(\eta_x(\theta)) = \eta_x(\theta) = x^\top\theta$ and~$a''(\eta) \equiv 1$; 
the resulting definition of~$\Delta_k$ recovers both the squared prediction distances~\eqref{def:Delta-population} used in Section~\ref{sec:gaussian} and the population Newton decrements~$\bar\Delta_0,\bar\Delta_1$ (cf.~\eqref{eq:distances-asymp}) used in Section~\ref{sec:asymp}.
Otherwise, the new definition extends~\eqref{def:Delta-population} in a way slightly different from~\eqref{eq:distances-asymp}: 
whereas in~\eqref{eq:distances-asymp} we average marginally over~$z$, but separately in~$\nabla L_k(\theta^*)$ and~$\H_k(\theta)$,
in~\eqref{def:Delta-glm} we condition on~$x$ {\em before} averaging. 
In particular,~$\Delta_0$ (cf.~\eqref{def:Delta-glm}) can also be expressed as
\[
\Delta_0 = \E_0 \bigg[\frac{1}{n_0} \sum_{i=1}^{n_0} \wh\Delta_0(x_i)\bigg] 
= \E\left[ \| \wh \H_0(\theta_{1})^{\dagger/2} \nabla \wh L_0(\theta_{1}) \|^2 \right].
\]
In fact, the last expression generalizes the second definition in~\eqref{def:Delta-glm} to arbitrary~$n_k$. 
This allows to connect the two notions: clearly, $\Delta_0 \to \bar \Delta_0$ as~$n_0 \to \infty$.


\paragraph{Moment assumptions.}
The result presented next holds under weak moment assumptions about~$\Prob_0,\Prob_1$: boundedness of the normalized kurtosis~$\kappa_0(\theta_0),\kappa_1(\theta_1),\kappa_0(\theta_1),\kappa_1(\theta_0)$ 
and the moments
\[
\E_k^{\vphantom 2} [\wh \Delta_k^2] = \E_{(x,y) \sim \Prob_k}[\wh \Delta_k^2(x)] \quad (k \in \{0,1\}).
\]
Upon inspection, in the case of a linear model~$\E_{k}^{\vphantom 2}[\wh \Delta_k^2]$ reduces to the expected fourth power of the design marginal along the direction~$\theta_1 - \theta_0$. 
Thus, this is as well a fourth-order moment assumption.
Meanwhile, we only allow for a {\em moderate} level of misspecification.
Essentially, we do this by not allowing the misspecification bias (cf.~\eqref{eq:misspec-bias}) be too big in terms of separations~$\Delta_0, \Delta_1$. 
More precisely, we assume that, for~$k \in \{0,1\}$ and~$\sfC := 64 \cdot 10^4$, 
\begin{subequations}
\label{eq:small-bias}
\begin{align}
\label{eq:small-bias-i}
\E_k[\Bias_k(\theta_{0}|x)^2] &\le {\Delta_1^2}/{\sfC}, \\
\label{eq:small-bias-ii}
\E_k[\Bias_k(\theta_1|x)^2] &\le {\Delta_0^2}/{\sfC}.
\end{align}
\end{subequations}
Constant~$\sfC$ also appears in the statistical guarantees presented next.
Technically, it stems from a few sequential applications of the Paley-Zygmund inequality~\cite{paley1932note}, and we expect that~$\sfC$ can be significantly reduced with some care.

\paragraph{Fixed-confidence bound.}


The result we present next extends the sample complexity bound established in~Corollary~\ref{cor:complexity-gaussian} to the case of GLMs in the small-sample regime.  
For simplicity of discussion, we assume~$n_1 = n_2 [=n]$ and introduce~$\Delta = \min\{\Delta_0,\Delta_1\}$ in the rest of this section.
Recall that, by~Corollary~\ref{cor:complexity-gaussian}, for test~\eqref{eq:test_1} the sample complexity is 
\[
n = O\bigg(\min\bigg\{\frac{1}{\Delta^2}, \frac{\sqrt{r}}{\Delta} \bigg\}\bigg).
\]
Thus, requiring this sample complexity to fall into the range~$n  \le \min\{r_1,r_2\}$ in accordance with~\eqref{eq:small-sample-regime-glm} 
would imply that~$\Delta = \Omega(1/\sqrt{\max\{r_0,r_1\}})$, cf.~\eqref{eq:high-dim-range}, leaving us with a rank-independent sample complexity bound~$n = O(1/\Delta^2)$. 
The result presented next generalizes this conclusion to GLMs. 

\begin{theorem}
\label{th:glm-fixed-prob}
Under~\eqref{eq:small-sample-regime-glm} and~\eqref{eq:small-bias-i}, the type~I error probability for test~\eqref{eq:test-glm} is at most~$2/5$ when
\begin{equation}
\label{eq:condition_n}
n 
\ge 
\frac{2\sfC 
\max\big\{
\kappa_0^{\vphantom 2}(\theta_0^{\vphantom 2}), \,
\kappa_1^{\vphantom 2}(\theta_0^{\vphantom 2}), \,
4\E_{1}^{\vphantom 2}[\wh \Delta_1^2] \big\}}{\Delta_1^2}.
\end{equation} 
The type II error probability is at most~$2/5$ under~\eqref{eq:small-bias-ii} and the condition complementary to~\eqref{eq:condition_n}.
\end{theorem}
The proof of Theorem~\ref{th:glm-fixed-prob} is technical and we defer it to Appendix~\ref{app:proofs-glm}. 
Combining~\eqref{eq:condition_n} with~\eqref{eq:small-sample-regime-glm} we conclude that the sample complexity of testing by~\eqref{eq:test-glm} is~$O(1/\Delta^2)$ whenever~$\Delta = \min\{\Delta_0, \Delta_1\}$ satisfies~$\Delta = \Omega(1/\sqrt{\min\{r_0, r_1\}})$. 
Here the degradation of the minimal separation restriction from~$\Delta = \Omega(1/\sqrt{\max\{r_0, r_1\}})$ (cf.~\eqref{eq:high-dim-range}) results from the small-sample regime assumption in~\eqref{eq:small-sample-regime-glm}.


\paragraph{Boosting confidence via majority vote.}
Next we extend the construction to handle arbitrary sample size and reach arbitrary confidence.
W.l.o.g., we split each sample into $b$ non-overlapping blocks with uniform size~$m$ to be specified later.
For each block, we run the test~$\wh T_j$ defined in~\eqref{eq:test-glm}. 
Then we aggregate the binary outcomes of these tests via the majority-vote rule, i.e., run the test
\begin{equation}
\label{eq:test_voting}
\wh T = \ind \left\{ \textstyle \sum_{j \in [b]} \wh T_j \geq {b}/{2} \right\}.
\end{equation}
By Theorem~\ref{th:glm-fixed-prob}, taking 
\begin{equation}
\label{eq:choose_k}
m
= \left\lceil 
\frac{2\sfC 
\max\big\{
\kappa_0^{\vphantom 2}(\theta_0^{\vphantom 2}), \,
\kappa_1^{\vphantom 2}(\theta_0^{\vphantom 2}), \,
4\E_1^{\vphantom 2}[\wh \Delta_1^2]
\big\}}{\Delta_1^2}
\right\rceil.
\end{equation}
we guarantee that each test~$\wh T_j$ makes a mistake with probability at most~$2/5$. 
Applying Hoeffding's inequality to~$b$ independent Bernoulli random variables we decrease the probability of error exponentially fast in~$b$ and arrive at the following result (see Appendix~\ref{app:proofs-glm} for the detailed proof).
\begin{corollary}
\label{cor:voting}

Under the premise of Theorem \ref{th:glm-fixed-prob}, the test in~\eqref{eq:test_voting} with~$m = \lceil n/b \rceil$ given by~\eqref{eq:choose_k}~satisfies
\[
\Prob_{\cH_0}[\wh T = 1] 
\le C\exp\left(
-\frac{cn\Delta_1^2}{\max\big\{\kappa_0(\theta_0),\, \kappa_1(\theta_0),\, 4\E[\wh \Delta_1^2] \big\}} \right)
\] 
if~$m \le \min\{r_0, r_1\}$.
The matching bound for the type II error holds under the complementary conditions.
\end{corollary}

\paragraph{Adaptive test.}
Next we construct a counterpart of test~\eqref{eq:test-glm} adaptive to relative label variances~$\nu_0(\theta^*),\nu_1(\theta^*)$. 
In the same way as when using the estimate~\eqref{eq:resampling-simple} of the noise variance in Section~\ref{sec:gaussian-adaptive}, the adaptive test relies on resampling the labels conditionally on design:
it is given by
\begin{equation}
\label{eq:test-glm-adapt}
\wh T = \ind\left\{
\sum_{i=1}^{n_0} \big[ \rho_i^{(0)}(\theta^*) \big]^2 
				- n_0 \wh\nu_0(\theta^*|X^{(0)}) 
\geq 
\sum_{i=1}^{n_1} \big[ \rho_i^{(1)}(\theta^*)  \big]^2
- n_1 \wh\nu_1(\theta^*| X^{(1)}) \right\}
\tag{$\textsf{\textup{GLM$^{+}$}}$}
\end{equation}
where~$\nu_0(\theta^*)$ and~$\nu_1(\theta^*)$ are replaced with their estimates -- analogues of the variance estimates~\eqref{eq:resampling-simple}:
\[
\begin{aligned}
\wh\nu_k(\theta^*|X^{(k)}) = \frac{1}{2n_k} \sum_{i=1}^{n_k} \frac{\big(y_i^{(k)} - \wt y_i^{(k)}\big)^2}{a''(\eta_{x_i{}^{(k)}}(\theta^*))}.
%
\end{aligned}
\]
As in the case of the variance estimate~\eqref{eq:resampling-simple} for the adaptive test in Section~\ref{sec:gaussian-adaptive}, here we use Assumption~\ref{ass:resampling} to sample~
~$\wt Y^{(k)} = [\wt y_1^{(k)}; ...; \wt y_{n_k}^{(k)}]$ which is i.i.d.~with~$Y^{(k)}$~conditionally on~$X^{(k)}$. 
Now observe that, using~\eqref{eq:glm-loss-derivatives} and recalling the small-sample condition~\eqref{eq:small-sample-regime-glm} we can express~$\wh\nu_k(\theta^*|X^{(k)})$ as
\[
\wh \nu_{k}(\theta^*|X^{(k)}) = \frac{1}{2} \|\wh \H_{k}(\theta^*)^{\pinv/2}[\nabla \wh L_k(\theta^*|X^{(k)}) - \wt \nabla \wh L_k(\theta^*|X^{(k)})]\|^2  \quad (k \in \{0,1\})
\]
where~$\wt \nabla \wh L_k(\theta^*|X)$ is the gradient of empirical risk over the sample $(X^{(k)}, \wt Y^{(k)})$ with resampled labels.
(Note that this is different from~\eqref{eq:trace-estimate} where we resampled the whole sample.)
Our next result shows that test~\eqref{eq:test-glm-adapt} indeed manages to mimic the idealized test~\eqref{eq:test-glm}.

\begin{proposition}
\label{prop:glm-adapt}
Granted Assumption~\ref{ass:resampling}, the claim of Theorem~\ref{th:glm-fixed-prob} remains valid if test~\eqref{eq:test-glm} is replaced with test~\eqref{eq:test-glm-adapt}.
\end{proposition}

This result is proved in appendix. 
The high-level idea is to recycle the proof of Theorem~\ref{th:glm-fixed-prob} similarly to how Theorem~\ref{th:asymp-adapt} was reduced to Proposition~\ref{prop:asymp} in the case of test~\eqref{eq:test-asymp-adapt}. 
More precisely, we write the statistic whose sign is examined in~\eqref{eq:test-glm-adapt} as~$\wh S + \wh V_1 - \wh V_0$ where~$\wh S$ is the corresponding statistic in~\eqref{eq:test-glm} and~$\wh V_k := n_k [\wh\nu_k(\theta^*|X^{(k)}) - \nu_k(\theta^*)]$ are zero-mean fluctuations due to variance estimates.
The new fluctuation terms turn out to admit similar tail bounds to some of the terms already arising in the decomposition for~$\wh S$ in Theorem~\ref{th:asymp-adapt}. 
Now, under~$\cH_0$ (w.l.o.g.) we have that~$\E_{\cH_0}[\wh S] = -n_1\Delta_1 + R$, where the term~$R$ appears due to model misspecification and can be bounded as~$|R| \le c n_1 \Delta_1$ for some~$c \ll 1$ using~\eqref{eq:small-bias-i}. 
When combined together, these two facts allow to bound the type I error for~\eqref{eq:test-glm-adapt} by
\[
\Prob_{\cH_0} \left[\wh S - \E_{\cH_0}[\wh S] \ge \frac{(1-c)n_1 \Delta_1}{3} \right] + \Prob_{\cH_0} \left[-\wh V_0 \ge \frac{(1-c)n_1 \Delta_1}{3}  \right] + \Prob_{\cH_0} \left[\wh V_1 \ge \frac{(1-c)n_1 \Delta_1}{3} \right].
\]
Then we recycle the proof of Theorem~\ref{th:glm-fixed-prob} to bound the first term and similarly bound the new terms.


\section{Numerical experiments}
\label{sec:simul}

\begin{figure}[t!]
\center
\begin{minipage}{0.32\textwidth}
\centering
\includegraphics[width=1\textwidth,clip=true,angle=0]{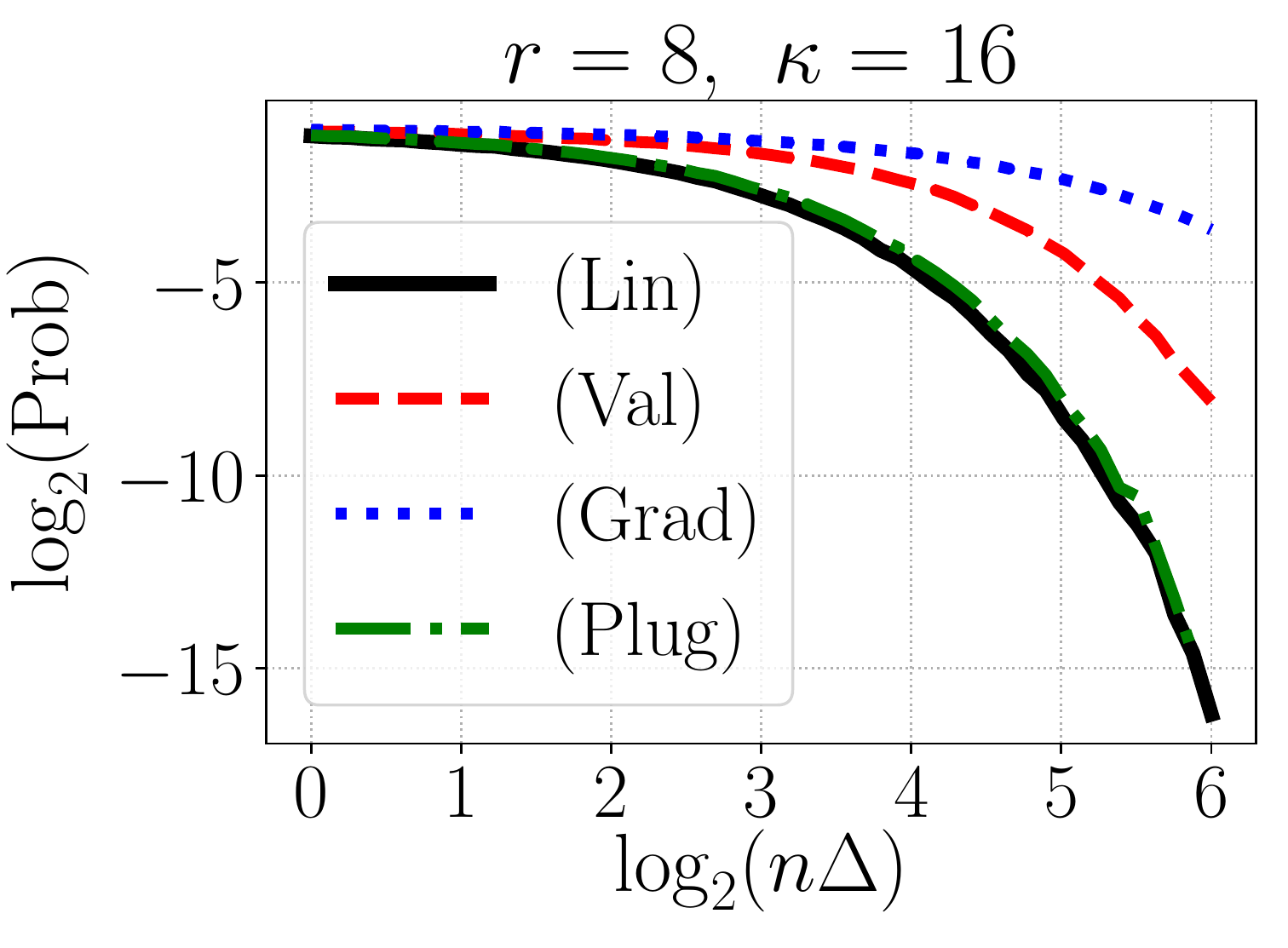}
\end{minipage}
\begin{minipage}{0.32\textwidth}
\centering
\includegraphics[width=1.04\textwidth,clip=true,angle=0]{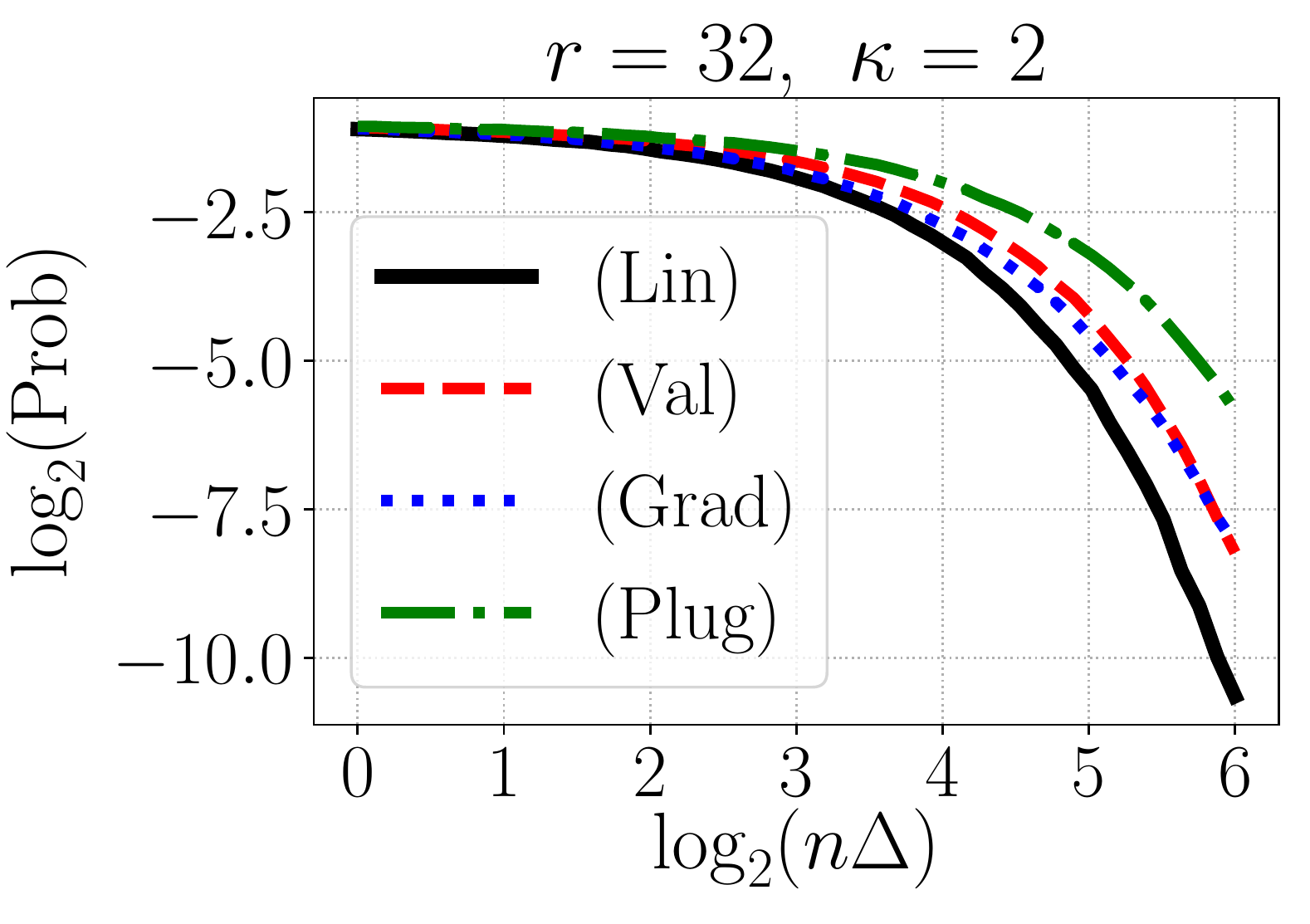}
\includegraphics[width=1.04\textwidth,clip=true,angle=0]{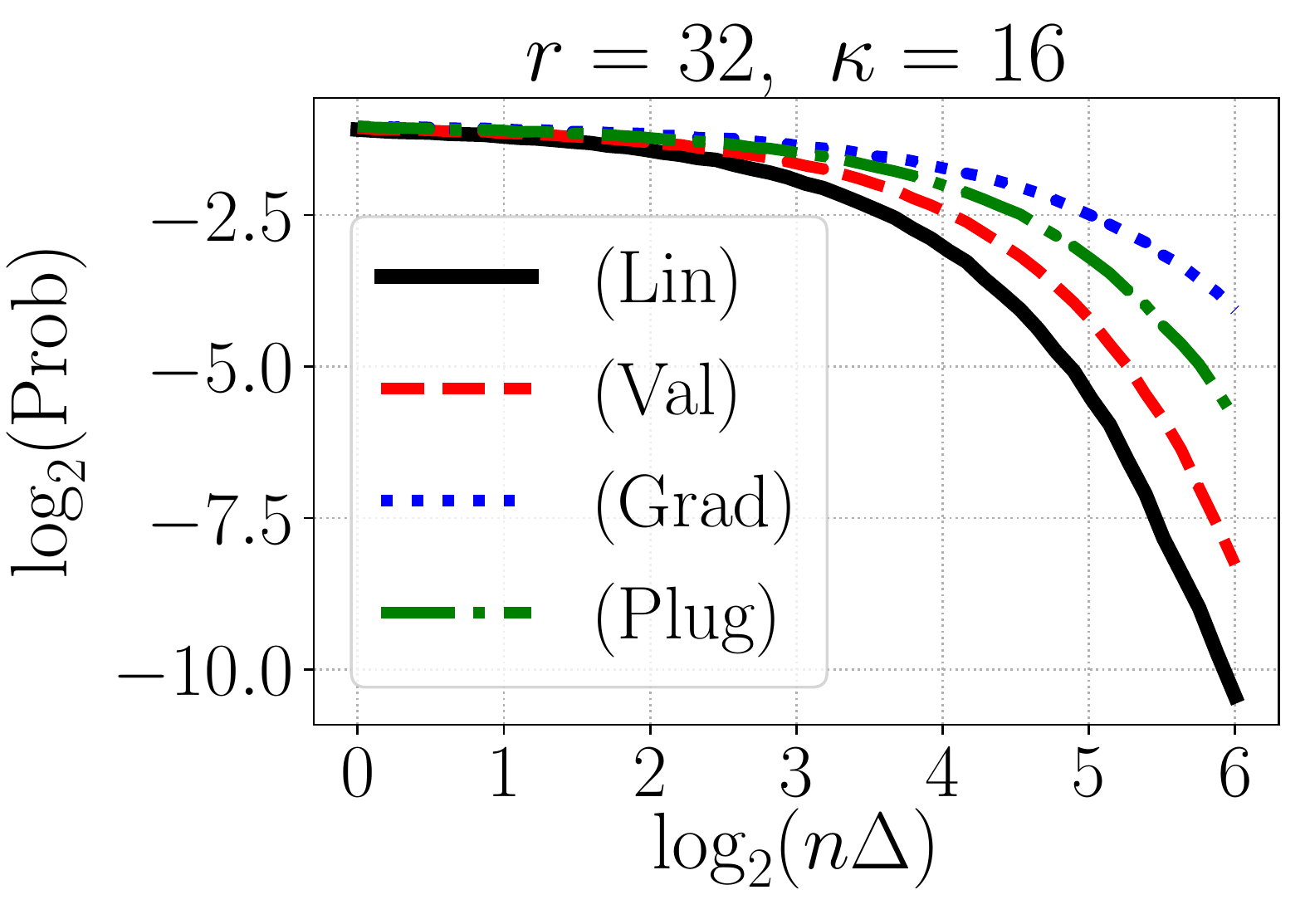}
\includegraphics[width=1.04\textwidth,clip=true,angle=0]{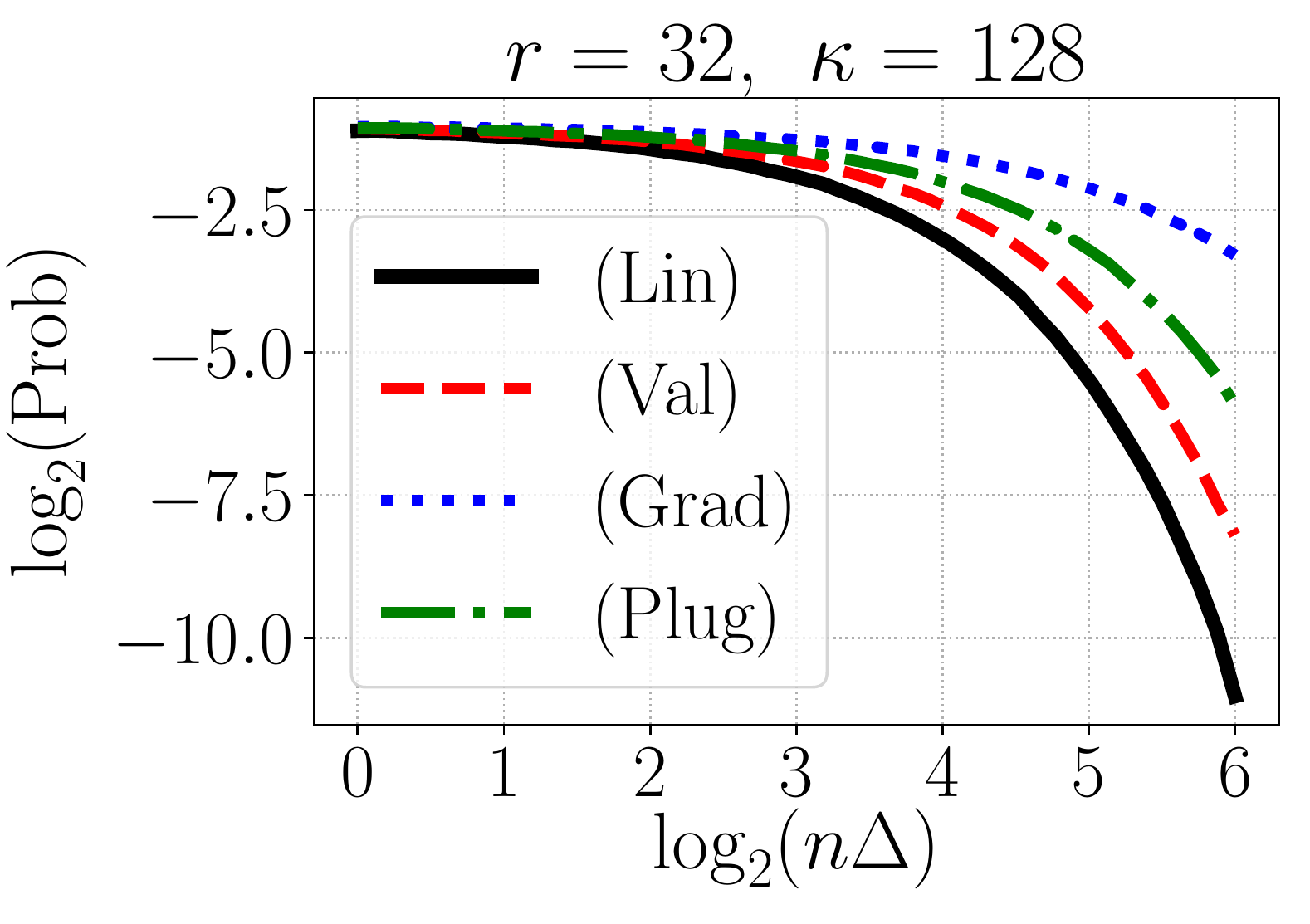}
\end{minipage}
\begin{minipage}{0.32\textwidth}
\centering
\includegraphics[width=0.97\textwidth,clip=true,angle=0]{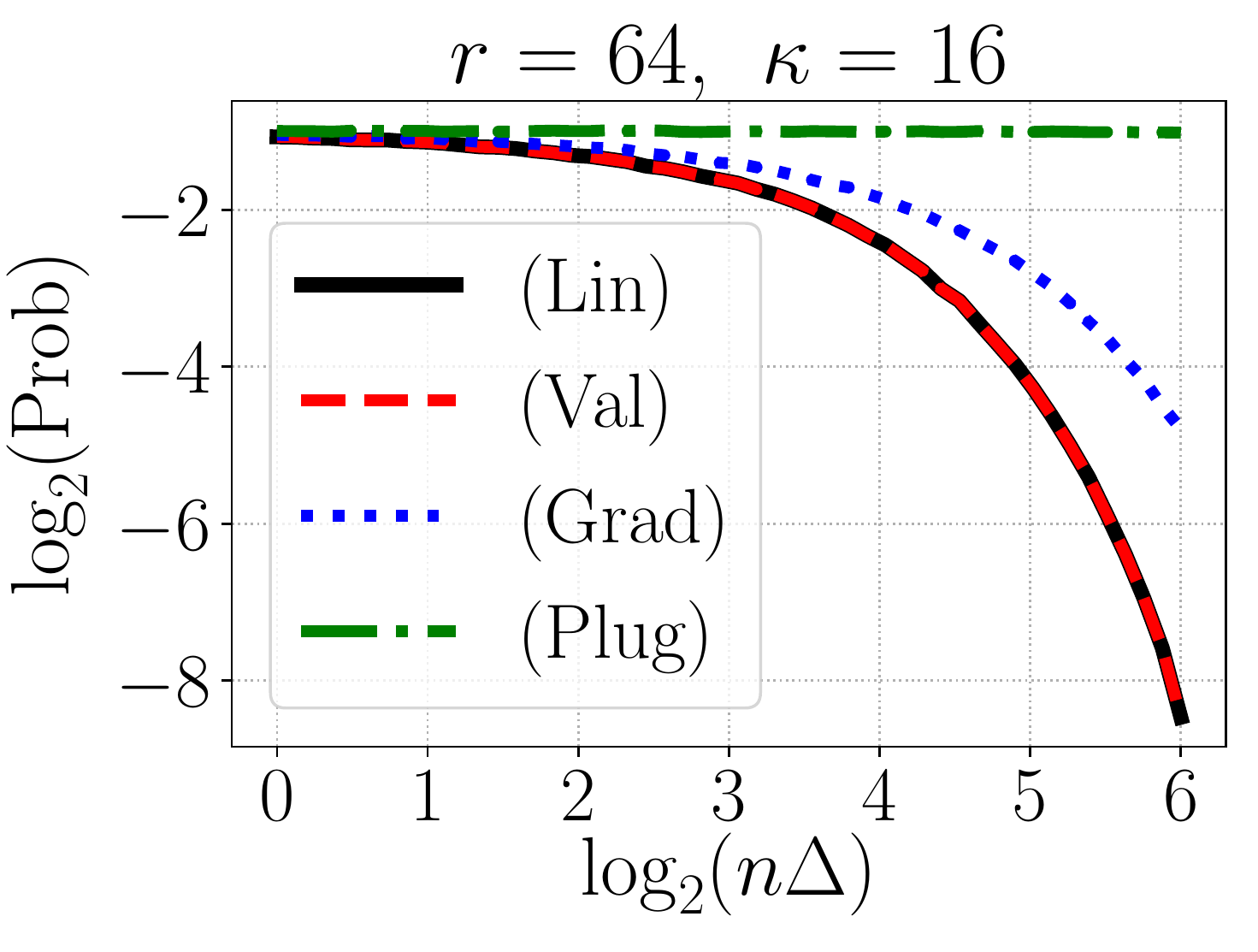}
\end{minipage}
\caption{Comparison of tests~\eqref{eq:test_1},~\eqref{eq:test_naive},~\eqref{eq:test_grad}, and~\eqref{eq:test_plug-in} 
in the well-specified linear model scenario. 
From up to down, we increase the condition number of the design covariance matrix on its span. 
From left to right we increase the covariance matrix rank.}
\label{fig:simul-gaussian}
\end{figure}

\paragraph{Contenders.}
We compare several tests in the scenario of a well-specified linear model (cf. Sec.~\ref{sec:gaussian}):
\begin{itemize}
\item[(a)]
our proposed test~\eqref{eq:test_1}; 
\item[(b)]
the test~\eqref{eq:test_naive} that compares the values of empirical risk against those of the population risk and reduces to~\eqref{eq:test_1-small-sample};
\item[(c)] 
the test~\eqref{eq:test_grad} that compares the~$\ell_2$-norms of the empirical risk gradients (cf.~Remark~\ref{rem:gradient-test}); 
\item[(d)] 
the plug-in test given by
\begin{equation}
\label{eq:test_plug-in}
\ind\{ \| Y^{(0)} - X^{(0)}\theta^*\|^2 + \| Y^{(1)} - X^{(1)} (\wh\theta_0 + \wh\theta_1 - \theta^*)\|^2 
\geq 
\| Y^{(0)} - X^{(0)}(\wh\theta_0 + \wh\theta_1 - \theta^*)\|^2 + \| Y^{(1)} - X^{(1)} \theta^{*}\|^2 \}
\tag{$\textsf{\textup{Plug}}$}
\end{equation}
where~$\wh \theta_0, \wh \theta_1$ are the least-squares estimates of~$\theta_0, \theta_1$. 
\end{itemize}
The motivation behind test~\eqref{eq:test_plug-in} is to mimic the conceptual test
\begin{equation}
\label{eq:test-oracle-sims}
\ind\{  \| Y^{(0)} - X^{(0)}\theta^*\|^2 + \| Y^{(1)} - X^{(1)} \bar\theta\|^2 \geq \| Y^{(0)} - X^{(0)} \bar\theta\|^2 + \| Y^{(1)} - X^{(1)} \theta^*\|^2\}
\end{equation}
which is unavailable since it uses the complementary model~$\bar\theta = \theta_0 + \theta_1 - \theta^*$.
In fact,~\eqref{eq:test-oracle-sims} is the likelihood-ratio test for the simplified testing problem in which~$\bar\theta$ is known, and one tests the two simple hypotheses:
$(\theta_0, \theta_1) = (\theta^*, \bar\theta)$ against~$(\theta_0, \theta_1) = (\bar\theta, \theta^*)$.
(Recall also the discussion after the formulation of Theorem~\ref{th:lower-gaussian} in Section~\ref{sec:gaussian-lower} and the proof of this theorem in Appendix~\ref{app:lower} for more details.) 
Hence, test~\eqref{eq:test-oracle-sims} is optimal in this problem by the Neyman-Pearson lemma. 
As such, one can hope that~\eqref{eq:test_plug-in} performs well in the initial testing problem (with unknown~$\bar\theta$) when the sample size is large enough for the plug-in estimate~$\wh\theta_0 + \wh \theta_1 - \theta^*$ of~$\bar\theta$ to be sufficiently close to~$\bar \theta$.
In fact, it is easily seen that test~\eqref{eq:test_plug-in} is optimal in the asymptotic regime considered in Section~\ref{sec:asymp} 
(i.e., when~$n_k \to \infty$ and~$\Delta_k \to 0$ such that~$n_k\Delta_k \to \lambda_k > 0$, and~$r_k$ is fixed). 
Indeed, in this limiting regime the additional terms in the statistic of~\eqref{eq:test_plug-in} (arising due to the fluctuations~$\wh\theta_0,\wh\theta_1$) vanish; 
on the other hand, the offset terms~$n_k\Delta_k$ (with which the fluctuations are compared) remain constant. 
\paragraph{Experimental setup and results.}
We take~$\Prob_0$~$\Prob_1$ as in~\eqref{eq:gaussian-distributions}, that is
$\Prob_k: x \sim \cN(0,\bSigma),$~$y|x \sim \cN(x^\top \theta_{k},1)$ for~$k \in \{0,1\}.$
The identical covariance matrices~$\bSigma_0 = \bSigma_1 = \bSigma$ are diagonal with size~$d$ and even rank~$r$ to be specified later. 
(Note that~$d \ge r$ does not influence the performance of any of the four tests, so we take~$d = 2r$ for simplicity.) 
Specifically,~$\bSigma_0$ has~$1$ in the first~$r/2$ positions on the diagonal and some~$\kappa > 0$, to be specified later, in the next~$r/2$ positions. 
Moreover, we take~$\theta_0 = 0$ and~$\theta_1 = [\sqrt{\Delta}; 0; ...; 0] \in \R^{d}$ for~$\Delta$ to be specified later; clearly, this indeed gives~$\|\bSigma^{1/2}(\theta_1 - \theta_0)\|^2 = \Delta$.
Finally, we fix the common sample size~$n = 64$, and change~$r \in \{8,32,64\}$ and~$\kappa~\in \{2, 16, 128\}$. 
For each pair~$(r,\kappa)$, we measure the type I error frequency for each test over~$T = 75\,000$ Monte-Carlo trials, with~$n\Delta$ on a logarithmic grid, and plot the resulting curve in the~$\log$-$\log$ scale.\\ 
The results of this experiment are presented in Figure~\ref{fig:simul-gaussian}.\footnote{Python codes for this experiment are available at~\url{https://github.com/ostrodmit/testing-without-recovery}.}
\begin{itemize}
\item
First, we see that test~\eqref{eq:test_1} has the best performance among the four tests; as expected from Theorem~\ref{th:upper-random-design}, its performance deteriorates as~$r$ grows. 
\item 
Second, the only test sensitive to the condition number~$\kappa$ is~\eqref{eq:test_grad}, as is expected since it is the only one of the four with non affine-equivariant statistic. 
As expected from our analysis in Appendix~\ref{app:grad}, the performance of~\eqref{eq:test_grad} degrades as the problem becomes more ill-conditioned.
\item
Third, as we change~$r$, we observe the following effect.
On one hand, the statistical performance of~\eqref{eq:test_naive} is rank-independent (recall~the discussion in the end of~Section~\ref{sec:discussion}). 
On the other hand,~\eqref{eq:test_1} performs the better the smaller is~$r \le n$, as the effect of noise reduction due to projection of residuals becomes more pronounced as we decrease the dimension~$r [\le n]$ of the signal subspace. As a result, there is a gap between the two curves that decreases as we increase the rank, and vanishes when~$r = n = 64$, i.e., when there is effectively no projection.
\item 
Finally, as we increase~$r$ the performance of~\eqref{eq:test_plug-in} deteriorates compared to~\eqref{eq:test_naive} and thus to~\eqref{eq:test_1} as well. 
This is because~\eqref{eq:test_plug-in} needs sufficient sample size to mimic the oracle test~\eqref{eq:test-oracle-sims}.

\end{itemize}

\section*{Acknowledgments}
A.~Javanmard is partially supported by the Sloan Research Fellowship in mathematics, an Adobe Data Science Faculty Research Award and the NSF CAREER Award DMS-1844481. M.~Ndaoud is partially supported by the James H. Zumberge Faculty Research and Innovation Fund and the NSF grant CCF-1908905.

\clearpage

\bibliographystyle{unsrt}
\bibliography{references}
 \clearpage
\appendix


\begin{center}
{\Large{{\bf Supplementary Materials for “Near-Optimal Procedures for Model Discrimination  with Non-Disclosure Properties”}}}
\end{center}

\section{Proofs for upper bounds in linear models}
\label{app:upper-linear}

\subsection{Proof of Proposition~\ref{prop:upper-fixed-design}}
Observing that~$\| X^{(1)}(\theta_0 - \theta_1)\|^2 = n_1\wh\Delta_1$ and recalling the prior decomposition of the testing error we have
\begin{align*}
&\Prob_{\cH_0}[\hat{T} = 1 | X^{(0)}, X^{(1)}] \\
&= \Prob \big[ \| \xi^{(0)}\|^2_{X^{(0)}} - \wh r_0  + \wh r_1 - \| \xi^{(1)}\|^2_{X^{(1)}} + 2\langle X^{(1)} (\theta_1 - \theta_0), \xi^{(1)} \rangle \geq n_1\wh \Delta_1 \big]  \\ 
&\leq \Prob \bigg[ \| \xi^{(0)}\|^2_{X^{(0)}} - \wh r_0 \geq \frac{n_1 \wh\Delta_1}{3}\bigg] 
	+ \Prob \bigg[ \wh r_1 - \| \xi^{(1)}\|^2_{X^{(1)}}  \geq \frac{n_1 \wh\Delta_1}{3}\bigg] 
	+ \Prob \bigg[\langle X^{(1)} (\theta_1-\theta_0),\xi^{(1)}\rangle \geq \frac{n_1 \wh\Delta_1}{6} \bigg] \\
&\leq \exp\bigg(-c \min\bigg\{n_1 \wh\Delta_1, \frac{n_1^2 \wh\Delta_1^2}{\wh    r_0}\bigg\}\bigg) + \exp\bigg(-\frac{c n_1^2 \wh\Delta_1^2}{\wh r_1}\bigg) + \exp(-cn_1 \wh\Delta_1)
\end{align*}
Here we omitted conditioning over~$X^{(0)}, X^{(1)}$ in the 2nd and 3rd lines for brevity.
The estimates in the last line rely on the standard Gaussian and chi-squared tail bounds (see~\cite[Lem.~1]{laurent2000adaptive}), namely
\begin{subequations}
\label{eq:tail-bounds}
\begin{align}
\label{eq:tail-bounds-chi-sq-right}
\Prob[\chi_s^2 - s \ge t]  &\le \exp(-c\min\{t, {t^2}/{s}\}), \\
\label{eq:tail-bounds-chi-sq-left}
\Prob[s - \chi_s^2 \ge t]  &\le \exp(-ct^2/s), \\ 
\label{eq:tail-bounds-normal}
\Prob[\cN(0,1) \ge t] 	   &\le \exp(-ct^2/2),
\end{align}
\end{subequations}
applied (conditionally on~$X^{(0)}, X^{(1)}$) to~$\| \xi^{(0)}\|^2_{X^{(0)}} \sim \chi_{\wh r_0}^2$,~$\| \xi^{(1)}\|^2_{X^{(1)}} \sim \chi_{\wh r_1}^2$, and~$\langle X^{(1)} (\theta_1 - \theta_0), \xi^{(1)} \rangle \sim \cN(0, n_1\wh \Delta_1)$. \qed

\subsection{Proof of Theorem~\ref{th:upper-random-design}}

Let~$\wh S = \wh S_0 - \wh S_1$ be the statistic whose sign is examined in~$\wh T$, i.e.,~$\wh T = \ind\{\wh S_0 \ge \wh S_1 \}$. 
Note that~$\wh r_k := \rank(\wh\bSigma_k) = \rank(\bPi_{X^{(k)}})$; 
moreover,~$\wh r_k = \min\{n_k, r_k\}$~{\em a.s.}
Under~$\cH_0$ (w.l.o.g.) we have
\[
\wh S_0 = \|\bPi_{X^{(0)}} \xi^{(0)}\|^2 - \rank(\bPi_{X^{(0)}}) = \|\bPi_{X^{(0)}} \xi^{(0)}\|^2 - \wh r_0.
\]
Clearly,~$\|\bPi_{X^{(0)}} \xi^{(0)}\|^2 \sim \chi_{\wh r_0}^2$ conditionally on~$X^{(0)}$ since~$\xi^{(0)}$ is independent of~$X^{(0)}$. 
Similarly,
\[
\begin{aligned}
\wh S_1 &= \|\bPi_{X^{(1)}} [\xi^{(1)}  + X^{(1)} (\theta_1 - \theta_0)]\|^2 - \rank(\bPi_{X^{(1)}}) \\
        &= \|\bPi_{X^{(1)}} \xi^{(1)}\|^2 - \wh r_1 + \| X^{(1)} (\theta_1 - \theta_0)\|^2 + 2 \langle \xi^{(1)}, X^{(1)} (\theta_1 - \theta_0) \rangle\\
        &= \|\bPi_{X^{(1)}} \xi^{(1)}\|^2 - \wh r_1 + n_1 \wh \Delta_1 + 2\langle \xi^{(1)}, X^{(1)}(\theta_1-\theta_0)\rangle \\
        &= \wh S_1^o + n_1 \wh \Delta_1 + 2\zeta^{(1)},
\end{aligned}
\]
where we defined~$\wh S_1^o := \|\bPi_{X^{(1)}} \xi^{(1)}\|^2 - \wh r_1$ (distributed as~$\chi_{\wh r_1}^2$ conditionally on~$X^{(1)}$) and~$\zeta^{(1)} := \langle \xi^{(1)}, X^{(1)} (\theta_1-\theta_0)\rangle$ (distributed as~$\cN(0,n_1\wh \Delta_1)$ conditionally on~$X^{(1)}$). 
Hence, using~\eqref{eq:tail-bounds} we get
\[
\begin{aligned}
&\Prob_{\cH_0}[\wh S > 0|X^{(0)},X^{(1)}] \\
&\le \Prob\bigg[\wh S_0 \ge \frac{n_1 \wh \Delta_1}{3} \; \bigg| \; X^{(0)}, X^{(1)} \bigg] 
	+ \Prob\bigg[-\wh S_1^o \ge  \frac{n_1 \wh\Delta_1}{3} \; \bigg| \; X^{(1)} \bigg] 
	+ \Prob\bigg[\zeta^{(1)} \ge \frac{n_1 \wh \Delta_1}{6} \; \bigg| \; X^{(1)} \bigg] \\
&\le 2\exp(-c n_1 \wh \Delta_1) 
    +2\exp\bigg(-\frac{c n_1^2 \wh \Delta_1^2}{\max\{\wh r_0, \wh r_1\}} \bigg). 
\end{aligned}
\]
Note that~$\bar r = \max\{\wh r_0, \wh r_1\}$ {\em a.s.}.
It remains to marginalize over~$X^{(0)},X^{(1)}$.
To this end, observe that
\[
\wh \Delta_1 = \frac{1}{n_1} \sum_{i=1}^{n_1} \langle x_i^{(1)}, \theta_1 - \theta_0 \rangle^2,
\]
where~$x_i^{(1)} \sim \cN(0,\bSigma_1)$ are~i.i.d.
As a result,~$n_1\wh \Delta_1/\Delta_1 \sim \chi_{n_1}^2$ and
$
\Prob [ n_1\wh \Delta_1 \le {n_1 \Delta_1}/{2} ] \le \exp(-cn_1),
$
cf.~\eqref{eq:tail-bounds-chi-sq-left}.
Thus,
$
\Prob_{\cH_0}[\wh S > 0] \le \exp(-cn_1) + 2\exp\left(-c n_1 \Delta_1 \right) 
    +2\exp\left(-{c n_1^2 \Delta_1^2}/{\bar r}\right)
$
as claimed.
\qed

\subsection{Proof of Corollary~\ref{cor:complexity-gaussian}}  

For simplicity we assume~$\Delta_0 = \Delta_1 = \Delta$; the general case is similar. 
Recall that~$\Delta^2 \le 1$. 
Recall also that~$r := \max \{r_0, r_1\}$ and~$\bar r := \min\{r,n\}$. 

\proofstep{1}. We first look at the reverse direction: assume that
\begin{equation}
\label{eq:upper-random-design-log}
cn \Delta \min\left\{1,\frac{n \Delta}{\bar r}\right\} \ge \log(1/\delta)
\end{equation}
which, by~\eqref{eq:upper-random-design}, corresponds to~$\max\{\Prob_{\cH_0}[\wh T = 1],\Prob_{\cH_1}[\wh T = 0]\} \le 2C \delta = O(\delta)$.
Choosing the second term in the minimum in~\eqref{eq:upper-random-design-log}, we get
\[
n \Delta \ge \sqrt{\log(1/\delta) \bar r/c} = \sqrt{\log(1/\delta) \min\{r,n\}/c},
\]
whence
\begin{equation}
\label{eq:complexity-gaussian-suffices}
n = \Omega \left( \min\left\{ \frac{\log(1/\delta)}{\Delta^2}, \frac{\sqrt{r\log(1/\delta)}}{\Delta} \right\} \right).
\end{equation}
On the other hand, choosing the first term in the minimum in~\eqref{eq:upper-random-design-log} yields~$n = \Omega(\log(1/\delta)/\Delta)$, which corresponds to the missing (so far) term in the sample complexity bound~\eqref{eq:complexity-gaussian} (recall that~$\Delta \le 1$).
Thus,~\eqref{eq:complexity-gaussian} is indeed {\em necessary} to guarantee that~$\max\{\Prob_{\cH_0}[\wh T = 1],\Prob_{\cH_1}[\wh T = 0]\} = O(\delta)$ as per~\eqref{eq:upper-random-design}. 

\proofstep{2}.
Conversely, assume that
\begin{equation}
\label{eq:complexity-gaussian-suffices-rw}
cn\Delta \ge \min\left\{\frac{\log(1/\delta)}{\Delta}, \sqrt{r\log(1/\delta)} + \log(1/\delta) \right\}, 
\end{equation}
which corresponds to the sufficient sample size as per~\eqref{eq:complexity-gaussian}.
Meanwhile,~\eqref{eq:upper-random-design-log} can be rewritten as
\begin{equation}
\label{eq:upper-random-design-log-rw}
\begin{aligned}
c n\Delta 
\ge \max\left\{1, \frac{\bar r}{n\Delta} \right\} \log(1/\delta)
&= \max\left\{1, \min\left\{\frac{1}{\Delta}, \frac{r}{n\Delta}\right\}\right\} \log(1/\delta)\\
&= \min\left\{\frac{1}{\Delta}, \max\left\{1, \frac{r}{n\Delta}\right\}\right\} \log(1/\delta),
\end{aligned}
\end{equation}
where in the last step we used that~$\Delta \le 1$. 
Clearly, the first case for the minimum in~\eqref{eq:complexity-gaussian-suffices-rw} is identical to the first case for the minimum in~\eqref{eq:upper-random-design-log-rw}. 
On the other hand, the second case for the minimum in~\eqref{eq:complexity-gaussian-suffices-rw} implies that~$n\Delta \ge \log(1/\delta)$ and~$n^2\Delta^2 \ge r \log(1/\delta)$, i.e., the second case in~\eqref{eq:upper-random-design-log-rw}.
\qed

\subsection{Proof of Proposition~\ref{prop:upper-fixed-design-adapt}}
Recall the decomposition of the type I error (all probabilities are conditional on~$X^{(0)},X^{(1)}$):
\begin{align}
&\; 
\Prob_{\cH_0}[\hat{T} = 1|X^{(0)}, X^{(1)}] \notag\\
&= \Prob \bigg[ \frac{1}{\wh \tau_0} \|\xi^{(0)}\|^2_{X^{(0)}} - \wh r_0  + \wh r_1 - \frac{1}{\wh \tau_1}\| \xi^{(1)} \|^2_{X^{(1)}} + \frac{2}{\wh \tau_1 \sigma_1} \langle X^{(1)} (\theta_1 - \theta_0), \xi^{(1)} \rangle \geq \frac{n\wh \Delta_1}{\wh\tau_1 \sigma_1^2} \bigg]  \notag\\ 
&= \Prob \bigg[ \frac{\wh \tau_1}{\wh \tau_0} \|\xi^{(0)}\|^2_{X^{(0)}} - \wh\tau_1 \wh r_0  + \wh\tau_1 \wh r_1 - \|\xi^{(1)}\|^2_{X^{(1)}} + \frac{2}{\sigma_1} \langle X^{(1)} (\theta_1 - \theta_0), \xi^{(1)} \rangle \geq \frac{n\wh \Delta_1}{\sigma_1^2} \bigg]  \notag\\
%
&\leq \Prob \bigg[ \wh r_1 - \|\xi^{(1)}\|^2_{X^{(1)}}  \geq \frac{n \wh\Delta_1}{5 \sigma_1^2}\bigg] 
	+ \Prob \bigg[\frac{1}{\sigma_1} \langle X^{(1)} (\theta_1 - \theta_0), \xi^{(1)} \rangle \geq \frac{n \wh\Delta_1}{10 \sigma_1^2} \bigg] \notag\\
&\quad 
	+ \Prob \bigg[ \frac{\wh\tau_1}{\wh \tau_0} \left( \|\xi^{(0)}\|^2_{X^{(0)}} - \wh r_0 \right) \geq \frac{n \wh\Delta_1}{5 \sigma_1^2}\bigg]
	+ \Prob \bigg[ (\wh \tau_1 - 1) \wh r_1 \ge \frac{n \wh\Delta_1}{5 \sigma_1^2} \bigg] 
	+ \Prob \bigg[ \wh \tau_1 \bigg(\frac{1}{\wh \tau_0} - 1\bigg) \wh r_0 \ge \frac{n \wh\Delta_1}{5 \sigma_1^2} \bigg].
\label{eq:gaussian-adapt-expansion}
\end{align}
As in the proof of Proposition~\ref{prop:upper-fixed-design}, the sum of the terms in the penultimate line of~\eqref{eq:gaussian-adapt-expansion} is bounded by
\[
2\exp\bigg(-\frac{c n \wh\Delta_1}{\sigma_1^2} \min \bigg\{ 1, \frac{n \wh\Delta_1}{\wh r_1 \sigma_1^2} \bigg\}\bigg).
\]
Using that~$\wh \tau_0,$~$\wh \tau_1$ and~$\| \xi^{(0)}\|^2_{X^{(0)}} - \wh r_0$ are conditionally independent and~$\Prob[\wh \tau_k \notin [\tfrac{1}{2}, 2]]  \le C\exp(-cn)$ for~$k \in \{0,1\}$ due to~\eqref{eq:var-ratio-tail-bound}, applying~\eqref{eq:tail-bounds} and the union bound, we bound the first term in the last line of~\eqref{eq:gaussian-adapt-expansion} by
\[
\begin{aligned}
\Prob \bigg[ \| \xi^{(0)}\|^2_{X^{(0)}} - \wh r_0 \geq \frac{n \wh\Delta_1}{5 \sigma_1^2}\bigg] + C\exp(-cn)
\le  
\exp\bigg(-\frac{c n \wh\Delta_1}{\sigma_1^2} \min\bigg\{1, \frac{n \wh\Delta_1}{\wh r_0 \sigma_1^2}\bigg\}\bigg) 
+ C\exp(-cn).
\end{aligned}
\]
Moreover, applying~\eqref{eq:var-ratio-tail-bound} to~$\wh \tau_1$ and noting that~$\wh r_1 \le n$~{\em a.s.}, we bound the next term as
\[
\Prob \bigg[ (\wh \tau_1 - 1)  \ge \frac{n \wh\Delta_1}{5 \wh r_1 \sigma_1^2} \bigg] 
\le 
2\exp \bigg( - \frac{c_3 n^2 \wh \Delta_1}{\wh r_1 \sigma_1^2} \min \bigg\{ 1, \frac{n \wh \Delta_1}{\wh r_1 \sigma_1^2}\bigg\} \bigg) 
\le 
2\exp \bigg( - \frac{c_3 n \wh \Delta_1}{\sigma_1^2} \min \bigg\{ 1, \frac{n \wh \Delta_1}{\wh r_1 \sigma_1^2} \bigg\} \bigg).
\]
In order to estimate the last term in the right-hand side of~\eqref{eq:gaussian-adapt-expansion}, we use that~$\Prob[\wh \tau_1 \ge 2] \le 2\exp(-cn)$ by~\eqref{eq:var-ratio-tail-bound} and that~$\wh \tau_0$ and~$\wh \tau_1$ are conditionally independent. 
We then have
\[
\begin{aligned}
\Prob \bigg[ \wh \tau_1 \bigg(\frac{1}{\wh \tau_0} - 1\bigg) \ge \frac{n \wh\Delta_1}{5 \wh r_0 \sigma_1^2} \bigg]	
&\le 
2\exp(-cn) + \Prob \bigg[ \bigg(\frac{1}{\wh \tau_0} - 1\bigg) \ge \frac{n \wh\Delta_1}{10 \wh r_0 \sigma_1^2} \bigg] \\
&\le 
2\exp(-cn) + \Prob \bigg[ \wh \tau_0 \le 1 - \frac{1}{2}\min \bigg\{ \frac{n \wh\Delta_1}{10 \wh r_0 \sigma_1^2}, 1\bigg\} \bigg] \\
&\le 
4\exp(-c_4 n) + 2\exp \bigg( - \frac{c_4 n \wh \Delta_1}{\sigma_1^2} \min \bigg\{ 1, \frac{n \wh \Delta_1}{\wh r_0 \sigma_1^2} \bigg\} \bigg)
\end{aligned}
\]
where in the second step we used that~$\tfrac{1}{1+t} \le 1-\tfrac{1}{2}\min[t,1]$ for~$t \ge 0$, and in the end we used~\eqref{eq:var-ratio-tail-bound}.
\qed

\section{Proof of Theorem~\ref{th:lower-gaussian}}

\label{app:lower} 

Our goal in this section is to prove Theorem~\ref{th:lower-gaussian}. 
We proceed in four steps correspondingly implemented in Sections~\ref{sec:lower-fixed-design-reduction}--\ref{sec:lower-orths}.
For simplicity, we introduce the following concise notation for the samples:
\[
(X,Y) = (X^{(0)},Y^{(0)}), \quad (X',Y') = (X^{(1)},Y^{(1)}).
\] 
Similarly, we denote~$\xi := \xi^{(0)} [= Y - X\theta_0]$ and~$\xi' := \xi^{(1)} [= Y' - X'\theta_1]$.

\subsection{Reduction to fixed-design setup}

\label{sec:lower-fixed-design-reduction}

We first show that one can pass to the fixed-design setup -- namely, focus on proving the bound
\begin{equation}
\label{eq:lower-fixed-design}
\underset{\wh T}{\inf} \underset{(\theta_0,\theta_1) \in \wh \Theta}{\sup} 
\left\{ \Prob_{\cH_0}[\wh T = 1|X,X'] + \Prob_{\cH_1}[\wh T = 0|X,X'] \right\}
\ge  
C \exp \left(-c n\Delta \min \left\{ 1, \frac{n\Delta}{\rmin} \right\} \right).
\end{equation}
Here the infimum is over all measurable maps~$(\theta^*, X, Y, X', Y') \to \{0,1\}$ as in~\eqref{eq:lower-gaussian}, the supremum over the set
\begin{equation}
\label{def:set-empirical}
\wh \Theta = \wh \Theta(32\Delta) := \{(\theta_0,\theta_1): \max\{\wh \Delta_0, \wh \Delta_1\} \geq 32\Delta\}
\end{equation} 
with~$\wh \Delta_{0},\wh \Delta_{1}$ defined in~\eqref{def:Delta-empirical}; finally,
$
\wh\rmin := \min\{ \rank(\wh \bSigma_0), \rank(\wh \bSigma_1) \}
$
with~$\wh \bSigma_k = \frac{1}{n} X^{(k)\top} X^{(k)}$. 
Moreover, when proving~\eqref{eq:lower-fixed-design}, we shall assume~w.l.o.g. that~$r_0 \leq r_1$, and recall that~$\rmin \geq 28$.
\newcommand{\II}{\mathsf{\Gamma}}
Finally, from now on we assume that 
\begin{equation}
\label{eq:evals-condition}
\lambda_{\min} (\wh \bSigma_{0,\II}) , \lambda_{\max} (\wh \bSigma_{1,\II}) \in 
\left[ \frac{1}{16}, \frac{49}{16} \right],
\end{equation}
where $\wh \bSigma_{k,\II}$ is the top left submatrix of $\wh \bSigma_k$ (for~$k \in \{0,1\}$) at the intersection of the first~$\lfloor\rmin/4 \rfloor$ rows and columns; this implies, in particular, that
\[
\rank(\wh \bSigma_{0,\II}) = \left\lfloor {\rmin}/{4} \right\rfloor.
\]
This is a non-restrictive assumption: indeed, in Section~\ref{sec:lower-orths} we shall prove that
\begin{equation}
\label{eq:evals-condition-prob}
\Prob\left\{ \lambda_{\min} (\wh \bSigma_{0,\II}) , \lambda_{\max} (\wh \bSigma_{1,\II}) \in \left[ \frac{1}{16},  \frac{49}{16} \right] \right\} \ge \frac{1}{8}. 
\end{equation}
When combined with~\eqref{eq:lower-fixed-design},~\eqref{eq:evals-condition}
will readily yield~\eqref{eq:lower-gaussian}.
Indeed, by Markov's inequality
\[
\wh \Delta_0 \le 32\Delta_0 \;\; \text{and} \;\; \wh \Delta_1 \le 32\Delta_1
\] 
simultaneously with prob.~$\ge 15/16$ over~$X,X'$. 
Under the joint event (and so w.p.~$\ge 1-7/8-1/16 = 1/16$) one has that~$\wh\Theta(32\Delta) \subseteq \Theta(\Delta)$, i.e., the supremum in~\eqref{eq:lower-gaussian} is over a larger set than in~\eqref{eq:lower-fixed-design}; thus, 
\[
\begin{aligned}
&\underset{\wh T}{\inf} \underset{(\theta_0,\theta_1) \in \Theta(\Delta)}{\sup} \left\{ \Prob_{\cH_0}[\wh T=1] + \Prob_{\cH_1}[\wh T=0] \right\} \\
\ge 
&\frac{1}{16} \underset{\wh T}{\inf} \underset{(\theta_0,\theta_1) \in \wh \Theta(32\Delta)}{\sup} \left\{ \Prob_{\cH_0}[\wh T=1|X,X'] + \Prob_{\cH_1}[\wh T=0|X,X'] \right\} \\
\ge
&\frac{C}{16} \exp\left(-c n\Delta \min \left\{ 1, \frac{n\Delta}{\rmin}\right\} \right),
\end{aligned}
\]
which is~\eqref{eq:lower-gaussian}. 
Thus, we can indeed focus on proving~\eqref{eq:lower-fixed-design} with fixed design matrices~$X,X'$ obeying~\eqref{eq:evals-condition}.
This is our next goal; Secs.~\ref{sec:lower-easy}--\ref{sec:lower-hard} are dedicated to reaching it, and in Sec.~\ref{sec:lower-orths} we shall verify~\eqref{eq:evals-condition-prob}.

\paragraph{What is ahead.}
Note that~\eqref{eq:lower-fixed-design} comprises two bounds simultaneously: the rank-independent bound~$\exp(-cn\Delta)$ and the second bound which depends on~$\rmin$.
These two bounds shall be proved correspondingly in Sections~\ref{sec:lower-easy} and~\ref{sec:lower-hard}. 
It turns out that for the first bound it suffices to consider a further simplified testing problem, where one observes {\em both models}~$\{\theta^*, \bar\theta\}$ rather than just~$\theta^*$,
where
\[
\bar{\theta} := \theta_0 + \theta_1 - \theta^*
\]
is the complementary to~$\theta^*$ model in~$\{\theta_0,\theta_1\}$.
The problem then reduces to discriminating between the two {\em simple} hypotheses~$(\theta_0, \theta_1) = (\theta^*, \bar\theta)$ and~$(\theta_0, \theta_1) = (\bar\theta, \theta^*)$.
Intuitively, the absence of the rank in the corresponding sample complexity bound can be explained by the fact that the testing problem reduces to a one-dimensional one since~$\theta_1-\theta_0$ is known.
On the other hand, capturing the dependency on~$\rmin$ requires to take into account that the complementary model~$\bar\theta$ actually is {\em not} observed. We capture this by putting a Gaussian prior on~$\bar \theta$ and carefully bounding the Bayes risk.

\subsection{Rank-independent bound}
\label{sec:lower-easy}
Let~$\bar{\theta} := \theta_0 + \theta_1 - \theta^*$, so that~$\{\theta^*, \bar\theta\} = \{\theta_0, \theta_1\}$ as the two sets regardless of the hypothesis. 
Now, observe that the two hypotheses~$\cH_0: \theta_0 = \theta^*$ and~$\cH_1: \theta_1 = \theta^*$ can be reformulated as the hypotheses~$(\theta_0, \theta_1) = (\theta^*, \bar\theta)$
and~$(\theta_0, \theta_1) = (\bar\theta, \theta^*)$
about the aggregated parameter~$(\theta_0, \theta_1)$. 
Since~$\bar\theta$ is unknown, we are dealing with a composite testing problem, and thus cannot apply the Neyman-Pearson lemma. To put us in the two-point hypothesis testing situation, we can simply fix a value of~$\bar\theta$, and assume that it is also known -- in other words, consider the following testing problem:
\begin{equation}
\label{eq:two-point-problem}
\boxed{
\begin{aligned}
\textit{Given the data}\;
n,\wh r_0, \wh r_1,(\theta^*, \bar\theta) \in \wh\Theta, \; 
(X,Y) \in \Prob_0^{\otimes n},\; (X',Y') \sim \Prob_1^{\otimes n}, \\
\textit{test} \;\;
\bar\cH_0: \; (\theta^*, \bar\theta) = (\theta_0, \theta_1) \; \text{\em against} \;\;
\bar\cH_1:  (\theta^*, \bar\theta) = (\theta_1, \theta_0).
\end{aligned}
\tag{$\textsf{2-point}$}
}
\end{equation}
Clearly,~\eqref{eq:two-point-problem} is a problem of testing between two {\em simple} hypotheses about the unknown parameter~$(\theta_0, \theta_1) \in \wh\Theta$. 
By the Neyman-Pearson lemma, the sum of type I and II error probabilities for any test in~\eqref{eq:two-point-problem} is lower-bounded by that for the likelihood-ratio test, and this test writes
\begin{equation}
\label{eq:test-oracle}
\wh T_{\LR} = \ind \left\{  \| Y - X\theta^*\|^2 + \| Y' - X' \bar\theta\|^2 \geq \| Y - X\bar\theta\|^2 + \| Y' - X' \theta^*\|^2 \right\}.
\end{equation}
On the other hand,~\eqref{eq:two-point-problem} cannot be harder than the initial testing problem considered in~\eqref{eq:lower-fixed-design}. 
Indeed, any test admissible in~\eqref{eq:lower-fixed-design}
is a mapping~$(\theta^*, \bar \theta, X,Y, X', Y') \to \{0,1\}$ that does {\em not} depend on~$\bar \theta$, whereas in~\eqref{eq:two-point-problem} one is also allowed to use the tests that do depend on~$\bar \theta$ and thus are not~$(\theta^*,X,Y, X', Y')$-measurable. 
As such, the left-hand side of~\eqref{eq:lower-fixed-design} admits the bound
\[
\begin{aligned}
\underset{\wh T: (\theta^*,X,Y,X',Y') \to \{0,1\}}{\inf} \;
	&\underset{(\theta_0, \theta_1) \in \wh\Theta}{\sup} 
\left\{\Prob_{\cH_0}[\wh T = 1|X,X'] + \Prob_{\cH_1}[\wh T = 0|X,X'] \right\} \\
\geq 
	&\underset{(\theta^*, \bar\theta) \in \wh\Theta}{\sup} 
\left\{ \Prob_{\bar\cH_0} [\wh T_{\LR} = 1|X,X'] + \Prob_{\bar\cH_1} [\wh T_{\LR} = 0|X,X'] \right\}, 
\end{aligned}
\]
with~$\wh T_{\LR}$ being the test given by~\eqref{eq:test-oracle}. 
Fixing~$(\theta^*, \bar\theta) \in \wh\Theta = \wh \Theta(32\Delta)$, under~$\bar \cH_0$ we get
\begin{align*}
&\Prob_{\bar\cH_0}[\wh T_{\LR} = 1|X,X'] \\ 
&\geq
\Prob \left\{ \| Y - X \theta_0\|^2 + \| Y' - X' \theta_1\|^2 \geq \| Y - X \theta_1\|^2 + \| Y' - X' \theta_0\|^2 \; \middle| \; X,X' \right\} \\
&= 
\Prob \left\{ 2\langle \xi,X(\theta_0-\theta_1) \rangle + 2\langle \xi', X'(\theta_0-\theta_1) \rangle \geq \|X(\theta_0 - \theta_1)\|^2 + \| X'( \theta_0 - \theta_1)\|^2 \middle| X,X'  \right\} \\
&\geq 
C\exp\big(-c\max \{\|X(\theta_0-\theta_1)\|^{2}, \|X'(\theta_0-\theta_1)\|^{2}\}\big) 
=
C\exp\big(-cn\max\{\wh \Delta_0,\wh\Delta_1\}\big).
\end{align*}
Here we first used that~$\langle \xi^{(k)},X^{(k)}(\theta_0-\theta_1) \rangle \sim \cN(0,\|X^{(k)}(\theta_0-\theta_1)\|^2)$ for~$k \in \{0,1\}$
conditionally on~$(X,X')$, and then applied the lower bound for the Gaussian tails (see~\cite[Eq.~7.1.13]{stegun}):
\[
\Prob[\cN(0,1) \ge t] \ge C\exp(-ct^2), \quad \forall t \ge 0.
\]
By symmetry, we also have the same bound for~$\Prob_{\bar\cH_1}[\wh T_{\LR} = 0|X,X']$. 
Now, taking~$(\theta_*, \bar\theta)$ on the boundary of~$\wh\Theta(32\Delta)$, i.e., such that~$\max[\wh \Delta_0, \wh \Delta_1] = 32\Delta$, we arrive at 
\begin{equation}
\label{eq:lower-easy}
\begin{aligned}
\underset{\wh T: (\theta^*,X,Y,X',Y') \to \{0,1\}}{\inf} \;
	\underset{(\theta_0, \theta_1) \in \wh\Theta}{\sup} 
\left\{\Prob_{\cH_0}[\wh T = 1|X,X'] + \Prob_{\cH_1}[ \wh T = 0|X,X'] \right\} 
\ge  C \exp (-c' n\Delta).
\end{aligned}
\end{equation}
This is precisely the rank-independent part of~\eqref{eq:lower-fixed-design}.
Next we deal with the rank-dependent part.

\subsection{Rank-dependent bound}
\label{sec:lower-hard}
We first observe that in the cases~$n\Delta \le 1$ and~$n\Delta \ge \rmin/4$, 
the bound~\eqref{eq:lower-fixed-design} follows from~\eqref{eq:lower-easy} up to a change in constants. 
Hence, it suffices to focus on the range 
\begin{equation}
\label{eq:Delta-range-for-complex-lower}
1 \leq n\Delta \leq {\rmin /4}.
\end{equation} 
Instead of fixing the pair~$(\theta^*,\bar\theta)$, we now only fix~$\theta^*$ 
and put on~$\bar \theta$ a Gaussian prior centered at~$\theta^*$ with covariance corresponding to the constraint~$(\theta^*,\bar \theta) \in \wh\Theta$ (cf.~\eqref{eq:prior-definition}).
 Note that, by symmetry, this is equivalent to putting the same prior on~$(\theta_0, \theta_1)$, or the centered prior with the same covariance on~$\theta_1 - \theta_0$ (cf.~\eqref{eq:prior-prediction-residual}, but the exact expression is not important for what follows in step~\proofstep{1} below). 
Technically, such prior violates the constraint~$(\theta_0, \theta_1) \in \wh \Theta$, so we have to verify that the event
\begin{equation}
\label{eq:constraint-set-event}
\cC := \{(\theta_0, \theta_1) \in \wh \Theta\}
\end{equation}
has sufficient mass under this prior. (Note that this {\em concentration step} is common in the literature on Bayesian lower bounds, see, e.g.,~\cite[Sec.~5.4]{johnstone2017gaussian}.)
To that end, we further proceed in three steps.

\proofstep{1}: \textbf{{\em Reduction to the unconstrained Bayes risk.}}
We let $\pi$ be any prior on~$(\theta_0,\theta_1) \in \R^d$, and let~$\pi_{\cC}$ be the conditioning of~$\pi$ to the event~$\cC$ in~\eqref{eq:constraint-set-event}:
\[
\piC(\cE) = \frac{\pi (\mathcal{\cE} \cap \cC)}{\pi(\cC)}, \quad \forall \cE \subseteq \R^{2d}.
\]
In what follows, we identify events over~$(\theta_0, \theta_1)$ with subsets of~$\R^{2d}$; in particular,~$\cC$ is identified with~$\wh \Theta = \wh\Theta(32\Delta)$.
Now, since~$\pi_\cC$ is supported on~$\cC$, we can bound the left-hand side of~\eqref{eq:lower-fixed-design} as
\begin{equation*}
\begin{aligned}
&\underset{\wh T}{\inf} 
\underset{(\theta_0,\theta_1) \in \wh\Theta}{\sup} \Prob_{\cH_0}[\wh T = 1|X,X'] + \Prob_{\cH_1}[\wh T = 0|X,X'] 
\ge 
\underset{\wh T}{\inf} \; \E_{\pi_\cC}\Big[ \Prob_{\cH_0}[\wh T = 1|X,X'] + \Prob_{\cH_1}[\wh T = 0|X,X']\Big]
\end{aligned}
\end{equation*}
where $\E_{\pi_\cC}$ is the expectation over~$(\theta_0, \theta_1) \sim \pi_{\cC}$, and the infimum on the right is over the same set of admissible tests~$(\theta^*, X,Y,X',Y') \to \{0,1\}$ as on the left.
Now, let~$\wh T_{\pi_{\cC}}$ be the corresponding Bayes test (i.e., the one on which the Bayes risk is attained):
\[
\begin{aligned}
&\underset{\wh T}{\inf} \; \E_{\pi_\cC}\Big[ \Prob_{\cH_0}[\wh T = 1|X,X'] + \Prob_{\cH_1}[\wh T = 0|X,X']\Big] 
= 
\E_{\pi_{\cC}} \left[\Prob_{\cH_0}[\wh T_{\pi_{\cC}} = 1|X,X'] + \Prob_{\cH_1}[\wh T_{\pi_{\cC}} = 0|X,X'] \right].
\end{aligned}
\]
On the other hand, 
\begin{align*}
\E_{\pi} \left[ \Prob_{\cH_0}[ \wh T_{\pi_\cC} = 1|X,X'] 
	+ \Prob_{\cH_1}[\wh T_{\pi_\cC} = 0|X,X'] \right] 
&\le 
\E_{\pi_{\cC}} \left[ \Prob_{\cH_0}[\wh T_{\pi_\cC} = 1|X,X'] 
	+ \Prob_{\cH_1}[  \wh T_{\pi_\cC} = 0|X,X'] \right] \\
	&\quad+ 2\pi(\R^{2d} \setminus \cC).
 \end{align*}
Finally, since~$\wh T_{\pi_{\cC}}$ is generally not the Bayes test for~$\pi$, we have
\[
\begin{aligned}
&\underset{\wh T}{\inf} \; \E_{\pi}\Big[ \Prob_{\cH_0}[ \wh T = 1|X,X'] + \Prob_{\cH_1}[\wh T = 0|X,X']\Big] 
\le 
\E_{\pi}\left[\Prob_{\cH_0}[ \wh T_{\pi_\cC}=1|X,X'] 
	+ \Prob_{\cH_1}[ \wh T_{\pi_\cC}=0|X,X']\right].
\end{aligned}
\]
We conclude that
\begin{equation}
\label{eq:unconstrained}
\begin{aligned}
\underset{\wh T}{\inf} 
\underset{(\theta_0,\theta_1) \in \wh\Theta}{\sup} \Prob_{\cH_0}[\wh T = 1|X,X'] + \Prob_{\cH_1}[\wh T = 0|X,X'] 
&\geq 
\underset{\wh T}{\inf} \, \E_{\pi} \left[\Prob_{\cH_0}[\wh T = 1|X,X'] + \Prob_{\cH_1}[\wh T = 0|X,X'] \right] \\
	&\quad- 2\pi(\R^{2d} \setminus \cC);
\end{aligned}
\end{equation}
thus, it suffices to bound the Bayes risk associated with~$\pi$ from below and~$\pi(\R^{2d} \setminus \cC)$ from above.

\proofstep{2}: \textbf{\em Choosing a prior and showing concentration.}
\newcommand{\cI}{{\mathcal{I}}}
Let us choose~$\pi$ as follows:~$\theta^*$ is fixed at an arbitrary value, and~$\theta_1 - \theta_0$ is a vector supported on~$\cI := \{1, \lfloor r/4 \rfloor \}$ and such that
\begin{equation}
\label{eq:prior-definition}
(\theta_1-\theta_0)_{\cI} \sim \mathcal{N}(0, \rho (\wh \bSigma_{0,\II})^{-1}) 
\;\; \text{with} \;\;
\rho :=  \frac{64\Delta}{\lfloor \rmin/4 \rfloor},
\end{equation}
where~$(u)_{\cI} = [u_1; ...; u_{\lfloor r/4 \rfloor}]$. 
By writing~$\cC = \wh\Theta(32\Delta)$, cf.~\eqref{def:set-empirical}, as~$\cC = \cC_0 \cup \cC_1$ with~$\cC_k = \{\wh \Delta_k \geq 32\Delta\}$ for~$k \in \{0,1\}$, and recalling that~$\wh \Delta_0 = \|\wh \bSigma_0^{1/2} (\theta_1 - \theta_0)\|^2 = \frac{1}{n}\|X (\theta_1 - \theta_0)\|^2$,  we observe that 
\[
\pi(\R^{2d} \setminus \cC) 
\leq \pi(\R^{2d} \setminus \cC_0) 
= \pi ( \{ \wh \Delta_0 \leq 32 \Delta \} ).
\]
Clearly,
\begin{equation}
\label{eq:prior-prediction-residual}
\tfrac{1}{\sqrt{n}}X(\theta_1-\theta_0) \sim \mathcal{N}\left(0, \rho \bPi_{X_\cI}\right),
\end{equation}
where
\begin{equation}
\label{eq:projector-I}
\bPi_{X_\cI} = X_\cI^{\vphantom\top} (X_\cI^\top X_\cI^{\vphantom\top})^{-1} X_\cI^\top = \frac{1}{n} X_\cI^{\vphantom\top} \wh\bSigma_{0,\II}^{-1} X_\cI^\top,
\end{equation}
$X_\cI$ being the submatrix of $X$ collecting the columns indexed by $\cI$.
Therefore, we have~$\wh \Delta_0 \sim \rho \chi^2_{\lfloor \rmin/4 \rfloor}$ and~$\wh\Delta_0 / (32\Delta) \sim \tfrac{2}{\lfloor \rmin/4 \rfloor} \chi^2_{\lfloor \rmin/4 \rfloor}$. In combination with the bound~\eqref{eq:tail-bounds-chi-sq-left} for the left tail of~$\chi^2_{s}$ this gives
\[
\pi ( \R^{2d} \setminus \cC ) \leq \exp \left(-c \rmin \right)
\]
for some~$c > 0$.
Now, recall that we are in the regime~$n\Delta \le \rmin/4$, cf.~\eqref{eq:Delta-range-for-complex-lower}, and thus
\[
\exp(-c\rmin) \le \exp\left(-\frac{ c'n^2\Delta^2}{\rmin}\right).
\]
Recalling the result of~\proofstep{1} (namely~\eqref{eq:unconstrained}) we conclude that, in order to establish the~$\rmin$-dependent part of~\eqref{eq:lower-fixed-design}, it suffices to lower-bound the {\em largest} of the two Bayes risks associated with~$\pi$ and~$\wt\pi$ as
\begin{align}
\max \Big\{ 
&\underset{\wh T}{\inf} \left\{ \E_{\pi} \left[\Prob_{\cH_0}[\wh T = 1|X,X'] + \Prob_{\cH_1}[\wh T = 0|X,X'] \right] \right\}, \notag \\
&\underset{\wh T}{\inf} \left\{ \E_{\wt\pi} \left[\Prob_{\cH_0}[\wh T = 1|X,X'] + \Prob_{\cH_1}[\wh T = 0|X,X'] \right] \right\} 
\Big\} \geq 3\exp \left(-\frac{c'n^2\Delta^{2}}{\rmin}\right).
\label{eq:bayes-risk-alternative}
\end{align}

In the next step we prove~\eqref{eq:bayes-risk-alternative}; it is in this step that we use condition~\eqref{eq:evals-condition} on the design matrices.

\proofstep{3}: \textbf{\em Lower-bounding the Bayes risk.}
For~$\pi$ chosen in~\proofstep{2}, the Neyman-Pearson lemma gives
\[
\begin{aligned}
&\underset{\wh T}{\inf} \E_{\pi} \left[\Prob_{\cH_0}[\wh T = 1|X,X'] + \Prob_{\cH_1}[ \wh T = 0|X,X'] \right]  
\ge 
\E_{\pi} \left[\Prob_{\cH_0}[\wh T_{\MLR} = 1|X,X'] + \Prob_{\cH_1}[\wh T_{\MLR} = 0|X,X'] \right]
\end{aligned}
\]
where~$\wh T_{\MLR}$ is the {\em marginal likelihood ratio} test associated with~$\pi$ and given by
\[
\wh T_{\MLR} = \ind \left\{  Q_{\cH_0}(Y,Y' | X,X') \leq Q_{\cH_1}(Y,Y' | X,X') \right\};
\] 
here~$Q_{\cH_0}(Y,Y' | X,X')$ (resp.,~$Q_{\cH_1}(Y,Y' | X, X')$) is the marginal (over~$\bar \theta \sim \pi$) likelihood of~$(Y,Y')$ under~$\cH_0$ (resp.,~$\cH_1$) conditionally on~$(X,X')$.
Now, recalling that~$\cH_1$ reads~$(\theta^*, \bar\theta) = (\theta_1, \theta_0)$ we get
\[
\begin{aligned}
&Q_{\cH_1}^2(Y,Y' | X,X') 
= 
\frac{
\exp\big(-\|Y' - X'\theta^{*}\|^2 - (Y- X\theta^*)^\top\left(\Id + n\rho\bPi_{X_\cI}\right)^{-1}(Y - X\theta^*) \big)
}{
C_n \det\left(\Id + n \rho \bPi_{X_\cI}\right)
},
\end{aligned}
\]
where~$\bPi_{X_\cI}$ is given in~\eqref{eq:projector-I}, and~$C_n$ is the normalization constant depending solely on~$n$.
Indeed, conditionally on~$(X,X')$ under~$\cH_1$ we have that~$Y' \sim \cN(X\bar\theta,\Id)$,~$X\bar\theta \sim \cN(X\theta^*,n\rho \bPi_{X_\cI})$
by~\eqref{eq:prior-prediction-residual}, and~$Y - X\bar\theta \sim \cN(0,\Id)$ independently; hence,~$Y \sim \cN(X\theta^*, \Id + n\rho \bPi_{X_\cI})$ marginally over~$\bar\theta$.
Similarly,
\[
Q_{\cH_0}^2(Y,Y' | X,X') 
= \frac{
\exp\big(-\|Y-X\theta^{*}\|^2-(Y'-X'\theta^*)^\top(\Id + n\rho\bOmega_{X,X'})^{-1}(Y' - X' \theta^*) \big)
}{
C_n \det (\Id + n \rho \bOmega_{X,X'})
},
\]
where
\[
\bOmega_{X,X'} := X'_{\cI} (X_{\cI}^\top X_{\cI}^{\vphantom\top})^{-1}(X'_\cI)^\top 
= \frac{1}{n} X'_{\cI} \wh\bSigma_{0,\II}^{-1} (X'_\cI)^\top.
\]
Hence,~$\wh T_{\MLR}$ writes
\[
\wh T_{\MLR} = \ind\{  \|\bA^{1/2}(Y- X\theta^*)\|^2 - \log\det(\Id-\bA) 
\geq \|\bB^{1/2}(Y'- X'\theta^*)\|^2 - \log\det(\Id-\bB)\},
\]
where 
\[
\bA := \Id - (\Id + n\rho\bPi_{X_I})^{-1} 
\quad \text{and} \quad 
\bB := \Id - (\Id + n\rho\bOmega_{X,X'})^{-1}.
\]
Let us now lower-bound the marginal over~$\pi$ type I error~$\E_{\pi} \Prob_{\cH_0}[\wh T_{\MLR} = 1|X,X']$ (clearly this would also give a lower bound for the expected sum of errors).
To this end, observe that
\[
\bA = \frac{n\rho}{1+ n\rho} \bPi_{X_\cI}
\] 
and, under~$\cH_0$ and conditionally on~$(X,X')$, we have
\[
\|\bA^{1/2}(Y- X\theta^*)\|^2 \sim \frac{n\rho}{1+n\rho}\|\bPi_{X_\cI} \xi\|^2  \quad \text{and} \quad \|\bB^{1/2}(Y'- X'\theta^*)\|^2 \sim n\rho \| \bOmega_{X,X'}^{1/2} \zeta \|^2,
\]
where~$\xi := Y-X\theta_0 \sim \cN(0,\Id)$, and~$\zeta$ is the whitening of~$Y' - X'\theta_0 = \xi' + X'(\theta_1 - \theta_0)$; thus,~$\zeta,\xi$ are independent and have the same distribution.
By rotational invariance, we may assume w.l.o.g.~that
\[
n\rho \bOmega_{X,X'} = \sum_{j=1}^{\lfloor \rmin /4 \rfloor} \mu_j \ue_j^{\vphantom\top} \ue_j^\top,
\]
where~$\ue_j$ is the~$j$-th canonical vector of~$\R^n$; here we used that~$\rank(\bOmega_{X,X'}) \le \lfloor \rmin /4 \rfloor$. 
As a result, 
\[
\E_\pi\Prob_{\cH_0}\left[\wh T_{\MLR} = 1|X,X'\right] 
=
\Prob\left[ 
\|\bPi_{X_{\cI}}{\xi}\|^2 
\geq 
\frac{1+ n\rho}{n\rho}\left(\rmin \log(1+n\rho)/4 + \sum_{j=1}^{\lfloor \rmin /4 \rfloor} \mu_j\zeta_j^2 - \log(1+\mu_j) \right) \right].
\]
Now, by the standard deviation bound for the generalized chi-squared distribution (\cite[Lemma~1]{laurent2000adaptive})
\[
\sum_{j=1}^{\lfloor\rmin/4 \rfloor} \mu_j\zeta_j^2 \leq \sum_{j=1}^{\lfloor\rmin/4 \rfloor} \mu_j + 2\sqrt{\sum_{j=1}^{\lfloor\rmin/4 \rfloor}\mu_j^2}
\] 
with probability~$c > 0$. 
Moreover,~$\mu_j - \log(1+\mu_j) \leq \mu_j^2$. Hence
\[
\E_\pi \Prob_{\cH_0}\left[\wh T_{\MLR} = 1|X,X'\right] 
\geq 
c\Prob\left[ \|\bPi_{X_\cI} \xi\|^2 \geq \frac{1+n\rho}{n\rho} \left( \rmin \log(1+n\rho)/4 + \sum_{j=1}^{\rmin/4} \mu_j^2 + 2\sqrt{\sum_{j=1}^{\rmin/4}\mu_j^2}\right) \right].
\]
Moreover,~$\log(1+n\rho) - \frac{n\rho}{1+n\rho} \leq (n\rho)^2$ and~$n\rho \le 1$ (cf.~\eqref{eq:Delta-range-for-complex-lower}), so we arrive at
\begin{equation}
\label{eq:before-commutativity-argument}
\E_\pi\Prob_{\cH_0}\left[\wh T_{\MLR} = 1|X,X'\right] 
\geq 
c\Prob\left[ \|\bPi_{X_\cI} \xi\|^2   - \rmin /4 \geq \frac{2}{n\rho}\left( \sum_{j=1}^{\rmin /4}\mu_j^2+2\sqrt{\sum_{j=1}^{\rmin /4}\mu_j^2}\right) +  \rmin n\rho/2 \right].
\end{equation}
Now, observe that  that~$\|\bOmega_{X,X'}\| \le 49$. 
Indeed, due to~\eqref{eq:evals-condition} we have that

\[
\|\bOmega_{X,X'}\| \le \frac{\lambda_{\max}(\wh \Sigma_{1,\II})}{\lambda_{\min}(\wh \Sigma_{0,\II})} \leq  49,
\] 
and thus~$\mu_j \le 49 n\rho$ for all~$j$. 
Returning to~\eqref{eq:before-commutativity-argument}, and recalling the definition of~$\rho$ (cf.~\eqref{eq:prior-definition}) we get
\[
\begin{aligned}
\E_\pi\Prob_{\cH_0}\left[\wh T_{\MLR} = 1|X,X'\right] 
&\geq 
c\Prob\left[ \|\bPi_{X_\cI} \xi\|^2 - \rmin /4 \geq 100\sqrt{\rmin}+1250\rmin n\rho \right] \\
&=
c\Prob\left[ \|\bPi_{X_\cI} \xi\|^2-\rmin /4\geq 100\sqrt{\rmin}+33 \cdot 10^4 n \Delta \right].
\end{aligned}
\]
%
We are now in the position to apply a lower bound for the right tail of the~$\chi^2_s$ distribution with~$s \ge 2$ degrees of freedom (\cite[Prop.3.1]{inglot}):\footnotemark
\begin{equation*}
\Prob \left[ \chi_s^2 - s \geq u \right] \geq \frac{c(u+s)}{\sqrt{s}(u+\sqrt{s})} \exp \left(-\frac{Cu^2}{s}\right),
\end{equation*}
This bound, when applied with~$s = \rmin /4$ and~$u = 100\sqrt{\rmin}+33 \cdot 10^4 n \Delta $,  suffices for our purposes. 
Indeed, recall that we are in the regime~$n \rho \le 1$, i.e.,~$n\Delta = O(\rmin)$. Thus, we need to apply the bound in the range~$c_1 \sqrt{s} \le u \le c_2 s$ where it becomes
\begin{equation}
\label{eq:chi-squared-lower-estimates}
\begin{aligned}
\Prob \left[ \chi_s^2 - s \geq u \right] \geq 
\frac{c\sqrt{s}}{u} \exp \left(-\frac{Cu^2}{s}\right) 
\ge c\exp \left(-\frac{2Cu^2}{s}\right) 
&\ge c\exp \left(-\frac{C'(\sqrt{\rmin} + n\Delta)^2}{\rmin}\right) \\
&\ge c' \exp \left(-\frac{C'' n^2 \Delta^2}{\rmin}\right).
\end{aligned}
\end{equation}
This proves the~$\rmin$-dependent part of~\eqref{eq:lower-fixed-design}. 
Recalling the discussion in Sec.~\ref{sec:lower-fixed-design-reduction} we see that, in order to conclude the proof of Theorem~\ref{th:lower-gaussian}, it only remains to show~\eqref{eq:evals-condition-prob}; this is our goal in the next section.

\footnotetext{One can verify that this bound is tight (matching~\eqref{eq:tail-bounds-chi-sq-right}) when~$u = O(k)$, which is our case; however, it proves to be loose for larger deviations.
Meanwhile, recent work~\cite[Corollary~3]{zhang2018non} establishes the sharpness of~\eqref{eq:tail-bounds} in all regimes.}

\subsection{Proof of~\eqref{eq:evals-condition-prob}}
\label{sec:lower-orths}
We shall rely on a classical result about random Gaussian matrices. 
For integer~$p < n$, let~$X \in \R^{n \times p}$ have i.i.d.~$\cN(0,1)$ entries. 
Then it is known (\cite[Theorem II.13]{davidson2001local}) that, for any $t>0$,
\begin{equation*}
\max \Big\{ \,
\Prob \left\{ \lambda_{\min}(X) \leq (1-t)\sqrt{n} - \sqrt{p} \right\}, \;\;
\Prob \left\{ \lambda_{\max}(X) \geq (1+t)\sqrt{n} + \sqrt{p} \right\} 
\Big\} 
\leq \exp \left(-nt^2/2 \right).
\end{equation*}
Noting that $|\cI| = \lfloor \rmin /4 \rfloor \leq n/4$ and taking $t= 1/4$, we conclude that
\[
\frac{1}{16} \leq  \lambda_{\min} (\wh \bSigma_{0,\II}) , \lambda_{\max} (\wh \bSigma_{1,\II}) \leq \frac{49}{16}
\]
with probability at least $1-2e^{-n/32} \geq 1/8$ as we are assuming $n \geq \rmin \geq 28$.


\section{Proofs for asymptotic results}
\label{app:upper-general}

\subsection{Proof of Proposition~\ref{prop:asymp}}
For convenience we define~$\Mopt_k := \M_k(\theta_k)$ for~$k \in \{0,1\}$.
As in the proof of Theorem~\ref{th:upper-random-design}, we decompose the statistic whose sign is examined in~\eqref{eq:test-asymp} as~$\wh S = \wh S_0 - \wh S_1$. 
W.l.o.g., we analyze the type II error probability~$\Prob_{\cH_1}[\wh T = 0];$ thus,~$\theta^* = \theta_1$.

\proofstep{1}.
Recall that~$\Mopt_1 := \M_1(\theta_1)$. 
We first study
\[
\begin{aligned}
\wh S_1 := n_1 \| \wh \H_1(\theta_1)^{\dagger/2} \nabla \wh L_1(\theta_1) \|^2 - \Tr[\Mopt_1]
\end{aligned}
\]
and the complementary to it term within~$\wh S_0$. 
Let us first define random variables
\[
Z(\theta) := \nabla \wh L_0(\theta) - \nabla L_0(\theta), \quad
Z'(\theta) := \nabla \wh L_1(\theta) - \nabla L_1(\theta);
\]
note also that~$\nabla L_0(\theta_0) = \nabla L_1(\theta_1) = 0$.
By the central limit theorem (CLT) we have the following convergence in distribution as~$n_0, n_1 \to \infty$: 
\[
\sqrt{n_0} Z(\theta) 
\weakto \cN(0,\G_0(\theta)), \quad 
\sqrt{n_1} Z'(\theta) 
\weakto \cN(0,\G_1(\theta)),
\]
with~$\G_0(\theta),\G_1(\theta)$ defined in~\eqref{eq:matrices-population}.
Define~$\chi^2_{\M}$ as the generalized chi-square law -- precisely, the law of~$\|\xi\|^2$ with~$\xi \sim \cN(0,\M)$ for PSD matrix~$\M$.
Recall that~$\E[\chi^2_{\M}] = \Tr[\M]$, and~\cite[Lemma~1]{laurent2000adaptive} gives
\begin{subequations}
\label{eq:generalized-chi-square}
\begin{align}
\label{eq:generalized-chi-square-right}
\Prob[\chi^2_{\M} - \Tr[\M] \ge t] &\le\exp \left(-c\min\left\{\frac{t}{\|\M\|},\frac{t^2}{\Tr[\M^2]}\right\}\right),\\
\label{eq:generalized-chi-square-left}
\Prob[\chi^2_{\M} - \Tr[\M] \le -t] &\le\exp \left(-\frac{ct^2}{\Tr[\M^2]}\right).
\end{align}
\end{subequations}
At fixed~$\theta$ we have~$n_0 \| \H_0(\theta)^{\dagger/2} Z(\theta) \|^2 \weakto \chi^2_{\M_0(\theta)}$ and~$n_1 \| \H_1(\theta)^{\dagger/2} Z(\theta) \|^2 \weakto \chi^2_{\M_1(\theta)},$ that is,
\[
n_1 \| \H_1(\theta_1)^{\dagger/2} \nabla \wh L_1(\theta_1) \|^2 \weakto \chi^2_{\Mopt_1}
\]
since~$\nabla L_{1}(\theta_1) = 0$.
Clearly, the fully empirical counterpart of this quantity~$n_1 \| \wh \H_1(\theta_1)^{\dagger/2} \nabla \wh L_1(\theta_1) \|^2$ has
the same asymptotic distribution. 
Indeed, at any fixed~$\theta$ one has~$\wh \H_k(\theta) \to \H_k(\theta)$ in probability as~$n_k \to \infty$ by the law of large numbers; in particular,~$\wh \H_1(\theta_1) \to \H_1(\theta_1)$. 
Thus, by the matrix version of Slutsky's theorem~(\cite{van2000asymptotic}) we have that~$\wh\H_1(\theta_1)^{\dagger/2} \nabla \wh L_1(\theta_1) \weakto \cN(0, \Mopt_1)$; as a result, indeed,
\[
n_1 \| \wh \H_1(\theta_1)^{\dagger/2} \nabla \wh L_1(\theta_1) \|^2 \weakto \chi^2_{\Mopt_1}.
\]
Using tail bounds~\eqref{eq:generalized-chi-square-right}--\eqref{eq:generalized-chi-square-left}, for any fixed~$t > 0$  we have
$
\Prob [\wh S_1 \ge t] 
\to 
\Prob[\chi^2_{\Mopt_1} - \Tr[\Mopt_1] \ge t]
$
as~$n_1 \to \infty$. 
Now, recall that we are in the regime~ (in particular,~$n_k \bar\Delta_k \to \lambda_k$ for~$ k \in \{0,1\}$.
Thus, choosing~$t = \lambda_0/3,$ and applying~\eqref{eq:generalized-chi-square-right}, we arrive at 
\begin{equation}
\label{eq:limit-1}
\lim_{[...]} \Prob \left[ \wh S_1 \ge \frac{n_0\bar\Delta_0}{3} \right] \le 
\exp\left(-{c \lambda_0} \min \left\{ \frac{1}{\|\Mopt_1\|}, \frac{\lambda_0}{\Tr[(\Mopt_1)^2]} \right\} \right),
\end{equation}
where the limit is in the regime specified in the premise of the theorem (we will keep this notation).

\proofstep{2}.
We now consider~$\wh S_0$ which we decompose as 
\[
\begin{aligned}
\wh S_0 
&:= n_0 \|\wh\H_0(\theta_1)^{\dagger/2} \nabla \wh L_0(\theta_1) \|^2 - \Tr[\M_0(\theta_1)] \\
&=
n_0 \| \wh\H_0(\theta_1)^{\dagger/2} Z(\theta_1) \|^2 -\Tr[\M_0(\theta_1)]  + n_0\bar\Delta_0 
+ 2 n_0 \langle \wh\H_0(\theta_1)^{\dagger/2} Z(\theta_1), \wh\H_0(\theta_1)^{\dagger/2} \nabla L_0(\theta_1) \rangle.
\end{aligned}
\]
We have~$n_0 \| \wh\H_0(\theta_1)^{\dagger/2} Z(\theta_1) \|^2 - \Tr[\M_0(\theta_1)] \weakto \chi^2_{\M_0(\theta_1)}$ by CLT combined with Slutsky's theorem (cf.~\proofstep{1}). Hence, using~\eqref{eq:generalized-chi-square-left}, we have
\begin{equation}
\label{eq:limit-2}
\begin{aligned}
\lim_{[...]}
\Prob \left[n_0 \| \wh\H_0(\theta_1)^{\dagger/2} Z(\theta_1) \|^2 - \Tr[\M_0(\theta_1)] \le - \frac{n_0\bar\Delta_0}{3}\right] 
&\le \exp\left(-\frac{c \lambda_0^2}{\Tr[\M_0^2(\theta_1)]} \right).
\end{aligned}
\end{equation}
Now, by CLT combined with Slutsky's theorem,~$\sqrt{n_0} \wh\H_0(\theta_1)^{\dagger/2} Z(\theta_1) \weakto \cN(0,\M_0(\theta_1))$.
Due to that, and since~$\|\wh\H_0(\theta_1)^{\dagger/2} \nabla L_0(\theta_1)\|^2 / \bar\Delta_0 \to 1$ in probability, we have that
\[
\sqrt{n_0} \langle \wh\H_0(\theta_1)^{\dagger/2} Z(\theta_1), \wh\H_0(\theta_1)^{\dagger/2} \nabla L_0(\theta_1) \rangle \weakto \cN(0,v_0) \;\; \text{with} \;\; v_0 \le \|\M_0(\theta_1)\|\bar\Delta_0.
\]
Therefore by~\eqref{eq:tail-bounds-normal}, 
\begin{equation}
\label{eq:limit-3}
\lim_{[...]} 
	\Prob\left[ 2n_0 \langle \wh\H_0(\theta_1)^{\dagger/2} Z(\theta_1), \wh\H_0(\theta_1)^{\dagger/2} \nabla L_0(\theta_1) \rangle 
\le - \frac{n_0\bar\Delta_0}{3} \right] 
=\exp\left(\frac{-\lambda_0}{36\|\M_0(\theta_1)\|}\right).
\end{equation}
Finally, combining~\eqref{eq:limit-1}--\eqref{eq:limit-3} through the union bound, and observing that~$\Tr[\M_0(\theta_1)] / \Tr[\Mopt_0] \to 1$ and~$\|\M_0(\theta_1) \| / \|\Mopt_0\| \to 1$ since~$\M_0(\cdot)$ is continuous at~$\theta_0$ and~$\bar\Delta_0 \to 0$, we arrive at the claimed bound: 
\begin{equation}
\lim_{[...]} \Prob_{\cH_1}[\wh T = 0] 
\le 
C \exp\bigg(-c\lambda_0 \min \bigg\{ \frac{1}{\|\Mopt_0\|}, \; \frac{1}{\|\Mopt_1\|}, \; \frac{\lambda_0}{\Tr[(\Mopt_0)^2]}, \; \frac{\lambda_0}{\Tr[(\Mopt_1)^2]} \bigg\} \bigg).
\end{equation}
The type I error bound follows by symmetry.
\qed

\subsection{Proof of Theorem~\ref{th:asymp-adapt}}

W.l.o.g., we again focus on the type II error. 
Following the argument in the main text after the theorem statement (modulo the switch of indices due to the different hypothesis), it remains to verify that
\[
\wh \cT_1 \weakto \chi_{\M_1(\theta_1)}^2, 
\quad 
\wh \cT_0 \weakto \chi_{\M_0(\theta_1)}^2.
\] 
Then the result will follow by applying the generalized chi-squared bounds~\eqref{eq:generalized-chi-square} with~$t = \lambda_0/3$ and replacing~$\M_0(\theta_1)$ with~$\M_0(\theta_0)$ by continuity. To verify the weak convergence, by using~\eqref{eq:trace-estimate} we write
\[
\wh\cT_0 
:= 
\frac{n_0}{2} \big\| \wh\H_0(\theta_1)^{{\dagger}/{2}}  \big[ Z(\theta_1) - \wt Z(\theta_1) \big] \big\|^2,
	\quad (k \in \{0,1\}),
\]
where~$\wt Z(\theta)$ is an independent copy of~$Z(\theta)$ from the previous proof. 
As before,~$\wh \H_0(\theta_1) \to \H_0(\theta_1)$ in probability by the law of large numbers, whereas~$\sqrt{n_0} Z(\theta) \weakto \cN(0,\G_0(\theta))$ by CLT, thus
\[
Z(\theta) - \wt Z(\theta) \weakto \cN(0,2\G_0(\theta))
\] 
by independence.  
By Slutsky's theorem, this results in~$\wh \cT_0 \weakto \chi_{\M_0(\theta_1)}^2$. 
The case of~$\wh \cT_1$ is similar, and the result follows. \qed

\section{Proofs for generalized linear models}
\label{app:proofs-glm}

\subsection{Proof of Theorem~\ref{th:glm-fixed-prob}}
W.l.o.g.~we analyze the type II error~$\Prob_{\cH_1}[\wh T = 0]$ (so that~$\theta^* = \theta_1$) and use the abridged notation:
\[
\begin{aligned}
(x_i, y_i) := (x_i^{(0)}, y_i^{(0)}), \quad (x_i', y_i') := (x_i^{(1)}, y_i^{(1)}), \\
\wh\Delta := \wh\Delta_0(x), \quad 
\Delta := \Delta_0, \quad 
\nu := \nu_0, \quad 
\kappa := \kappa_0.
\end{aligned}
\]
\proofstep{0}.
We first focus on the ``shifted'' term in~\eqref{eq:test-glm},
\[
\wh S_0 := \sum_{i=1}^{n} \frac{(a'(\eta_{i}(\theta_1)) - y_i)^2}{a''(\eta_{i}(\theta_1))} - n \nu(\theta_1),
\]
where~$\eta_i(\theta) := \eta_{x_i}(\theta)$. 
By adding and subtracting~$\E[y_i|x_i]$ under the square, we get
\begin{align}
\wh S_0
&= 
   \sum_{i=1}^{n}  \frac{(y_i - \E[y_i|x_i])^2}{a''(\eta_{i}(\theta_1))} -  \nu(\theta_1) 
	+  
	\frac{(a'(\eta_{i}(\theta_1)) - \E[y_i|x_i])^2}{a''(\eta_{i}(\theta_1))} 
    - 	
    \frac{2(a'(\eta_{i}(\theta_1)) - \E[y_i|x_i])(y_i - \E[y_i|x_i])}{a''(\eta_{i}(\theta_1))} 
    \notag\\
&=  \wh Q_0 + \wh R_0 + \wh U_0.
\label{eq:glm-decomposition}
\end{align}
The term~$\wh Q_0$ corresponds to the centered~$\chi^2$-type statistic.
The term~$\wh R_0$ is non-random when conditioned on~$X = [x_1; ...; x_n],$ and incorporates the offset of the statistic due to~$\Delta$. 
The cross-term~$\wh U_0$ gives smaller-order fluctuations. 
Representing the ``unshifted'' term~$\wh S_1 = \wh S_0 - \wh S$ in a similar manner, where~$\wh S$ is the whole statistic whose sign is examined in~\eqref{eq:test-glm}, we decompose~$\wh S$ as follows:
\[
\wh S
= \wh S_0 - \wh S_1
= [\wh Q_0 - \wh Q_1] + [\wh R_0 - \wh R_1] + [\wh U_1 - \wh U_0].
\]
Here~$\wh Q_1, \wh R_1, \wh U_1$ are the counterparts of~$\wh Q_0, \wh R_0, \wh U_0$ which we do not write down explicitly here.

\proofstep{1}.
Let us define
\[
\omega_i := \frac{y_i - \E[y_i|x_i]}{\sqrt{a''(\eta_{i}(\theta_1))}}.
\]
For~$\wh Q_0 = \sum_{i=1}^n \omega_i^2 - n \nu(\theta_1)$, we have~$\E[\wh Q_0] = 0$ and
$
\Var[\wh Q_0] = \sum_{i=1}^n \Var[\omega_i^2] = n [\kappa(\theta_1) - \nu^2(\theta_1)],
$
where we used independence of~$\omega_i$'s. 
Hence, by Chebyshev's inequality
\[
\begin{aligned}
\Prob\left[  |\wh Q_0 | \le \frac{n \Delta}{c} \right] 
	&\ge 1-\frac{c^2 [\kappa(\theta_1) - \nu^2(\theta_1)]}{n \Delta^2},
\end{aligned}
\]
where~$c > 0$ will be chosen later. In particular, the right-hand side is~$\ge 39/40$ whenever
\[
n\Delta^2 \ge {40 c^2\kappa(\theta_1)}. 
\]
As a result, repeating the same analysis for~$\wh Q_1$, we get
\[
\Prob\left[  |\wh Q_0 - \wh Q_1| \le \frac{2n \Delta}{c} \right] \ge \frac{19}{20}
\]
whenever 
\begin{equation}
\label{eq:n-kurtosis-apx}
n \ge \frac{40 c^2 \max\left\{\kappa_1(\theta_1), \kappa_0(\theta_1)\right\}}{\Delta^2}.
\end{equation}
On the other hand, we verify by simple calculations that~\eqref{eq:n-kurtosis-apx} with~$c = 120$ is indeed implied by~\eqref{eq:condition_n}.\\

\proofstep{2}. 
Our next goal is to similarly bound~$\wh R_0 - \wh R_1$ from below by~$\Omega(n \Delta^2)$.
To this end, we have
\[
\begin{aligned}
\wh R_0 
&= \sum_{i=1}^{n} \wh \Delta(x_i) + \frac{(a'(\eta_{i}(\theta_0)) - \E[y_i|x_i])^2}{a''(\eta_{i}(\theta_1))}
	+ 2 \cdot \frac{a'(\eta_{i}(\theta_1)) - a'(\eta_{i}(\theta_0))}{\sqrt{a''(\eta_{i}(\theta_1))}} \cdot \frac{a'(\eta_{i}(\theta_0)) - \E[y_i|x_i]}{\sqrt{a''(\eta_{i}(\theta_1))}} \\
&\ge \sum_{i=1}^{n} \frac{\wh \Delta(x_i)}{2} - \frac{(a'(\eta_{i}(\theta_0)) - \E[y_i|x_i])^2}{a''(\eta_{i}(\theta_1))} 
= \sum_{i=1}^{n} \frac{\wh \Delta(x_i)}{2} - \Bias_0(\theta_1|x_i),
\end{aligned}
\]
where we used the expressions~\eqref{def:Delta-glm} and~\eqref{eq:misspec-bias}:
\[
\wh\Delta(x) := \frac{[a'(\eta_x(\theta_0)) - a'(\eta_x(\theta_{1}))]^2}{a''(\eta_x(\theta_{1}))},
\quad
\Bias_0(\theta|x) := \frac{(a'(\eta_{x}(\theta_0) - \E[y|x])^2}{a''(\eta_{x}(\theta))}
\]
and in the second line we used that~$-2ab \le a^2/2 + 2b^2$. 
Repeating the analysis for~$\wh R_1$, we arrive at
\[
\wh R_0  - \wh R_1 
\ge \sum_{i=1}^{n} \frac{\wh \Delta(x_i)}{2} - \Bias_0(\theta_1|x_i) - \Bias_1(\theta_1|x_i').
\]
Now observe that, due to the small misspecification bias assumption~\eqref{eq:small-bias-ii},
\[
\E\left[\sum_{i=1}^{n} \Bias_0^2(\theta_1|x_i) + \Bias_1(\theta_1|x_i') \right] 
\le n[\Bias_0(\theta_1) + \Bias_1(\theta_1)] 
\le \frac{n\Delta}{400}.
\]
Hence, by Markov's inequality, 
\[
\Prob\left[\sum_{i=1}^{n} \Bias_0(\theta_1|x_i) + \Bias_1(\theta_1|x_i') \le \frac{n\Delta}{40} \right] \ge \frac{9}{10}.
\]
On the other hand, by the Paley-Zygmund inequality~(\cite{paley1932note}), for arbitrary~$t \in [0,1]$ we have
\[
\Prob\left[\sum_{i=1}^{n} \wh \Delta(x_i) \ge t {n \Delta} \right] 
\ge (1-t)^2 \frac{n^2 \Delta^2}{\E\Big[\big(\sum_{i=1}^{n} {\wh \Delta(x_i)}\big)^2\Big]}.
\]
The denominator satisfies
\[
\E\left[\left(\sum_{i=1}^{n} {\wh \Delta(x_i)}\right)^2\right] 
= 
n \E[\wh \Delta^2] + n(n-1) \Delta^2 
= 
n^2 \Delta^2 + n \Var[\wh \Delta] 
\le 
\frac{81}{80} n^2 \Delta^2.
\]
where in the end we used that~$n \Delta^2 \ge 80 \Var[\wh \Delta]$, cf.~\eqref{eq:condition_n}. 
Thus, by choosing~$t = 1/10$ we arrive at
\[
\Prob\left[\sum_{i=1}^{n} \wh \Delta(x_i) \ge \frac{n \Delta}{10} \right]  \ge \frac{8}{10}.
\]
Combining this with the previous result through the union bound, we finally arrive at the bound
\[
\Prob\left[\wh R_0  - \wh R_1 \ge \frac{n \Delta}{40} \right]  \ge \frac{7}{10}.
\]

\proofstep{3}.
It remains to upper-bound the cross-term~$\wh U_0-\wh U_1$. 
First, note that condition~\eqref{eq:condition_n} implies
\begin{equation}
\label{eq:n-lower-bound-for-last-step}
\begin{aligned}
\sqrt{n}\Delta
&\geq 40 c\sqrt{20\max \left\{\kappa_{0}(\theta_{1}), \kappa_1(\theta_{1}), 4\E[\wh \Delta^2]\right\}}\\
\end{aligned}
\end{equation}
with~$c = 120$ -- i.e., with the same constant as in step~\proofstep{1}.
Now, let us consider the following event:
\[ 
\mathcal{E} = \left\{ \sum_{i=1}^n \hat{\Delta}^2(x_i) + \Bias_0^2(\theta_1|x_i) + \Bias_1^2(\theta_1|x_i') \leq {10n\E[\wh\Delta^2]} \right\}.
\]
Using Markov's inequality and~\eqref{eq:small-bias-ii}, one can verify that~$\Prob\left[ \mathcal{E} \right] \geq 49/50$.
Moreover, introducing
\[
\alpha^2_0(\theta_1|x) := \frac{(a'(\eta(\theta_1)) - \E[y|x])^4}{a''(\eta(\theta_1))^2} \leq 4 \hat{\Delta}^2(x) + 4\Bias_0^2(\theta_1|x),
\]
we have that
\begin{align*}
&\Prob\left[  |\wh U_0 - \wh U_1| \le \frac{n \Delta}{c} \right]  \\ 
&\ge\Prob\left[ |\wh U_0 | + |\wh U_1 | \le \frac{n \Delta}{c} \right]  \\
&\ge
	\Prob\left[ |\wh U_0| \le \frac{n \Delta}{2c} \right] \Prob\left[ |\wh U_1 | \le \frac{n \Delta}{2c} \right] \\
&\ge 
	\left(\frac{49}{50} - \Prob\left[ |\wh U_0 |\ind_{\mathcal{E}} \ge \frac{n \Delta}{2c} \right] \right)\left(\frac{49}{50} - \Prob\left[ |\wh U_1|\ind_{\mathcal{E}} \ge \frac{n \Delta}{2c} \right] \right)\\
&\stackrel{(i)}{\ge} 
	\left(\frac{49}{50}-\frac{16c^2 \E\left[\ind_{\mathcal{E}}\sum_{i=1}^n \alpha_{0}(\theta_1|x_i) \nu_0(\theta_1|x_i)\right]}{n^2 \Delta^2}\right)
	\cdot \left(\frac{49}{50}-\frac{16c^2 \E\left[\ind_{\mathcal{E}}\sum_{i=1}^n \beta_{1}(\theta_1|x_i') \nu_1(\theta_1|x_i')\right]}{n^2 \Delta^2}\right)\\
&\stackrel{(ii)}{\ge} 
	\left( \frac{49}{50}-\frac{16c^2 \sqrt{\kappa_{0}(\theta_1) \E\left[\ind_{\mathcal{E}}\sum_{i=1}^n \alpha^2_{0}(\theta_1|x_i)\right] }}{n^{3/2} \Delta^2}\right)
	\cdot \left( \frac{49}{50}-\frac{16c^2 \sqrt{\kappa_{1}(\theta_1) \E\left[\ind_{\mathcal{E}}\sum_{i=1}^n \beta^2_{1}(\theta_1|x_i')\right] }}{n^{3/2} \Delta^2}\right) \\
&\stackrel{(iii)}{\ge} 
	\left(\frac{49}{50}-\frac{64c^2\sqrt{10\kappa_1(\theta_1)\E[\wh\Delta^2]})}{n\Delta^2} \right) \cdot \left(\frac{49}{50}-\frac{16c^2\sqrt{10\kappa_0(\theta_1)\E[\wh\Delta^2]}}{n\Delta^2} \right)
\stackrel{(iv)}{\ge} \frac{19}{20}. 
\end{align*}
Here we proceeded as follows:
\begin{itemize}
\item
in~$(i)$ we recognized~$\alpha^2_0(\theta_1|x)$ and~$\Bias_1^2(\theta_1|x')$ and then used Markov's inequality in the form
\[
\Prob[|X| \ge \varepsilon] \le \E[X^2]/\varepsilon^2;
\]
\item
in~$(ii)$ we used the Cauchy-Schwarz inequality twice: first on the sum under the expectation to obtain the product of two random variables, and then on the expectation of the product;
\item
in (iii) we used the definition of the event~$\mathcal{E}$;
\item
in (iv) we used~\eqref{eq:condition_n} and performed a direct calculation. 
\end{itemize}
By putting~$c = 120$ and combining the results of~\proofstep{1}--\proofstep{3} via the union bound, we bound the probability of type II error by~$2/5$ as required.
The bound for the type II error follows by symmetry.
\qed

\subsection{Proof of Corollary \ref{cor:voting}}
By Theorem \ref{th:glm-fixed-prob}, we have that~$\Prob_{\cH_0}[\wh T_j = 1] \le 2/5$ for all $j \in [b]$. 
Hence
\begin{align*} 
\Prob_{\cH_0}[\wh T = 1] 
&= 
\Prob_{\cH_0}\left[ \sum_{j \in [b]} \hat{T}_j \geq \frac{b}{2} \right]
\leq 
\Prob_{\cH_0}\left[ \sum_{j \in [b]} \hat{T}_j - \E_0[\hat{T}_j] \geq \frac{b}{10}\right]
\leq  
e^{-cb},
\end{align*}
for some constant $c>0$, where the final step is by Hoeffding's inequality. Plugging in the bound~\eqref{eq:choose_k} for~$m = \lceil n/b \rceil$, we verify the claim.
The matching inequality for $\Prob_{\cH_1}[\wh T = 0]$ follows by symmetry.
\qed


\subsection{Proof of Proposition~\ref{prop:glm-adapt}}
When analyzing the error probabilities for test~\eqref{eq:test-glm-adapt}, we can decompose the statistic almost in the same way as for~\eqref{eq:test-glm}, with two additional terms corresponding to the fluctuations of~$\wh \nu_0(\theta^*)$ and~$\wh \nu_1(\theta^*)$ around their expectations. More precisely, representing as~$\wh S = \wh S_0 - \wh S_1$ the statistic whose sign is examined in~\eqref{eq:test-glm-adapt} (in the same way as in the proof of Theorem~\ref{th:glm-fixed-prob}), we decompose~$\wh S_0$ as
\[
\wh S_0 = \wh Q_0 + \wh R_0 - \wh U_0 - \wh V_0 \;\; \text{with} \;\; \wh V_0 := n_0 [\wh \nu_0(\theta^*) - \nu_0(\theta^*)],
\]
where the first three terms are the same as in~\eqref{eq:glm-decomposition}. Thus, defining~$\wh V_1 := n_1 [\wh\nu_1(\theta^*) - \nu_1(\theta^*)]$, we have
\[
\wh S
= [\wh Q_0 - \wh Q_1] + [\wh R_0 - \wh R_1] + [\wh U_1 - \wh U_0] + [\wh V_1 - \wh V_0]. 
\]
We now show that the additional fluctuation term~$|\wh V_1 - \wh V_0|$ admits essentially the same bound as~$|\wh U_1 - \wh U_0|$. 
Indeed, in the notation of the proof of Theorem~\ref{th:glm-fixed-prob}, by Chebyshev's inequality we have
\[
\begin{aligned}
\Prob\left[ |\wh V_0| \ge \frac{n\Delta}{c} \right] 
\le \frac{c^2 \Var[n\wh\nu_0(\theta^*)]}{n^2 \Delta^2} 
= \frac{c^2}{n \Delta^2} \Var\left[\frac{(y - \wt y)^2}{2a''(\eta_{x}(\theta^*))} \right]
= \frac{c^2}{n \Delta^2} \left( \E\left[\frac{(y - \wt y)^4}{4a''(\eta_{x}(\theta^*))^2}\right] - \nu_0(\theta^*) \right);
\end{aligned}
\]
here~$(x,y) \sim \Prob_0$ and~$\wt y$ is independent of~$y$ (with the same distribution) conditionally on~$x$, cf.~Assumption~\ref{ass:resampling}.
Furthermore, 
\[
\E\left[\frac{(y - \wt y)^4}{4a''(\eta_{x}(\theta^*))^2}\right] 
= \E\left[\frac{(y - \E[y|x] - (\wt y - \E[\wt y|x]) )^4}{4a''(\eta_{x}(\theta^*))^2}\right].
\]
Thus, for independent~$\xi := y - \E[y|x]$ and~$\wt \xi := \wt y - \E[y|x]$ (so that~$\E[\xi|x] = \E[\wt \xi|x] = 0$) we have
\[
\E\left[\frac{(y - \wt y)^4}{4a''(\eta_{x}(\theta^*))^2}\right] 
= \E\left[\frac{\E[(\xi - \wt\xi)^4|x]}{4a''(\eta_{x}(\theta^*))^2}\right] 
= \E\left[\frac{2\E[\xi^4|x] + 6\E[\xi^2|x]^2}{4a''(\eta_{x}(\theta^*))^2}\right] 
= \frac{\kappa_0(\theta^*)}{2} + \frac{3\nu_0^2(\theta^*)}{2}.
\]
As a result,
\[
\Prob\left[ |\wh V_0| \ge \frac{n\Delta}{c} \right] 
\le \frac{c^2[\kappa_0(\theta^*) + \nu_0^2(\theta^*)]}{2n\Delta^2} 
\le \frac{c^2\kappa_0(\theta^*)}{n\Delta^2}.
\]
Proceeding for~$\wh V_1$ in a similar manner, we arrive at
\[
\Prob\left[ |\wh V_0 - \wh V_1| \ge \frac{2n\Delta}{c} \right] 
\le \frac{c^2[\kappa_0(\theta^*) + \kappa_1(\theta^*)]}{n\Delta^2},
\]
that is,
\[
\Prob_{\cH_1}\left[ |\wh V_0 - \wh V_1| \ge \frac{2n\Delta}{c} \right] 
\le \frac{c^2[\kappa_1(\theta_1) + \kappa_0(\theta_1)]}{n\Delta^2},
\]
As a result, under~\eqref{eq:n-kurtosis-apx} (assuming hypothesis~$\cH_1$ w.l.o.g.) the bound~$|\wh V_0 - \wh V_1| \le 2n\Delta/c$ holds with fixed probability. Thus, we can combine this result with steps~\proofstep{1}-\proofstep{3} in the proof of Theorem~\ref{th:glm-fixed-prob}, (slightly increasing the constants), and arrive at a fixed type II error probability bound for~\eqref{eq:test-glm-adapt}.
\qed
\section{Suboptimality of test~\eqref{eq:test_grad}}
\label{app:grad}

For simplicity, we let~$n_0 = n_1 = n$ and consider the setup of a well-specified linear model with fixed design (cf.~Sec.~\ref{sec:gaussian}); the example can easily be extended for  random design. We use the abridged notation~$(X,Y) = (X^{(0)},Y^{(0)})$ and~$(X',Y') = (X^{(1)},Y^{(1)})$ for the two samples, and we also let~$\xi = Y - X\theta_0$ and~$\xi' = Y' - X'\theta_1$. 
Recall that, by Proposition~\ref{prop:upper-fixed-design}, our proposed test in~\eqref{eq:test_1} has type I error probability
\[
\Prob_{\cH_0}[\wh T_{\Lin} = 1|X, X'] \le C \exp \left(-c n_1 \wh\Delta_1  \min \left\{ 1, \frac{n_1\wh\Delta_1}{\wh r} \right\} \right)
\]
and the complementary bound for type II error (with the replacement~$1 \mapsto 0$ of the subscript).
Here the (empirical) separation measures~$\wh\Delta_0,\wh\Delta_1$ are given by
\[
\wh\Delta_0 = \|\wh\bSigma_0^{1/2}(\theta_1 - \theta_0)\|^2,
\quad
\wh\Delta_1 = \|\wh\bSigma_1^{1/2}(\theta_1 - \theta_0)\|^2,
\]
and~$\wh r = \max \{ \wh r_0, \wh r_1 \}$, where~$\wh r_k = \rank(\wh\bSigma_k)$.
We now construct a problem instance in which the gradient test~$\wh T_{\Grad}$, given by~\eqref{eq:test_grad}, is suboptimal.
To this end, let~$\wh \bSigma_0$ and~$\wh \bSigma_1$ have the same rank~$\wh r$, smallest non-zero eigenvalue~$1/\kappa$ and other eigenvalues equal to~$1$, for~$\kappa > 0$ to be fixed later. 
Now, let us take~$\theta_1 - \theta_0$ as the unit vector in the direction corresponding to the smallest eigenvalue of~$\wh \bSigma_1$. 
Then~$\wh \Delta_1 = {1}/{\kappa} \le \wh \Delta_0.$
On the other hand, the quantities~$\wh \delta_k := \|\wh\bSigma_k(\theta_1 - \theta_0)\|^2$,~$k \in \{0,1\}$, satisfy
\[
\wh \delta_1 = \frac{\wh \Delta_1}{\kappa} = \frac{1}{\kappa^2} \le \wh \delta_0.
\]
Now, observe that~$\wh\delta_1$ controls the sum of errors of the gradient norm test~\eqref{eq:test_grad} in a similar manner as~$\wh\Delta_1$ controls that for test~\eqref{eq:test_1}. Indeed, proceeding as in the proof of Proposition~\ref{prop:upper-fixed-design}, we have
\[
\begin{aligned}
&\Prob_{\cH_0}[\wh T_{\Grad} =1 | X, X'] \\
&= \Prob 
\Big[ 
\underbrace{\tfrac{1}{n}\| X^{\top}\xi\|^2 - \E[\tfrac{1}{n}\| X^{\top}\xi\|^2|X]}_{U}  
+ 
\underbrace{\E[\tfrac{1}{n}\| (X')^\top \xi'\|^2|X'] - \tfrac{1}{n} \| (X')^\top \xi'\|^2}_{V}    
+ 
\underbrace{2\langle \wh \bSigma_1 (\theta_1 - \theta_0), (X')^\top \xi' \rangle}_{W} \\ &\quad\quad\geq n \wh \delta_1 
\Big].
\end{aligned}
\]
As our goal is to show suboptimality, let us bound~$\Prob_{\cH_0}[\wh T_{\Grad}=1 | X,X']$ from below. Since~$\tfrac{1}{n}XX^\top$ has the same eigenvalues as~$\wh \bSigma_0 = \frac{1}{n}X^\top X$, and similarly for~$\frac{1}{n} X' (X')^\top$, we have that~$\tfrac{1}{n}\| X^{\top}\xi\|^2$ is distributed as~$\frac{1}{\kappa}\zeta^2 + \chi_{\wh r-1}^2$, where~$\chi_{\wh r-1}^2$ is chi-squared with~$\wh r-1$ degrees of freedom, and~$\zeta \sim \cN(0,1)$ is independent from~$\chi_{\wh r-1}^2.$ 
Thus, we estimate the (conditional to~$X$) variance of~$U$ above as
\[
\Var[U|X] \le \wh r( \wh r-1) + 3\wh r - (\wh r-1)^2 \le 4\wh r,
\]
where we used that~$\E[\zeta^4] = 3$. Similarly,~$\Var[V|X'] \le 4 \wh r$. Finally, by the definition of~$\theta_1 - \theta_0$ we find 
\[
\Var[W|X'] = 4 \big\| X' \wh\bSigma_1(\theta_1 - \theta_0) \big\|^2 
= \frac{4}{\kappa^2} \|X' (\theta_1 - \theta_0)\|^2
= \frac{4n \wh\Delta_1}{\kappa^2} = \frac{4n \wh\delta_1}{\kappa}. 
\]
Thus,
\[
\Var[\tfrac{1}{\kappa} \zeta^2 - \tfrac{1}{\kappa} + V + W] 
\le \Var[U + V + W] 
\le 3\left(8\wh r + \frac{4n\wh\delta_1}{\kappa} \right). 
\]
By Chebyshev's inequality, any centered random variable~$Z$ satisfies
$
\Prob[Z \le 2\sqrt{\Var[Z]}] \ge 1/2. 
$
Applying this to~$Z =\tfrac{1}{\kappa} \zeta^2 - \tfrac{1}{\kappa} + V + W$, we arrive at
\[
\begin{aligned}
\Prob_{\cH_0}[\wh T_{\Grad}=1 | X, X'] 
&\ge \frac{1}{2} \Prob\left[\chi_{\wh r-1}^2 - (\wh r-1) \ge n\wh\delta_1 - 10\left(\sqrt{\wh r + \frac{n\wh\delta_1}{\kappa}}\right)\right] \\
&\ge \frac{1}{2} \Prob\left[\chi_{\wh r-1}^2 - (\wh r-1) \ge \frac{1}{\kappa} \left( n\wh\Delta_1- 10\sqrt{\kappa^2 \wh r + n\wh\Delta_1}\right)\right] 
\end{aligned}
\]
Now, for~$\wh r$ large enough, consider the regime
\begin{equation}
\label{eq:hard-range}
40\sqrt{\wh r} \le \frac{n\wh\Delta_1}{\kappa} \le 3\wh r  
\quad \iff \quad
40\sqrt{\wh r} \le \frac{n}{\kappa^2} \le 3\wh r.
\end{equation}
By simple algebra, in this regime we have that~$10\sqrt{\kappa^2 \wh r + n\wh\Delta_1} \le n\wh\Delta_1/2$, thus
\[
\Prob_{\cH_0}[\wh T_{\Grad}=1 | X, X']  
\ge \frac{1}{2} \Prob\left[\chi_{\wh r-1}^2 - (\wh r-1) \ge \frac{n\wh\Delta_1}{2\kappa} \right].
\]
On the other hand, in this regime we can apply the lower bound for the right tail probability of the chi-squared distribution (cf.~\eqref{eq:chi-squared-lower-estimates}), which gives
\[
\Prob_{\cH_0}[\wh T_{\Grad}=1 | X, X']  
\ge c\exp\left( -\frac{C}{\kappa^2} \frac{n^2\wh\Delta_1^2}{\wh r}\right).
\]
Now, allowing~$\wh r$ to grow, let us put~$\kappa = \sqrt{\wh r}$ and~$n = 40\wh r^{3/2}$ 
(this is on the ``left edge'' of~\eqref{eq:hard-range}). Then
\[
\Prob_{\cH_0}[\wh T_{\Grad}=1 | X, X']  \ge c.
\]
On the other hand, in this regime we have~$n\wh\Delta_1 \asymp \wh r$ up to a constant, whence, returning to test~\eqref{eq:test_1}, 
\[
\Prob_{\cH_0}[\wh T_{\Lin} = 1|X, X'] + \Prob_{\cH_1}[\wh T_{\Lin} = 0|X, X'] 
\le C \exp \bigg(-\frac{cn^2\wh\Delta_1^2}{\wh r}\bigg) \le C\exp(-c\wh r).
\]
Thus, we see that there is a widening gap -- as~$\wh r$ grows -- between the error bounds for the tests~$\wh T_{\Lin}$ and~$\wh T_{\Grad}$ (and moreover,~$\wh T_{\Grad}$ is not even consistent in this regime, whereas~$\wh T_{\Lin}$ still is).
\qed
\section{Why empirical prediction scores do not violate non-disclosure}
\label{app:fisher}

Recall that, when discussing the {\bf non-disclosure mechanism} behind our proposed testing protocol in Section~\ref{sec:protocol}, we replaced the data~$
\big( Z^{(0)},Z^{(1)}, \ell, \wh S_0(\theta^*), \wh S_1(\theta^*) \big)$ with~$\left( Z^{(0)},Z^{(1)}, \ell \right)$, thus neglecting the additional statistical information about~$\Prob_0, \Prob_1$ contained in the prediction scores. 
We justify this step by comparing the traces of Fisher information matrices corresponding to the samples and prediction scores.
For the sake of concreteness we focus on the well-specified linear model setup as in Section~\ref{sec:gaussian} with $\ell$ being the squared loss.
For simplicity, we assume~$n_0 = n_1 = n$, and the design is non-random and the same in both samples,\footnote{One can also let~$X^{(0)}$ and~$X^{(1)}$ be different: the interested reader may verify that, modulo a few adjustments. The argument we are about to present will essentially remain the same.} so that~$Y^{(k)} \sim \cN(X\theta_k,\Id_n)$ fully specifies the distributions of the samples. 
We also define~$\xi^{(k)} := Y^{(k)} - X \theta_k \sim \cN(0,\Id_n)$.
Given our focus on prediction rather than estimation, we shall pass to the ``prediction parameters''~$\mu_k := X\theta_k \in \R^n$ living in the subset~$\col(X) \subseteq \R^n$ of dimension~$\wh r$. 
Finally, w.l.o.g.~we grant the null hypothesis, so that,~$\mu^* = \mu_0$ is the mean parameter corresponding to~$\theta^* = \theta_0$. 

We shall now focus on the Fisher information about~$\mu^* = \mu_0$ in~$\wh S_0(\theta_0)$ and~$\wh S_1(\theta_0)$. 
Let us first observe that~$\wh S_k(\theta) = \| \bPi_{X} [Y^{(k)} - X\theta] \|^2$ can be reparametrized in terms of~$\mu := X\theta$ as follows:
\[
\wh S_k(\theta) 
= \|\bPi_{X}[\xi^{(k)}]  + X(\theta_k - \theta)\|^2 
= \|\zeta^{(k)}  + \mu_k - \mu\|^2 =:  \wh T_k(\mu),
\]
where~$\zeta^{(k)} := \bPi_X \xi^{(k)} \sim \cN(0,\bPi_X)$. 
Now, consider the corresponding Fisher information matrices:
\[
\mathfrak{I}_{\wh T_k(\mu_0)}(\mu) := \E\left[ \nabla \cL_{\wh T_k(\mu_0)}(\mu) \otimes \nabla \cL_{\wh T_k(\mu_0)}(\mu) \right],
\]
where $v\otimes u = vu^{{\sf T}}$ and ~$\cL_{\wh T_k(\mu_0)}$ is the log-likelihood of~$\wh T_k(\mu_0)$.
The functional
\[
\Tr[\mathfrak{I}_{\wh T_k(\mu_0)}(\mu)] = 
\E\left[ \Big\| \nabla \cL_{\wh T_k(\mu_0)}(\mu)\Big\|^2\right]
\]
evaluated at~$\mu = \mu_0$ quantifies the amount of information about~$\mu^* = \mu_0$ contained in~$\wh T_k (\mu_0)$, and can be compared with~$\Tr[\mathfrak{I}_{\wh Z_k}(\mu_0)]$, where~$\mathfrak{I}_{\wh Z_k}(\mu)$ corresponds to the likelihood of~$\wh Z_k$. 
In addition, by rotational invariance of Gaussian distribution we have
\[
{\wh T_k}(\mu) 
= \|\omega^{(k)} + \|\mu - \mu_k\| \ue_1\|^2,
\] 
where~$\ue_1$ is the first canonical vector in~$\R^{\wh r}$, and~$\omega^{(k)} \sim \cN(0,\Id_{\wh r})$. 
Now, by the monotonicity property of Fisher information matrices (see, e.g.,~\cite[Theorem~2.86]{schervish1995theory}) we have that
\[
\mathfrak{I}_{\wh T_k(\mu_0)}(\mu) \preccurlyeq 
\mathfrak{I}_{\omega^{(k)} + \|\mu_0 - \mu_k\| \ue_1}(\mu). 
\]
Furthermore,~$\omega^{(k)} + \|\mu - \mu_k\| \ue_1$ has independent entries, of which only the first one depends on~$\mu$, so
\[
\mathfrak{I}_{\omega^{(k)} + \|\mu_0 - \mu_k\| \ue_1}(\mu) 
= 
\mathfrak{I}_{\omega_1^{(k)} + \|\mu_0 - \mu_k\|}(\mu).
\]
Thus,~$\Tr[\mathfrak{I}_{\wh T_k(\mu_0)}(\mu)] \le \Tr[\mathfrak{I}_{\omega_1 + \|\mu_0 - \mu_k\|}(\mu)]$, where~$\omega_1 \sim \cN(0,1)$. 
Finally, the latter trace can be computed explicitly. 
Indeed, the log-likelihood of the observation~$\omega_1 + \|\mu_0 - \mu_k\| \sim \cN(\|\mu_0 - \mu_k\|,1)$ at arbitrary~$\mu$ (i.e., with the actual data-generating distribution corresponding to~$\mu = \mu_0$) is given by
\[
\cL_{\omega_1 + \|\mu_0 - \mu_k\|} (\mu) = -\tfrac{1}{2} (\omega_1 + \|\mu_0 - \mu_k\|  - \| \mu - \mu_k\|)^2. 
\]
Taking the gradient (in~$\mu$) through the composition formula, we get
\[
\begin{aligned}
\nabla \cL_{\omega_1 + \|\mu_0 - \mu_k\|} (\mu) 
&= 
( \omega_1 + \|\mu_0 - \mu_k\|  - \| \mu - \mu_k\| ) \, \nabla (\|\mu - \mu_k\|) \\
&= 
\frac{\omega_1 + \|\mu_0 - \mu_k\|  - \| \mu - \mu_k\|}{\|\mu - \mu_k\|}
(\mu - \mu_k) 
= 
(\omega_1 + \|\mu_0 - \mu_k\|  - \| \mu - \mu_k\|) \, \ue_{\mu -\mu_k},
\end{aligned}
\]
where~$\ue_{\mu - \mu_k}$ is the unit-norm vector in the direction of~$\mu - \mu_k$. 
In particular, at~$\mu = \mu_0$ we get~$\nabla \cL_{\omega_1 + \|\mu_0 - \mu_k\|} (\mu_0) = \omega_1 \ue_{\mu_0 - \mu_k}$, thus~$\| \nabla \cL_{\omega_1 + \|\mu_0 - \mu_k\|} (\mu_0) \|^2 = \omega_1^2$ and
$
\Tr[\mathfrak{I}_{\omega_1 + \|\mu_0 - \mu_k\|}(\mu_0)] 
= 1$. 
We conclude that
\[
\Tr[\mathfrak{I}_{\wh T_k(\mu_0)}(\mu_0)] \le 1, \quad k \in \{0,1\}.
\]

Let us now contrast this with the amount of information about~$\mu_0$ that is already contained in the samples themselves. On one hand,~$Y^{(1)}$ does not reveal any information about~$\mu_0$ as it only depends on~$\mu_1$; hence,~$\Tr[\mathfrak{I}_{Y^{(1)}}(\mu_0)] = 0$. On the other hand,~$\Tr[\mathfrak{I}_{Y^{(0)}}(\mu_0)] = n$ by a straightforward calculation in the Gaussian sequence model. 
Thus, we conclude that each prediction scores contains information about~$\mu^*$ roughly corresponding to a single additional data point -- just as we claimed.  
\qed

\end{document}